\documentclass{article}

\usepackage{nips15submit_e,times}
    \usepackage{authblk}
    
\usepackage{hyperref}
\usepackage{url}
\usepackage{wrapfig}
\usepackage{mathrsfs}
\usepackage{epsfig}
\usepackage{graphics}
\usepackage{epsf}
\usepackage{amsmath, amsthm, amssymb, multirow, paralist}
\usepackage{url}
\usepackage{color}
\usepackage{algorithm,algorithmic}
       \usepackage{enumitem}
\usepackage{booktabs}

    \newtheorem{theorem}{Theorem}
    
   \newtheorem{lemma}{Lemma}
    
    \newtheorem{definition}{Definition}
     \newtheorem{assumption}{Assumption}

    \def \R {\mathbb{R}}

    \def \v {\mathbf{v}}
    
    \def \x {\mathbf{x}}
    
    \def \E {\mathrm{E}}
    
    \def \x {\mathbf{x}}
    \def \xt {\widetilde{\x}}
    \def \p {\mathbf{p}}
    \def \a {\mathbf{a}}
    \def \diag {\mbox{diag}}
    
    \def \b {\mathbf{b}}

    \def \z {\mathbf{z}}
    \def \s {\mathbf{s}}
    
    \def \y {\mathbf{y}}

    \def \G {\mathcal{G}}
    \def \g {\mathbf{g}}
    \def \F {\mathcal{F}}
    \def \I {\mathbb{I}}

    \def \xh {\widehat{\x}}

    \def \N {\mathbb{N}}
    \def \sgd {\textsc{Sgd}}
     \def \ada {\textsc{AdaGrad}}
         \def \shb {\textsc{Shb}}
         \def \snag {\textsc{Snag}}

    \nipsfinalcopy 

    \title{Universal  Stagewise Learning for Non-Convex Problems with  Convergence on  Averaged Solutions}

\author[1]{Zaiyi Chen }
\author[2 ]{Zhuoning Yuan}
\author[3]{Jinfeng Yi}
\author[3]{Bowen Zhou}
\author[1]{Enhong Chen}
\author[2 ]{Tianbao Yang}

\affil[1]{School of Computer Science, University of Science and Technology of China}
\affil[2]{Department of Computer Science, University of Iowa}
\affil[3]{JD AI Research, China}
\affil[ ]{\texttt{czy6516@mail.ustc.edu.cn, tianbao-yang@uiowa.edu}}

\begin{document}
\maketitle 

\vspace*{-0.2in}
\begin{center}{March 6, 2019}\end{center} 
\vspace*{0.2in}

\begin{abstract}
Although stochastic gradient descent (\sgd) method and its variants (e.g., stochastic momentum methods, \ada) are the choice of algorithms for solving non-convex problems (especially deep learning), there still remain big gaps between the theory and the practice with many questions unresolved. For example, there is still a lack of theories of convergence for {\sgd} and its variants that use stagewise step size and return an averaged solution in practice. In addition, theoretical insights of why adaptive step size of {\ada} could improve non-adaptive step size of {\sgd} is still missing for non-convex optimization.   This paper aims to address these questions and fill the gap between theory and practice. We propose a universal stagewise optimization framework for a broad family of {\bf non-smooth non-convex} (namely weakly convex) problems with the following key features: (i) at each stage any suitable stochastic convex optimization algorithms (e.g., {\sgd}  or \ada)  that return an averaged solution can be employed for minimizing a regularized convex problem; (ii) the step size is decreased in  a stagewise manner; (iii)  an averaged solution  is returned as the final solution that is selected from all stagewise averaged solutions with sampling probabilities  {\it increasing} as the stage number. Our theoretical results of stagewise {\ada}  exhibit its adaptive convergence, therefore shed insights on its faster convergence for problems with sparse stochastic gradients than stagewise \sgd. To the best of our knowledge, these new results are the first of their kind for addressing the unresolved issues of existing theories  mentioned earlier. Besides theoretical contributions, our empirical studies show that our stagewise {\sgd} and {\ada}  improve the generalization performance of existing variants/implementations of {\sgd} and {\ada}.
\end{abstract}

\section{Introduction}
Non-convex optimization has recently received increasing attention due to its popularity in emerging machine learning tasks, particularly for learning deep neural networks. One of the keys to the success of deep learning for big data problems is the employment of simple stochastic algorithms such as {\sgd} or {\ada}~\cite{DBLP:conf/nips/KrizhevskySH12,Dean:2012:LSD:2999134.2999271}. Analysis of these stochastic algorithms for  non-convex optimization is an important and interesting research topic, which already attracts much attention from the community  of theoreticians~\cite{DBLP:journals/siamjo/GhadimiL13a,DBLP:journals/mp/GhadimiL16,DBLP:journals/mp/GhadimiLZ16,yangnonconvexmo,sgdweakly18,ada18bottou,ada18orabona}. However, one issue that has been largely ignored in existing theoretical results is that the employed algorithms in practice usually differ from their plain versions that are well understood in theory.  Below, we will mention several important heuristics used in practice that have not been well understood for non-convex optimization, which motivates this work. 

First,  a heuristic for setting the step size in training deep neural networks is to change it in a stagewise manner from a large value to a small value (i.e., a constant step size is used in a stage for a number of iterations and is decreased for the next stage)~\cite{DBLP:conf/nips/KrizhevskySH12}, which  lacks theoretical analysis to date.
In existing literature~\cite{DBLP:journals/siamjo/GhadimiL13a,sgdweakly18},  {\sgd} with an iteratively decreasing step size or a small constant step size has been well analyzed for non-convex optimization problems with guaranteed convergence to a stationary point. For example, the existing theory usually suggests an iteratively decreasing step size proportional to $1/\sqrt{t}$ at the $t$-th iteration or a small constant step size, e.g., proportional to $\epsilon^2$ with $\epsilon\ll 1$ for finding an $\epsilon$-stationary solution whose gradient's magnitude (in expectation) is small than $\epsilon$.

Second, the averaging heuristic is usually used in practice, i.e., an averaged solution is returned for prediction~\cite{bottou-2010-large}, which could yield improved stability and generalization~\cite{DBLP:conf/icml/HardtRS16}. However, existing theory for many stochastic non-convex optimization algorithms only provides guarantee on a uniformly sampled solution or a non-uniformly sampled solution with {\it decreasing} probabilities for latest solutions~\cite{DBLP:journals/siamjo/GhadimiL13a,yangnonconvexmo,sgdweakly18}. In particular,  if an iteratively decreasing step size proportional to $1/\sqrt{t}$ at the $t$-th iteration is employed, the convergence guarantee was provided for a random solution that is non-uniformly selected from all iterates with a sampling probability proportional to $1/\sqrt{t}$ for the $t$-th iterate. This means that the latest solution always has the smallest probability to be selected as the final solution, which contradicts to the common wisdom. 
 If a small constant step size is used, then usually a uniformly sampled solution is returned with convergence guarantee. However, both options are seldomly used in practice.

A third common heuristic in practice is to use  adaptive coordinate-wise step size of {\ada}~\cite{Dean:2012:LSD:2999134.2999271}. Although adaptive step size has been well analyzed for convex problems (i.e., when it can yield faster convergence than {\sgd})~\cite{duchi2011adaptive,SadaGrad18},   it still remains an mystery for non-convex optimization with missing insights from theory. Several recent studies have attempted to analyze {\ada} for non-convex problems~\cite{ada18bottou,ada18orabona,adamtypenoncx18,adagradmom18}. Nonetheless, none of them are able to exhibit the adaptive convergence of {\ada} to data as in the convex case and its advantage over {\sgd} for non-convex problems.

To overcome the shortcomings of existing theories for stochastic non-convex optimization, this paper analyzes new algorithms that employ some or all of these commonly used heuristics in a systematic framework, aiming to fill the gap between theory and practice. The main results and contributions are summarized below: 

\begin{itemize}[leftmargin=*]
\vspace*{-0.1in}
\item We propose a universal stagewise optimization framework for solving a  family of non-convex problems, i.e., weakly convex problems, which is broader than smooth non-convex problems and includes some non-smooth non-convex problems. At each stage,  any suitable  stochastic convex optimization algorithms (e.g., {\sgd}, {\ada}) with a constant step size  parameter  can be employed for optimizing a regularized convex problem with a number of iterations, which usually return an averaged solution. The step size parameter is decreased in a stagewise manner following a polynomial decaying scheme. 
\item We analyze several variants of the proposed framework by employing different basic algorithms,  including {\sgd}, stochastic heavy-ball (\textsc{Shb}) method,  stochastic Nesterov's accelerated gradient (\textsc{Snag}) method, stochastic alternating direction methods of multipliers (ADMM), and  {\ada}. We prove the convergence of their stagewise versions for an averaged solution that is randomly selected from all stagewise averaged solutions. 

\item To justify a heuristic approach that returns the last averaged solution in stagewise learning, we present and analyze a  non-uniform sampling strategy over stagewise averaged solutions with sampling probabilities increasing as the stage number. 

\item Regarding the convergence results, for stagewise {\sgd}, \textsc{Shb}, \textsc{Snag}, we establish the same order of iteration complexity for finding a nearly stationary point as the existing theories of their non-stagewise variants. For stagewise {\ada}, we establish an adaptive convergence  for finding a nearly stationary point, which is provably better than (stagewise) {\sgd}, \textsc{Shb}, and \textsc{Snag}  when the cumulative growth of stochastic gradient is slow. 
\item Besides theoretical contributions, we also empirically verify the effectiveness of the proposed stagewise algorithms. In particular, our empirical studies show that (i) the stagewise {\ada} dramatically improves the generalization performance of existing variants of {\ada}, (ii)  stagewise {\sgd, \textsc{Shb}, \textsc{Snag}} also outperform their plain variants with an iteratively decreasing step size; (iii) the proposed stagewise algorithms achieve similar if not better generalization performance than their heuristic variants implemented in existing libraries on standard benchmark datasets. 
\end{itemize}

\section{Related Work}
We review some theoretical results for stochastic non-convex optimization in this section.

{\sgd} for unconstrained smooth non-convex problems was first analyzed  by Ghadimi and Lan~\cite{DBLP:journals/siamjo/GhadimiL13a}, who established an $O(1/\epsilon^4)$ iteration complexity for finding an $\epsilon$-stationary point $\x$ in expectation satisfying $\E[\|\nabla f(\x)\|]\leq \epsilon$, where $f(\cdot)$ denotes the objective function. As mentioned earlier, the returned solution is either a uniformly sampled solution or a non-uniformly sampled one with sampling probabilities proportional to decreasing step size. Similar results were established for the stochastic momentum variants of {\sgd} (i.e.,  {\shb}, {\snag})  by~\cite{yangnonconvexmo,DBLP:journals/mp/GhadimiL16}. Recently,  SGD was also analyzed for (constrained) weakly convex problems,  whose objective function is non-convex and not necessarily smooth, by~Davis and Drusvyatskiy~\cite{sgdweakly18}. Since the objective function could be non-smooth, the convergence guarantee is provided on the magnitude of the Moreau envelope's subgradient with the same order of iteration complexity as in the smooth case. However, none of these studies provide results for algorithms that return an averaged solution. 

Although adaptive variants of {\sgd}, e.g., {\ada}~\cite{duchi2011adaptive}, \textsc{Adam}~\cite{kingma2014adam,sashank2018adam},  were widely used for training deep neural networks, there are few studies on theoretical analysis of these algorithms for non-convex problems. Several recent studies attempted to analyze {\ada} for non-convex problems~\cite{ada18bottou,ada18orabona,adamtypenoncx18,adagradmom18}. Ward et al.~\cite{ada18bottou} only analyzed a variant of {\ada}  that uses a global adaptive step size instead of coordinate-wise adaptive step size as in the original {\ada} used in practice. 
Li and Orabona~\cite{ada18orabona} gave two results about the convergence of variants of {\ada}. One is given in terms of asymptotic convergence for coordinate-wise adaptive step size, and another one is given in terms of non-asymptotic convergence for global adaptive step size. 
When we prepare this manuscript, we note that two recent studies~\cite{adamtypenoncx18,adagradmom18} appeared online, which also analyzed the convergence of {\ada} with coordinate-wise adaptive step size and its momentum variants.   Although all of these studies established an iteration complexity of $O(1/\epsilon^4)$ for different variants of {\ada} for finding an $\epsilon$-stationary solution of a stochastic non-convex optimization problem, none of them can exhibit the potential adaptive advantage of {\ada} over {\sgd} as in the convex case.  To the best of our knowledge, our result is the first one that explicitly shows that  coordinate-wise adaptive step size could yield faster convergence than using non-adaptive step size for non-convex problems similar to that in the convex case. 
 Besides that, these studies also suffer from the following shortcomings: (i) they all assume smoothness of the problem, while we consider non-smooth and non-convex problems;  (ii) their convergence is provided on a solution with minimum magnitude of gradient that is expensive to compute,  though their results also imply a convergence on a random solution selected from all iterates with decreasing sampling probabilities. In contrast, these shortcomings do not exist in this paper. 

Statewise step size has been employed in stochastic algorithms and analyzed for  convex optimization problems~\cite{DBLP:journals/jmlr/HazanK11a,ICMLASSG}. Hazan and  Kale~\cite{DBLP:journals/jmlr/HazanK11a} proposed an epoch-GD method for stochastic strongly convex problems, in which a stagewise step size is used that decreases geometrically and the number of iteration for each stage increases geometrically. Xu et al.~\cite{ICMLASSG} proposed an accelerated stochastic subgradient method for optimizing convex objectives satisfying a local error bound condition, which also employs a stagewise scheme with a constant number of iterations per-stage and geometrically decreasing stagewise step size. The difference from the present work is that they focus on convex problems.

The proposed stagewise algorithm is similar to several existing algorithms in design~\cite{ICMLASSG, davis2017proximally}, which are originated from the proximal point algorithm~\cite{rockafellar76}. I.e., at each stage a proximal strongly convex subproblem is formed and then a stochastic algorithm is employed for optimizing the proximal subproblem inexactly with a number of iterations. Xu et al.~\cite{ICMLASSG} used this idea for solving problems that satisfy a local error bound condition, aiming to achieve faster convergence than vanilla {\sgd}. Davis and  Grimmer~\cite{davis2017proximally} followed this idea to solve weakly convex problems. At each stage,  {\sgd} with decreasing step sizes for a strongly convex problem is employed for solving the proximal subproblem in these two papers. Our stagewise algorithm is developed following the similar idea. The key differences from \cite{ICMLASSG, davis2017proximally} are that (i) we focus on weakly convex problems instead of convex problems considered in \cite{ICMLASSG}; (ii) we use non-uniform sampling probabilities that increase as the stage number to select an averaged solution as the final solution, unlike the uniform sampling used in~\cite{davis2017proximally}; (iii) we present a unified algorithmic framework and convergence analysis, which enable one to employ any suitable stochastic convex optimization algorithms at each stage. It gives us several interesting variants including stagewise stochastic momentum methods, stagewise {\ada}, and stagewise stochastic ADMM. For stagewise {\ada} that employs {\ada} as the basic algorithm for solving the proximal subproblem, we derive an adaptive convergence that is faster than {\sgd} when the cumulative growth of stochastic gradients is slow.  

Finally, we refer readers to several recent papers for other algorithms for weakly convex problems~\cite{modelweakly18,Drusvyatskiy2018}. For example, Drusvyatskiy and  Paquette~\cite{Drusvyatskiy2018} studied a subclass of weakly convex problems whose objective consists of a composition of a convex function and a smooth map, and proposed  a prox-linear method that could enjoy a lower iteration complexity than $O(1/\epsilon^4)$ by smoothing the objective of each subproblem. Davis and  Drusvyatskiy~\cite{modelweakly18} studied a more general algorithm that successively minimizes a proximal regularized stochastic model of the objective function. When the objective function is smooth and has a finite form, variance-reduction based methods are also studied~\cite{Reddi:2016:SVR:3045390.3045425,DBLP:conf/cdc/ReddiSPS16,DBLP:conf/icml/ZhuH16,DBLP:conf/icml/Allen-Zhu17,DBLP:journals/corr/abs/1805.05411}, which have provable faster convergence than {\sgd} in terms of $\epsilon$. However, in all of these studies the convergence is provided on an impractical solution, which is either a solution that gives the minimum value of the (proximal) subgradient's norm~\cite{Drusvyatskiy2018} or on a uniformly sampled solution from all iterations~\cite{Reddi:2016:SVR:3045390.3045425,DBLP:conf/cdc/ReddiSPS16,DBLP:conf/icml/ZhuH16,DBLP:conf/icml/Allen-Zhu17}.


\section{Preliminaries}
The problem of interest in this paper is:
\begin{align}\label{eqn:prob}
\min_{\x\in \Omega} \phi(\x) = \E_\xi [\phi(\x;\xi)]
\end{align}
where $\Omega$ is a closed convex set, $\xi\in\mathcal U$ is a random variable, $\phi(\x)$ and $\phi(\x; \xi)$ are  non-convex functions, with the basic assumptions on the problem given in Assumption~\ref{ass:0}.  

To state the convergence property of an algorithm for solving the above problem. We need to introduce some definitions. These definitions can be also found in related literature, e.g., \cite{davis2017proximally,sgdweakly18}. In the sequel, we let $\|\cdot\|$ denote an Euclidean norm, $[S]=\{1,\ldots, S\}$ denote a set, and $\delta_{\Omega}(\cdot)$ denote the indicator function of the set $\Omega$.
\begin{definition}{(Fr\'echet subgradient)} For a non-smooth and non-convex function $f(\cdot)$, 
\begin{align*}
    \partial_F f(\x) = \bigg\{\v\in\R^d|f(\y)\geq f(x)+\langle \v,\y-\x\rangle + o(\|\y-\x\|),\; \forall \y\in \R^d\bigg\}
\end{align*}
denotes the Fr\'echet subgradient of $f$. 
\end{definition}
\begin{definition}{(First-order stationarity)} For problem~(\ref{eqn:prob}), a point $\x\in\Omega$ is a first-order stationary point if 
\begin{align*}
   0\in \partial_F (\phi + \delta_\Omega)(\x),
\end{align*}
where $\delta_\Omega$ denotes the indicator function of $\Omega$. Moreover, a point $\x$ is said to be $\epsilon$-stationary  if 
\begin{align}
\text{dist}(0, \partial_F (\phi + \delta_\Omega)(\x))\leq \epsilon. 
\end{align}
where $\text{dist}$ denotes the Euclidean distance from a point to a set. 
\end{definition}
\begin{definition}{(Moreau Envelope and Proximal Mapping)} For any function $f$ and $\lambda>0$, the following function is called a Moreau envelope of $f$
\begin{align*}
   f_\lambda(\x)= \min_{\z} f(\z) + \frac{1}{2\lambda}\|\z - \x\|^2.
\end{align*}
Further, the optimal solution to the above problem denoted by 
\[
\text{prox}_{\lambda f}(\x) =\arg\min_{\z} f(\z) + \frac{1}{2\lambda}\|\z - \x\|^2
\]is called a proximal mapping of $f$. 
\end{definition}

\begin{definition}{(Weakly convex)} A function  $f$ is $\rho$-weakly convex, if $f(\x) + \frac{\rho}{2}\|\x\|^2$ is convex. 
\end{definition}

It is known that if $f(\x)$ is $\rho$-weakly convex and $\lambda<\rho^{-1}$, then its Moreau envelope $f_{\lambda}(\x)$ is $C^1$-smooth with the gradient given by (see e.g.~\cite{sgdweakly18})
\begin{align*}
\nabla f_{\lambda}(\x) = \lambda^{-1}( \x  - \text{prox}_{\lambda f}(\x))
\end{align*}
A small norm of $\nabla f_{\lambda}(\x)$ has an interpretation that $\x$ is close to a point that is nearly stationary. In particular for any $\x\in\R^d$, let $\xh = \text{prox}_{\lambda f}(\x)$, then we have
\begin{align}\label{eqn:keye}
f(\xh) & \leq f(\x),\quad \|\x - \xh\| = \lambda \|\nabla f_{\lambda}(\x)\|,\quad\text{dist}(0, \partial f(\xh))\leq \|\nabla f_\lambda(\x)\|.
\end{align}
This means that  a point $\x$ satisfying $\|\nabla f_\lambda(\x)\|\leq \epsilon$  is close to a point in distance of $O(\epsilon)$ that is $\epsilon$-stationary. 

It is notable that for a non-smooth non-convex function $f(\cdot)$, there could exist  a sequence of solutions $\{\x_k\}$ such that $\nabla f_\lambda(\x_k)$ converges while $\text{dist}(0, \partial f(\x_k))$ may not converge~\cite{Drusvyatskiy2018}. To handle such a challenging issue  for non-smooth non-convex problems, we will follow existing works~\cite{modelweakly18,Drusvyatskiy2018,davis2017proximally} to prove the near stationarity in terms  of $\nabla f_\lambda(\x)$.  In the case when $f$ is smooth, $\|\nabla f_{\lambda}(\x)\|$ is closely related to the magnitude of the projected gradient $\mathcal G_\lambda(\x)$ defined below, which has been used as a criterion for constrained non-convex optimization~\cite{DBLP:conf/cdc/ReddiSPS16},
\begin{align}
\mathcal G_\lambda(\x) = \frac{1}{\lambda}(\x - \text{prox}_{\lambda \delta_\Omega}(\x - \lambda \nabla f(\x))).
\end{align}
It was shown that when $f(\cdot)$ is smooth with $L$-Lipschitz continuous gradient~\cite{Drusvyatskiy16a}:
\begin{align}
(1 - L \lambda)\|\mathcal G_\lambda(\x)\|\leq \|\nabla f_{\lambda}(\x)\|\leq (1+L\lambda) \|\mathcal G_\lambda(\x)\|, \forall \x\in\Omega.
\end{align}
Thus, the  near stationarity in terms  of $\nabla f_\lambda(\x)$ implies the  near stationarity in terms  of $\mathcal G_\lambda(\x)$ for a smooth function $f(\cdot)$.

Now, we are ready to state the basic assumptions of the considered problem~(\ref{eqn:prob}).
\begin{assumption}\label{ass:0}\hspace{\fill}
\begin{enumerate}[label=(A\arabic*)]
    \item\label{ass:1}{
        There is a measurable mapping $g: \Omega\times \mathcal U\rightarrow \R$ such that $\E_\xi [g(\x; \xi)] \in \partial_F \phi(\x)$ for any $\x\in\Omega$.
    }
    \item\label{ass:2}{
 For  any $\x\in\Omega$, $\E[\|g(\x; \xi)\|^2]\leq G^2$.   }
    \item\label{ass:3}{
        Objective function $\phi$ is $\mu$-weakly convex.
    }
    \item\label{ass:4}{
         there exists $\Delta_\phi>0$ such that $\phi(\x) - \min_{\z\in\Omega}\phi(\z)\leq \Delta_\phi$ for any $\x\in\Omega$.
    }
    
\end{enumerate}
\end{assumption}
{\bf Remark:} Assumption~1-\ref{ass:1},~1-\ref{ass:2} assume  a stochastic subgradient is available for the objective function and its Euclidean norm square is bounded in expectation, which are standard assumptions for non-smooth optimization. Assumption~\ref{ass:3} assumes weak convexity of the objective function, which is weaker than assuming smoothness. Assumption~\ref{ass:4} assumes that the objective value with respect to the optimal value is bounded. Below, we present some examples of objective functions in machine learning that are weakly convex. 

\paragraph{Ex. 1: Smooth Non-Convex Functions.} If $\phi(\cdot)$ is a $L$-smooth function (i.e., its gradient is $L$-Lipschitz continuous), then it is $L$-weakly convex. 

\paragraph{Ex. 2: Additive Composition.} Consider 
\begin{align}
\phi(\x)  = \E[f(\x; \xi)] + g(\x),
\end{align}
where $\E[f(\x; \xi)]$ is a $L$-weakly convex function, and $g(\x)$ is a closed convex function. In this case, $\phi(\x)$ is $L$-weakly convex. This class includes many interesting regularized problems in machine learning with smooth losses and convex regularizers. For smooth non-convex loss functions, one can consider truncated square loss for robust learning, i.e., $f(\x; \a, b) = \phi_\alpha((\x^{\top}a - b)^2)$, where $\a\in\R^d$ denotes a random data and $\b\in\R$ denotes its corresponding output, and $\phi_\alpha$ is a smooth non-convex truncation function (e.g., $\phi_\alpha(x) = \alpha \log(1 + x/\alpha), \alpha>0$). Such truncated non-convex losses have been considered in~\cite{loh2017statistical}. When $|x^2\phi''_\alpha|\leq k$ and $\|\a\|\leq R$, it was proved that $f(\x; \a, b)$ is a smooth function with Lipschitz continuous gradient~\cite{DBLP:journals/corr/abs-1805-07880}. For $g(\x)$, one can consider any existing convex regularizers, e.g., $\ell_1$ norm, group-lasso regularizer~\cite{Yuan06modelselection}, graph-lasso regularizer~\cite{kim2009multivariate}.

\paragraph{Ex. 3: Convex and Smooth Composition} Consider 
\[
\phi(\x; \xi)  = h(c(\x; \xi))
\]
where $h(\cdot): \R^m\rightarrow \R$ is  closed convex and $M$-Lipschitz continuous, and $c(\x; \xi): \R^d\rightarrow \R^m$ is nonlinear smooth mapping with $L$-Lipschitz continuous gradient. This class of functions has been considered in~\cite{Drusvyatskiy2018} and it was proved that $\phi(\x; \xi)$ is $ML$-weakly convex. An interesting example is phase retrieval~\cite{}, where $\phi(\x; a, b) = |(\x^{\top}\a)^2 - b|$. More examples of this class can be found in~\cite{modelweakly18}. 

\paragraph{Ex. 4: Smooth and Convex Composition} Consider 
\[
\phi(\x; \xi)  = h(c(\x; \xi))
\]
where $h(\cdot): \R\rightarrow \R$ is  a $L$-smooth function satisfying $h'(\cdot)\geq 0$, and $c(\x; \xi): \R^d\rightarrow \R$ is convex and $M$-Lipschitz continuous.  This class of functions has been considered in~\cite{DBLP:journals/corr/abs-1805-07880} for robust learning and it was proved that $\phi(\x; \xi)$ is $ML$-weakly convex. An interesting example is truncated Lipschitz continuous loss $\phi(\x; \a, b) = \phi_\alpha(\ell(\x^{\top}\a, b))$, where $\phi_\alpha$ is a smooth truncation function with $\phi'(\cdot)\geq 0$ (e.g., $\phi_\alpha=\alpha\log(1+x/\alpha)$) and $\ell(\x^{\top}\a, b)$ is a convex and Lipschitz-continuous function (e.g., $|\x^{\top}\a - \b|$ with bounded $\|\a\|$).

\paragraph{Ex. 5: Weakly Convex Sparsity-Promoting  Regularizers} Consider 
\[
\phi(\x; \xi)  =  f(\x; \xi) + g(\x),
\]
where $f(\x; \xi)$ is a convex or a weakly-convex function, and $g(\x)$ is a weakly-convex sparsity-promoting regularizer. Examples of weakly-convex sparsity-promoting regularizers include: 
\begin{itemize}
\item Smoothly Clipped Absolute Deviation (SCAD) penalty~\cite{doi:10.1198/016214501753382273}: $g(\x) = \sum_{i=1}^d g_\lambda(x_i)$ and 
\begin{align*}
g_\lambda(x) =\left\{\begin{array}{lc}\lambda |x|& |x|\leq \lambda \\  - \frac{x^2 - 2a\lambda |x| + \lambda^2}{2(a-1)} & \lambda<|x|\leq a\lambda\\ \frac{(a+1)\lambda^2}{2}& |x|>a\lambda\end{array} \right.
\end{align*}
where $a>2$ is fixed and $\lambda>0$. It can be shown that SCAD penalty is $(1/(a-1))$-weakly convex~\cite{loh2017statistical}.  
\item Minimax Convex Penalty (MCP)~\cite{citeulike:12268172}: $g(\x) = \sum_{i=1}^d g_\lambda(x_i)$ and
\begin{align*}
g_\lambda(x) = \text{sign}(x)\lambda\int_{0}^{|x|}\left(1 - \frac{z}{\lambda b}\right)_+ dz 
\end{align*} 
where $b>0$ is fixed and $\lambda>0$. MCP is $1/b$-weakly convex~\cite{loh2017statistical}. 
\end{itemize}

%

\section{Stagewise Optimization: Algorithms and Analysis}
In this section, we will present the proposed algorithms and the analysis of their convergence. We will first  present a Meta algorithmic framework  highlighting the key features of the proposed algorithms and then present several variants of the Meta algorithm by employing different basic algorithms.

The Meta algorithmic framework is described in Algorithm~\ref{alg:meta}.
There are several key features that differentiate Algorithm~\ref{alg:meta} from existing stochastic algorithms that come with theoretical guarantee. First, the algorithm is run with multiple stages.  At each stage, a stochastic algorithm ($\mathcal{\text{SA}}$) is called to optimize a proximal problem $f_s(\x)$ inexactly that consists of the original objective function and a quadratic term, which is guaranteed to be convex due to the weak convexity of $\phi$ and $\gamma< \mu^{-1}$. The convexity of $f_s$ allows one to employ any suitable  existing stochastic algorithms (cf. Theorem~\ref{thm:meta}) that have convergence guarantee for convex problems. It is notable that \textsc{SA} usually returns an averaged solution $\x_s$ at each stage.  
Second, a decreasing sequence of step size parameters $\eta_s$ is used.  At each stage, the \textsc{SA} uses a constant step size parameter $\eta_s$ and runs the updates for a number of $T_s$ iterations. We do not initialize $T_s$ as it might be adaptive to the data as in stagewise {\ada}. Third, the final solution is selected from the stagewise averaged solutions $\{\x_s\}$ with non-uniform sampling probabilities proportional to a sequence of non-decreasing positive weights $\{w_s\}$. In the sequel, we are particularly interested in $w_s = s^\alpha$ with $\alpha>0$. The setup of $\eta_s$ and $T_s$ will depend on the specific choice of \textsc{SA}, which will be exhibited later for different variants.

\begin{algorithm}[t]
    \caption{A Meta Stagewise Algorithm: Stagewise-SA}\label{alg:meta}
    \begin{algorithmic}[1]
        \STATE \textbf{Initialize:} a sequence of decreasing step size parameters $\{\eta_s\}$,  a sequence of non-decreasing positive weights $\{w_s\}$,  $\x_0\in\Omega$, $\gamma<\mu^{-1}$
        \FOR {$s = 1,\ldots, S$}
        \STATE Let $f_{s}(\cdot) = \phi(\cdot)+\frac{1}{2\gamma}\| \cdot-\x_{s-1}\|^2$
        \STATE $\x_{s} = \textsc{SA}(f_{s}, \x_{s-1}, \eta_s, T_s)$ // $\x_s$ is usually an averaged solution
        \ENDFOR
        \STATE \textbf{Return:} $\x_{\tau}$, $\tau$ is randomly chosen from $\{0, \ldots, S\}$ according to probabilities $p_\tau = \frac{w_{\tau+1}}{\sum_{k=0}^{S}w_{k+1}}, \tau=0, \ldots, S$.
    \end{algorithmic}
\end{algorithm}

To illustrate that Algorithm~\ref{alg:meta} is a universal framework such that any suitable \textsc{SA} algorithm can be employed, we present the following  result by assuming that  SA has an appropriate convergence for a convex problem. 
\begin{theorem} \label{thm:meta}
   Let $f(\cdot)$ be a convex function, $\x_* =\arg\min_{\x\in\Omega}f(\x)$ and $\Theta$ denote some problem dependent parameters.  Suppose for $\x_+ = \textsc{SA}(f, \x_0, \eta, T)$, we have
        \begin{align}\label{eqn:sa_converge}
        \E[f(\x_+) - f(\x_*)] \leq \varepsilon_1(\eta, T, \Theta)\|\x_0 - \x_*\|_2^2 + \varepsilon_2(\eta, T, \Theta)(f(\x_0) - f(\x_*)) + \varepsilon_3(\eta, T, \Theta).
        \end{align}
 Under assumption~1-\ref{ass:1},~\ref{ass:3} and~\ref{ass:4}, by running   Algorithm~\ref{alg:meta} with $\gamma=1/(2\mu)$, $w_s=s^\alpha, \alpha> 0$, and with $\eta_s, T_s$ satisfying $\varepsilon_1(\eta_s, T_s, \Theta)\leq 1/(48\gamma), \varepsilon_2(\eta_s, T_s, \Theta)\leq 1/2$, we have
    \begin{align*}
        \E\big[\|\nabla\phi_\gamma(\x_{\tau})\|^2\big] \leq \frac{32\Delta_\phi(\alpha+1)}{\gamma(S+1)}+ \frac{48 \sum_{s=1}^Sw_s \varepsilon_3(\eta_s, T_s,  \Theta)}{\gamma\sum_{s=1}w_s},
    \end{align*}
    where $\tau$ is randomly selected from $\{0, \ldots, S\}$ with probabilities $p_\tau\propto w_{\tau+1}, \tau=0, \ldots, S$. If $\varepsilon_3(\eta_s, T_s, \Theta)\leq c_3/s$ for some constant $c_3\geq 0$ that may depend on $\Theta$, we have
    \begin{align}\label{eqn:metabound}
        \E\big[\|\nabla\phi_\gamma(\x_{\tau})\|^2\big] \leq \frac{32\Delta_\phi(\alpha+1)}{\gamma(S+1)} + \frac{48c_3 (\alpha+1)}{\gamma(S+1)\alpha^{\mathbb I(\alpha<1)}}. 
    \end{align}
    
\end{theorem}

\vspace*{-0.1in}{\bf Remark:} It is notable that the convergence guarantee is provided on a stagewise average solution $\x_\tau$. To justify a heuristic approach that returns the final average solution for prediction, we analyze a new sampling strategy that samples a solution among all stagewise average solutions with sampling probabilities increasing as the stage number increases. This sampling strategy is better than uniform sampling strategy or a strategy  with decreasing sampling probabilities in the existing literature.   
The convergence upper bound in~(\ref{eqn:sa_converge}) of SA covers the results of a broad family of stochastic convex optimization algorithms. 
When $\varepsilon_2(\eta_s, T_s, \Theta)=0$ (as in {\sgd}), the upper bound can be improved by a constant factor. Moreover, we do not optimize the value of $\gamma$. Indeed, any $\gamma<1/\mu$ will work, which only has an effect on constant factor in the convergence upper bound. 

Next, we present several variants of the Meta algorithm by employing {\sgd}, {\ada}, and stochastic momentum methods as the basic SA algorithm, to which we refer as stagewise {\sgd},  stagewise {\ada}, and stagewise stochastic momentum methods, respectively. It is worth mentioning that one can follow similar analysis to analyze other stagewise algorithms by using  their basic convergence for stochastic convex optimization, including RMSProp~\cite{pmlr-v70-mukkamala17a}, AMSGrad~\cite{sashank2018adam}, which is omitted in this paper. 

\begin{proof} Below, we use $\E_s$ to denote expectation over randomness in the $s$-th stage given all history before $s$-th stage. 
 Define
 \begin{align}\label{eqn:z}
 \z_{s} = \arg\min_{\x\in\Omega} f_s(\x) = \text{prox}_{\gamma (\phi+\delta_\Omega)}(\x_{s-1})
 \end{align}
Then $\nabla\phi_\gamma(\x_{s-1})  =\gamma^{-1}(\x_{s-1} - \z_s)$. Then we have $\phi(\x_{s})\geq \phi(\z_{s+1})+\frac{1}{2\gamma}\|\x_{s}-\z_{s+1}\|^2$. Next, we apply Lemma~\ref{lem:sgd} to each call of {\sgd} in stagewise {\sgd}, 
\begin{align*}
\E[f_s(\x_s) - f_s(\z_s)]\leq  \underbrace{\varepsilon_1(\eta_s, T_s, \Theta)\|\x_{s-1} - \z_s\|_2^2 + \varepsilon_2(\eta_s, T_s, \Theta)(f_s(\x_{s-1}) - f_s(\z_s)) + \varepsilon_3(\eta_s, T_s, \Theta)}\limits_{\mathcal E_s}.
\end{align*}
Then 
\begin{align*}
    \E_s\bigg[\phi(\x_{s}) + \frac{1}{2\gamma}\|\x_{s}-\x_{s-1}\|^2 \bigg] &\leq f_{s}(\z_{s}) + \mathcal{E}_{s}\leq f_{s}(\x_{s-1}) + \mathcal{E}_{s}\\
    &\leq \phi(\x_{s-1})+\mathcal{E}_{s}
\end{align*}
On the other hand, we have that 
\begin{align*}
    \|\x_{s}-\x_{s-1}\|^2 =& \|\x_{s}-\z_{s}+\z_{s}-\x_{s-1}\|^2\\
    =& \|\x_{s}-\z_{s}\|^2+\|\z_{s}-\x_{s-1}\|^2 + 2\langle \x_{s}-\z_{s}, \z_{s}- \x_{s-1}\rangle\\
    \geq& (1-\alpha_{s}^{-1})\|\x_{s}-\z_{s}\|^2 + (1-\alpha_{s})\|\x_{s-1}-\z_{s}\|^2
\end{align*}
where the inequality follows from the Young's inequality with $0<\alpha_{s}<1$. Thus we have that
\begin{align}\label{eqn:bound0}
   &\E_s\bigg[ \frac{(1-\alpha_{s})}{2\gamma} \|\x_{s-1}-\z_{s}\|^2\bigg] \leq \E_s\bigg[ \phi(\x_{s-1}) -\phi(\x_{s})+\frac{(\alpha_{s}^{-1}-1)}{2\gamma}\|\x_{s}-\z_{s}\|^2+\mathcal{E}_{s} \bigg]\nonumber\\
    &\leq \E_s\bigg[\phi(\x_{s-1}) -\phi(\x_{s})+\frac{(\alpha_{s}^{-1}-1)}{\gamma(\gamma^{-1}-\mu)}(f_s(\x_s) - f_s(\z_s))+ \mathcal{E}_{s} \bigg]\nonumber\\
    &\leq \E_s\bigg[ \phi(\x_{s-1}) -\phi(\x_{s}) + \frac{\alpha_{s}^{-1}-\gamma\mu}{(1-\gamma\mu)}\mathcal{E}_{s}\bigg]\leq \E_s\bigg[\phi(\x_{s-1}) -\phi(\x_{s})\bigg]\nonumber\\
    &+ \E_s\bigg[\frac{\alpha_{s}^{-1}-\gamma\mu}{(1-\gamma\mu)}\{\varepsilon_1(\eta_s, T_s, \Theta)\|\x_{s-1} - \z_s\|^2 +  \varepsilon_2(\eta_s, T_s, \Theta)(f_s(\x_{s-1}) - f_s(\z_s)) + \varepsilon_3(\eta_s, T_s, \Theta)\}\bigg]
\end{align}
{Next, we bound $f_s(\x_{s-1}) - f_s(\z_s)$ given that $\x_{s-1}$ is fixed. According to the definition of $f_s(\cdot)$, we have
\begin{align*}
&f_s(\x_{s-1}) - f_s(\z_s)  = \phi(\x_{s-1}) - \phi(\z_s) - \frac{1}{2\gamma}\|\z_s - \x_{s-1}\|^2\\
&= \phi(\x_{s-1}) - \phi(\x_s) + \phi(\x_s) - \phi(\z_s) - \frac{1}{2\gamma}\|\z_s - \x_{s-1}\|^2\\
& = [\phi(\x_{s-1}) - \phi(\x_s) ] + \left[f_s(\x_s) - f_s(\z_s) + \frac{1}{2\gamma}\|\z_s - \x_{s-1}\|^2 - \frac{1}{2\gamma}\|\x_{s}- \x_{s-1}\|^2\right]  - \frac{1}{2\gamma}\|\z_s - \x_{s-1}\|^2\\
& \leq  [\phi(\x_{s-1}) - \phi(\x_s) ] + [f_s(\x_s) - f_s(\z_s)].
\end{align*}
Taking expectation over randomness in the $s$-th stage on both sides, we have 
\begin{align*}
&f_s(\x_{s-1}) - f_s(\z_s)\leq \E_s[\phi(\x_{s-1}) - \phi(\x_s) ] + \E_s[f_s(\x_s) - f_s(\z_s)]\\
& \leq \E[\phi(\x_{s-1}) - \phi(\x_s) ] +  \varepsilon_1(\eta_s, T_s, \Theta)\|\x_{s-1} - \z_s\|_2^2 + \varepsilon_2(\eta_s, T_s, \Theta)(f_s(\x_{s-1}) - f_s(\z_s)) + \varepsilon_3(\eta_s, T_s, \Theta).
\end{align*}
Thus, 
\begin{align*}
(1 - \varepsilon_2(\eta_s, T_s, \Theta))(f_s(\x_{s-1}) - f_s(\z_s))\leq  \E[\phi(\x_{s-1}) - \phi(\x_s) ]  +  \varepsilon_1(\eta_s, T_s, \Theta)\|\x_{s-1} - \z_s\|_2^2 +  \varepsilon_3(\eta_s, T_s, \Theta).
\end{align*}
Assuming that $\varepsilon_2(\eta_s, T_s, \Theta)\leq 1/2$, we have
\begin{align*}
\varepsilon_2(\eta_s, T_s, \Theta)(f_s(\x_{s-1}) - f_s(\z_s))\leq  \E_s[\phi(\x_{s-1}) - \phi(\x_s) ]  +  \varepsilon_1(\eta_s, T_s, \Theta)\|\x_{s-1} - \z_s\|_2^2 +  \varepsilon_3(\eta_s, T_s, \Theta).
\end{align*}
Plugging this upper bound into~(\ref{eqn:bound0}), we have
\begin{align}
   &\E_s\bigg[ \frac{(1-\alpha_{s})}{2\gamma} \|\x_{s-1}-\z_{s}\|^2\bigg] \leq \E_s\bigg[\phi(\x_{s-1}) -\phi(\x_{s})\bigg]\nonumber\\
    &+ \E_s\bigg[\frac{\alpha_{s}^{-1}-\gamma\mu}{(1-\gamma\mu)}\{2\varepsilon_1(\eta_s, T_s, \Theta)\|\x_{s-1} - \z_s\|^2 + \phi(\x_{s-1}) - \phi(\x_s) + 2\varepsilon_3(\eta_s, T_s, \Theta)\}\bigg]
\end{align}
By setting $\alpha_s = 1/2, \gamma = 1/(2\mu)$ and assuming $\varepsilon_1(\eta_s, T_s, \Theta)\leq 1/(48\gamma)$, we have
\begin{align*}
   \E_s\bigg[ \frac{1}{8\gamma} \|\x_{s-1}-\z_{s}\|^2 \bigg] &\leq4 \E_s\bigg[\phi(\x_{s-1}) -\phi(\x_{s})\bigg]+6\varepsilon_3(\eta_s, T_s, \Theta)\}
\end{align*}
}
Multiplying both sides by $w_s$, we have that
\begin{align*}
  & w_s\gamma\E_s[\|\nabla \phi_{\gamma}(\x_{s-1})\|^2]\leq \E_s\bigg[32w_s\Delta_{s}  + 48\varepsilon_3(\eta_s, T_s, \Theta)w_s\bigg]
    \end{align*}
 By summing over $s=1,\ldots,S+1$, we have
 \begin{align*}
 \sum_{s=1}^{S+1}w_s \E[\|\nabla \phi_{\gamma}(\x_{s-1})\|^2]\leq \E\bigg[\frac{32}{\gamma}\sum_{s=1}^{S+1}w_s\Delta_{s} + \frac{48}{\gamma}\sum_{s=1}^{S+1}w_s\varepsilon_3(\eta_s, T_s, \Theta) \bigg]
   \end{align*}
Taking the expectation w.r.t. $\tau\in \{0,\ldots, S\}$, we have that
\begin{align*}
   \E[\|\nabla\phi_{\gamma}(\x_{\tau})\|^2]] \leq \E\bigg[\frac{32\sum_{s=1}^{S+1}w_s\Delta_s}{\gamma\sum_{s=1}^{S+1}w_s} + \frac{48\sum_{s=1}^{S+1}w_s\varepsilon_3(\eta_s, T_s, \Theta))}{\gamma\sum_{s=1}^{S+1}w_s}\bigg]
\end{align*}
For the first term on the R.H.S, we have that
\begin{align*}
\sum_{s=1}^{S+1} w_s\Delta_s &= \sum_{s=1}^{S+1}w_s (\phi(\x_{s-1}) - \phi(\x_s)) = \sum_{s=1}^{S+1} (w_{s-1}\phi(\x_{s-1}) - w_s\phi(\x_s)) + \sum_{s=1}^{S+1}(w_s - w_{s-1})\phi(\x_{s-1})\\
&= w_0 \phi(\x_0) - w_{S+1}\phi(\x_{S+1}) +\sum_{s=1}^{S+1}(w_s - w_{s-1})\phi(\x_{s-1})\\
& =\sum_{s=1}^{S+1}(w_s - w_{s-1})(\phi(\x_{s-1}) - \phi(\x_{S+1}))\leq \Delta_\phi\sum_{s=1}^{S+1}(w_s - w_{s-1}) = \Delta_\phi w_{S+1}
\end{align*}
Then, 
\begin{align*}
   \E[\|\nabla\phi_{\gamma}(\x_\tau)\|^2] \leq \frac{32\Delta_\phi w_{S+1}}{\gamma\sum_{s=1}^{S+1}w_s} +\frac{48\sum_{s=1}^{S+1}w_s \varepsilon_3(\eta_s, T_s, \Theta)}{\gamma\sum_{s=1}^{S+1}w_s}
\end{align*}
The standard calculus tells that
\begin{align*}
\sum_{s=1}^Ss^{\alpha}&\geq \int_0^{S}x^\alpha d x= \frac{1}{\alpha+1}S^{\alpha+1}\\
\sum_{s=1}^Ss^{\alpha-1}&\leq SS^{\alpha-1} = S^\alpha, \forall\alpha\geq 1, \quad \sum_{s=1}^Ss^{\alpha-1}\leq \int_{0}^S x^{\alpha-1}d x = \frac{ S^\alpha}{\alpha}, \forall 0<\alpha<1
\end{align*}
Combining these facts and the assumption $\varepsilon_3(\eta_s, T_s, \Theta)\leq c/s$, we have that
\begin{align*}
 \E[\|\nabla\phi_{\gamma}(\x_{\tau})\|^2]\leq \left\{\begin{array}{cc}\frac{32\Delta_\phi(\alpha+1)}{\gamma(S+1)}   +  \frac{48 c (\alpha+1)}{\gamma(S+1)} & \alpha\geq 1\\\\ \frac{32\Delta_\phi(\alpha+1)}{\gamma(S+1)}  + \frac{48 c(\alpha+1)}{\gamma(S+1)\alpha} & 0<\alpha< 1
 \end{array}\right.
\end{align*}
In order to have $\E[\|\nabla\phi_{\gamma}(\x_{\tau})\|^2]\leq \epsilon^2$, we can set $S=O(1/\epsilon^2)$. The total number of iterations is
\begin{align*}
\sum_{s=1}^{S} T_s \leq \sum_{s=1}^{S}12\gamma s\leq 6\gamma S(S+1) = O(1/\epsilon^4)
\end{align*}

\end{proof}

Next, we present several variants of the Meta algorithm by employing {\sgd}, stochastic momentum methods, and {\ada} as the basic SA algorithm, to which we refer as stagewise {\sgd}, stagewise stochastic momentum methods, and stagewise {\ada}, respectively.

\subsection{Stagewise {\sgd}}
\begin{algorithm}[t]
    \caption{{\sgd}$(f, \x_1, \eta, T)$}\label{alg:sgd}
\begin{algorithmic}
    \FOR{$t=1,\ldots, T$}
    \STATE Compute a stochastic subgradient $\g_t$ for $f(\x_{t})$.
    \STATE $\x_{t+1} =\Pi_{\Omega}[\x_t-\eta \g_t$]
    \ENDFOR
    \STATE \textbf{Output}: $\widehat{\x}_T = \sum_{t=1}^T \x_t/T$
\end{algorithmic}
\end{algorithm}
In this subsection, we analyze the convergence of stagewise {\sgd}, in which {\sgd} shown in Algorithm~\ref{alg:sgd} is employed in the Meta framework. Besides Assumption~\ref{ass:0}, we impose the following assumption in this subsection. 

\begin{assumption}\label{ass:0-sgd}\hspace{\fill}
       the domain $\Omega$ is bounded, i.e., there exists $D>0$ such that $\|\x - \y\|\leq D$ for any $\x, \y\in\Omega$. 
\end{assumption}
It is worth mentioning that bounded domain assumption is imposed for simplicity, which  is usually assumed in convex optimization. For machine learning problems, one usually imposes some bounded norm constraint to achieve a regularization. Recently, several studies have found that imposing a norm constraint is more effective than an additive norm regularization term in the objective for deep learning~\cite{HenryRNELC18,DBLP:journals/corr/CDLCDACV}. Nevertheless, the bounded domain assumption is not essential for the proposed algorithm. We present a more involved analysis in the next subsection for unbounded domain $\Omega=\R^d$. The following is a basic convergence result of {\sgd}, whose proof can be found in the literature and is omitted. 

\begin{lemma}\label{lem:sgd}
For Algorithm~\ref{alg:sgd}, assume that $f(\cdot)$ is convex and $\E\|\g_t\|^2\leq G^2, t\in[T]$, then for any $\x\in\Omega$ we have
\begin{align*}
\E[f(\xh_T) - f(\x)]\leq \frac{\|\x - \x_1\|^2}{2\eta T} + \frac{\eta G^2}{2}
\end{align*}
\end{lemma}
To state the convergence, we introduce a notation 
\begin{align}
\nabla\phi_\gamma(\x) = \gamma^{-1}(\x - \text{prox}_{\gamma(\phi+\delta_\Omega)}(\x)),
\end{align}
 which is the gradient of the Moreau envelope of the objective function $\phi+ \delta_{\Omega}$. The following theorem exhibits the convergence of stagewise {\sgd} 
\begin{theorem} \label{thm:sgd}
    Suppose Assumption~\ref{ass:0} and~\ref{ass:0-sgd} hold. By setting $\gamma=1/(2\mu), w_s=s^\alpha, \alpha> 0$, $\eta_{s} = c/s, T_s = 12\gamma s/c$ where $c>0$ is a free parameter, then stagewise {\sgd} (Algorithm~\ref{alg:meta} employing {\sgd}) returns a solution $\x_\tau$ satisfying 
    \begin{align*}
        \E\big[\|\nabla\phi_\gamma(\x_{\tau})\|^2\big] \leq \frac{16\mu\Delta_\phi(\alpha+1)}{S+1} +  \frac{24\mu c \hat G^2(\alpha+1)}{(S+1)\alpha^{\mathbb I(\alpha<1)}},
    \end{align*}
    where $\hat G^2 = 2G^2 + 2\gamma^{-2}D^2$, and $\tau$ is randomly selected from $\{0, \ldots, S\}$ with probabilities $p_\tau\propto w_{\tau+1}, \tau=0, \ldots, S$.
\end{theorem}
{\bf Remark:} To find a solution with $\E\big[\|\nabla\phi_\gamma(\x_\tau)\|^2\big]\leq \epsilon^2$, we can set $S=O(1/\epsilon^2)$ and the total iteration complexity is in the order of $O(1/\epsilon^4)$. The above theorem is essentially a corollary of Theorem~\ref{thm:meta} by applying~\ref{lem:sgd} to $f_s(\cdot)$ at each stage. We present a complete proof in the appendix. 


\subsection{Stagewise stochastic momentum (SM) methods}
In this subsection, we present stagewise stochastic momentum methods and their analysis. In the literature, there are two popular variants of stochastic momentum methods, namely, stochastic heavy-ball method ({\shb}) and stochastic Nesterov's accelerated gradient method ({\snag}). Both methods have been used for training deep neural networks~\cite{DBLP:conf/nips/KrizhevskySH12,DBLP:conf/icml/SutskeverMDH13}, and have been analyzed by~\cite{yangnonconvexmo} for non-convex optimization. To contrast with the results in~\cite{yangnonconvexmo}, we will consider the same unified stochastic momentum methods that subsume {\shb}, {\snag} and {\sgd} as special cases when $\Omega=\R^d$. The updates are presented in Algorithm~\ref{alg:sum}. 

To present the analysis of stagewise SM methods, we first provide a convergence result for minimizing $f_s(\x)$ at each stage. 



\begin{algorithm}[t]
    \caption{Unified Stochastic Momentum Methods: SUM$(f, \x_0, \eta, T)$}\label{alg:sum}
\begin{algorithmic}
    \STATE \textbf{Set} parameters:   
    $\rho\geq 0$ and $\beta\in (0,1)$.
    \FOR{$t=0,\ldots, T$}
    \STATE Compute a stochastic subgradient $\g_t$ for $f(\x_{t})$.
    \STATE $\y_{t+1} = \x_t-\eta \g_t$
    \STATE $\y_{t+1}^\rho = \x_t-\rho\eta\g_t$
    \STATE $\x_{t+1} = \y_{t+1} + \beta (\y_{t+1}^\rho - \y_{t}^\rho)$
    \ENDFOR
    \STATE \textbf{Output}: $\widehat{\x}_T = \sum_{t=0}^T \x_t/(T+1)$
\end{algorithmic}
\end{algorithm}

\begin{lemma}\label{lem:sum}
    For Algorithm~\ref{alg:sum}, assume  $f(\x) = \phi(\x) + \frac{1}{2\gamma}\|\x - \x_0\|^2$ is a $\lambda$-strongly convex function, $\g_t = \g(\x_t; \xi) + \frac{1}{\gamma}(\x_t - \x_0)$ where $\g(\x; \xi)\in\partial_F\phi(\x_t)$ such that $\E[\|\g(\x; \xi)\|^2]\leq G^2$,  and $\eta\leq (1-\beta)\gamma^2\lambda/(8\rho\beta+4)$, then we have that 
    \begin{align}\label{eqn:Es2}
        &\E[f(\widehat{\x}_T) - f(\x_*)]\leq \nonumber\\
       & \frac{(1-\beta)\|\x_0 - \x_*\|^2}{2\eta (T+1)} + \frac{\beta (f(\x_0)-f(\x_*))}{(1-\beta)(T+1)}  +  \frac{2\eta G^2(2\rho\beta+1) }{1-\beta} +\frac{4\rho\beta+4}{(1-\beta)}\frac{\eta}{\gamma^2}\|\x_0 - \x_*\|^2
    \end{align}
    where $\widehat{\x}_T = \sum_{t=0}^T\x_t/(1+T)$ and $\x_*\in\arg\min_{\x\in\R^d}f(\x)$.
\end{lemma}
{\bf Remark:} It is notable that in the above result, we do not use the bounded domain assumption since we consider $\Omega=\R^d$ for the unified momentum methods in this subsection. The key to  get rid of bounded domain assumption is by exploring the strong convexity of $f(\x)= \phi(\x) + \frac{1}{2\gamma}\|\x - \x_0\|^2$.

\begin{theorem} \label{thm:nsum}
    Suppose Assumption~\ref{ass:0} holds. By setting $\gamma=1/(2\mu), w_s=s^\alpha, \alpha > 0$, $\eta_{s} = (1-\beta)\gamma/(96s(\rho\beta+1))$, $T_s\geq 2304 (\rho\beta+1)s$, then we have 
\begin{align*}
    \E[\|\nabla\phi_{\gamma}(\x_\tau)\|^2] \leq  \frac{16\mu\Delta_\phi(\alpha+1)}{S+1} +  \frac{ (\beta G^2 + 96 G^2(2\rho\beta+1)(1-\beta))(\alpha+1)}{96(S+1)(2\rho\beta +1)(1-\beta)\alpha^{\mathbb I(\alpha<1)}},
\end{align*}
    where $\tau$ is randomly selected from $\{0, \ldots, S\}$ with probabilities $p_\tau\propto w_{\tau+1}, \tau=0, \ldots, S$.
\end{theorem}
{\bf Remark:} The bound in the above theorem is in the same order as that in Theorem~\ref{thm:sgd}.  The total iteration complexity for finding a solution $\x_\tau$ with $\E\big[\|\nabla\phi_\gamma(\x_\tau)\|^2\big]\leq \epsilon^2$ is $O(1/\epsilon^4)$. 

\subsection{Stagewise {\ada}}
In this subsection, we analyze  stagewise \textsc{AdaGrad} and  establish its adaptive complexity.  In particular, we consider the Meta algorithm that employs  \textsc{AdaGrad} in Algorithm~\ref{alg:adagrad}. The key difference of stagewise {\ada} from stagewise {\sgd} and stagewise SM is that the number of iterations $T_s$ at each stage is adaptive to the history of learning. It is this adaptiveness that makes the proposed stagewise {\ada} achieve adaptive convergence. It is worth noting that such adaptive scheme has been also considered in~\cite{SadaGrad18} for solving stochastic strongly convex problems. In contrast, we consider stochastic weakly convex problems. 
Similar to previous analysis of {\ada}~\cite{duchi2011adaptive,SadaGrad18}, we assume $\|g(\x;\xi)\|_\infty\leq G, \forall \x\in\Omega$ in this subsection. Note that this is stronger than Assumption~\ref{ass:0}-\ref{ass:1}. We formally state this assumption required in this subsection below. 
\begin{assumption}\label{ass:new}

       $\|g(\x; \xi)\|_\infty\leq G$ for any $\x\in\Omega$.
\end{assumption}


\begin{algorithm}[t]
    \caption{\textsc{AdaGrad}($ f, \x_0, \eta, *$)} \label{alg:adagrad}
    \begin{algorithmic}[1]
    \STATE \textbf{Initialize:} $\x_1 = \x_0$, $\g_{1:0}=[]$, $H_0\in\R^{d\times d}$
    \WHILE{{$T$ does not satisfy the condition in Theorem~\ref{thm:nadagrad}}}
    \STATE Compute a stochastic subgradient $\g_t $ for $f(\x_t)$
    \STATE Update $g_{1:t} = [g_{1:t-1}, \g(\x_t)]$, $s_{t,i}  = \|g_{1:t,i}\|_2$
    \STATE Set $H_t = H_0 + \diag(\s_t)$ and $\psi_t(\x) = \frac{1}{2}(\x-\x_1)^{\top}H_t(\x-\x_1)$
    \STATE Let $\x_{t+1}=\arg\min\limits_{\x\in\Omega}\eta \x^{\top}\left(\frac{1}{t}\sum_{\tau=1}^t\g_\tau\right) + \frac{1}{t}\psi_t(\x)$
    \ENDWHILE
    \STATE \textbf{Output}:  $\xh_T=\sum_{t=1}^{T}\x_t/T$
    \end{algorithmic}
    \end{algorithm}

The convergence analysis of stagewise {\ada} is build on the following lemma, which is attributed to~\cite{SadaGrad18}.
\begin{lemma}\label{lem:adagrad}
    Let $f(\x)$ be a convex function, $H_0=G I$ with $G\geq\max_t \|\g_t\|_\infty$,  and iteration number $T$ satisfy $T\geq M\max\{\frac{G+ \max_i\|g_{1:T,i}\|}{2},  \sum_{i=1}^d\|g_{1:T,i}\|\}$. 
    Algorithm~\ref{alg:adagrad} returns an averaged solution $\xh_T$ such that
    \begin{align}\label{eqn:Es}
        \E[ f(\xh_T)-f(\x_*)] \leq \frac{1}{M\eta}\|\x_{0}-\x_*\|^2 +\frac{\eta}{M},
    \end{align}
    where $\x_*=\arg\min_{\x\in\Omega}f(\x)$, $g_{1:t}=(\g(\x_1),\ldots, \g(\x_t))$ and $g_{1:t,i}$ denotes the $i$-th row of $g_{1:t}$.
\end{lemma}


The convergence property of stagewise \textsc{AdaGrad} is described by following theorem.
\begin{theorem} \label{thm:nadagrad}
    Suppose Assumption~\ref{ass:0}, Assumption~\ref{ass:0-sgd} and Assumption~\ref{ass:new} hold. By setting $\gamma=1/(2\mu), w_s=s^\alpha, \alpha >0$, $\eta_{s} = c/ \sqrt{s}$, $T_s\geq M_s\max\{\hat G+ \max_i\|g^s_{1:T_s,i}\|,  \sum_{i=1}^d\|g^s_{1:T_s,i}\|\}$ where $c>0$ is a free parameter, and $M_s\eta_s \geq  24\gamma $, then we have 
\begin{align*}
    \E[\|\nabla\phi_{\gamma}(\x_\tau)\|^2]     \leq& \frac{16\mu\Delta_\phi(\alpha+1)}{S+1} + \frac{4\mu^2c^2(\alpha+1)}{(S+1)\alpha^{\mathbb I(\alpha<1)}},
\end{align*}
    where $\hat G = G + \gamma^{-1}D$, and $g^s_{1:t, i}$ denotes the cumulative stochastic gradient of the $i$-th coordinate at the $s$-th stage. 
\end{theorem}
{\bf Remark:} It is obvious that the total number of iterations $\sum_{s=1}^ST_s$ is adaptive to the data. Next, let us present more discussion on the iteration complexity. Note that $M_s = O(\sqrt{s})$. By the boundness of stochastic gradient $\|g^{1:T_s, i}\|\leq O(\sqrt{T_s})$, therefore $T_s$ in the order of $O(s)$ will satisfy the condition in Theorem~\ref{thm:nadagrad}. Thus in the worst case, the iteration complexity for finding $ \E[\|\nabla\phi_{\gamma}(\x_\tau)\|^2]\leq \epsilon^2$ is in the order of $\sum_{s=1}^S O(s)\leq O(1/\epsilon^4)$. To show the potential advantage of adaptive step size as in the convex case, let us consider a good case when the cumulative growth of stochastic gradient is slow, e.g., assuming $\|g^s_{1:T_s,i}\|\leq O({T_s}^{\alpha})$ with $\alpha<1/2$. Then $T_s = O(s^{1/(2(1-\alpha))})$ will work, and then the total number of iterations $\sum_{s=1}^S T_s \leq S^{1+1/(2(1-\alpha))}\leq O(1/\epsilon^{2+1/(1-\alpha)})$, which is better than $O(1/\epsilon^4)$. Finally, we remark that the bounded domain assumption could be removed similar to last subsection. 

\vspace*{0.1in}

\subsection{Stagewise Stochastic ADMM for Solving Problems with Structured Regularizers}
In this subsection, we consider solving a regularized problem with a structured regularizer, i.e., 
\begin{align}\label{eqn:struc}
\min_{\x\in\Omega} \phi(\x): = \E[f(\x; \xi)] + \psi(A\x),
\end{align}
where $A\in\R^{d\times m}$, and $\psi(\cdot): \R^m\rightarrow\R$ is some convex structured regularizer (e.g., generalized Lasso $\psi(A\x) = \|A\x\|_1$). We assume that $\phi(\cdot)$ is $\mu$-weakly convex. Although Stagewise {\sgd} can be employed to solve the above problem, it is usually expected to generate a sequence of solutions that respect certain properties (e.g., sparsity) promoted by the regularizer. When $\E[f(\x; \xi)]$ is convex, the problem is usually solved by stochastic ADMM shown in Algorithm~\ref{alg:admm} (assuming $f(\x) = \E[f(\x; \xi)] + \psi(A\x)$), in which the following  steps are alternatively executed: 
\begin{align}
\x_{\tau+1} & = \arg\min_{\x\in\Omega}\partial f(\x_t, \xi_{\tau})^{\top}\x + \frac{\beta}{2}\left\|(A\x - \y_\tau)- \frac{1}{\beta}\lambda_\tau\right\|^2 + \frac{\|\x - \x_\tau\|^2_{C}}{\eta},\label{eqn:admm1}\\
\y_{\tau+1} & = \arg\min_{\x\in\Omega}\psi(\y) + \frac{\beta}{2}\left\|(A\x_{\tau+1} - \y)- \frac{1}{\beta}\lambda_t\right\|^2,\label{eqn:admm2}\\
\lambda_{\tau+1} & = \lambda_\tau - \beta(A\x_{\tau+1} - \y_{\tau+1}),\label{eqn:admm3}
\end{align}
where $\beta>0$ is the penalty parameter of ADMM, $\|\x\|^2_C =\x^{\top}C\x$, and $C = \alpha I - \eta\beta A^{\top}A\succeq I$ with some appropriate $\alpha>0$.

\begin{algorithm}[t]
\caption{SADMM($f, \x_0, \eta, \beta, t$)} \label{alg:admm}
\begin{algorithmic}[1]
\STATE \textbf{Input}: $\x_0\in\R^d$, a step size $\eta$, penalty parameter $\beta$, the number of iterations $t$ and a domain $\Omega$.
\STATE \textbf{Initialize:} $\x_1 = \x_0, \y_1 = A\x_1, \lambda_1=0$
\FOR{$\tau=1,\ldots, t$}
\STATE Update $
\x_{\tau+1}$ by~(\ref{eqn:admm1}) 
\STATE Update $\y_{\tau+1}$ by~(\ref{eqn:admm2}) 
\STATE Update  $\lambda_{\tau+1}$ by~(\ref{eqn:admm3})
\ENDFOR
\STATE \textbf{Output}:  $\xh_{t} = \sum_{\tau=1}^t\x_\tau/t$
\end{algorithmic}
\end{algorithm}

In order to employ SADMM for solving~(\ref{eqn:struc}) with a weakly convex objective, we use $\hat f_s(\cdot; \xi) =f(\cdot; \xi) + \frac{1}{2\gamma}\|\x - \x_{s-1}\|^2$ to define $f_s(\x) = \E[\hat f_s(\x; \xi)] + \psi(A\x)$ in the $s$-th call of SADMM in the Meta framework. 

A convergence upper bound  of stochastic ADMM for solving
\begin{align}\label{eqn:struc2}
\min_{\x\in\Omega}  f(\x) = \E[f(\x; \xi)] + \frac{1}{2\gamma}\|\x - \x_0\|^2 + \psi(A\x),
\end{align}
is given in the following lemma. 
\begin{lemma}(Corollary 3~\cite{DBLP:conf/nips/XuLLY17})\label{cor:G}
  For Algorithm~\ref{alg:admm}, assume  $f(\x)$ is a  convex function and $\psi(\cdot)$ is a  $\rho$-Lipschitz continuous convex function, $\g(\x_t; \xi)\in\partial_Ff(\x_t; \xi_t)$ is used in the update, $C = \alpha I - \eta \beta A^{\top}A\succeq I$, and Assumption~\ref{ass:0-sgd} holds. Then, 
\begin{align*}
\E[f(\xh_t) - f(\x_*)]\leq& \frac{\alpha\|\x_0 - \x_*\|_2^2}{2\eta t}+  \frac{\beta\|A\|_2^2\|\x_0 - \x_*\|_2^2}{2t} + \frac{\rho^2}{2\beta t} +  \frac{\eta \hat G^2}{2}  + \frac{\rho\|A\|_2D}{t}.
\end{align*}
where $\hat G = G + \gamma^{-1}D$.
\end{lemma}

\begin{theorem} \label{thm:nadmm}
    Suppose Assumption~\ref{ass:0} and Assumption~\ref{ass:0-sgd} hold and SADMM$(f_s, \x_{s-1}, \eta_s, \beta_s, T_s)$ is employed in the Meta Algorithm~\ref{alg:meta}. By setting $\gamma=1/(2\mu), w_s=s^\alpha, \alpha > 0$, $\eta_{s} = c_1/s$, $\beta_s = c_2s$, $T_s\geq 24s\gamma\max(\alpha/c_1, c_2\|A\|_2^2) $, where $c_1, c_2>0$, then we have 
\begin{align*}
    \E[\|\nabla\phi_{\gamma}(\x_\tau)\|^2] \leq  \frac{16\mu\Delta_\phi(\alpha+1)}{S+1} +  \frac{C(\alpha+1)}{(S+1)\alpha^{\I(\alpha<1)}},
\end{align*}
    where $\tau$ is randomly selected from $\{0, \ldots, S\}$ with probabilities $p_\tau\propto w_{\tau+1}, \tau=0, \ldots, S$, and $C$ is some constant depending on $c_1, c_2, \rho, G, D, \|A\|_2$.
\end{theorem}
{\bf Remark:} The above result can be easily proved. Therefore, the proof is omitted.

\section{Experiments}
\vspace*{-0.12in}
In this section, we present some empirical results to verify the effectiveness  of the proposed stagewise algorithms. We use two benchmark datasets, namely CIFAR-10 and CIFAR-100~\cite{cifar10} for our experiments. We implement the proposed stagewise algorithms in TensorFlow. We compare different algorithms for learning ResNet-20~\cite{he2016deep} with batch normalization~\cite{DBLP:conf/icml/IoffeS15} adopted after each convolution and before ReLU activation.

   \begin{figure}[t] 
\centering
	\hspace*{-0.1in} 
	\includegraphics[scale=0.22]{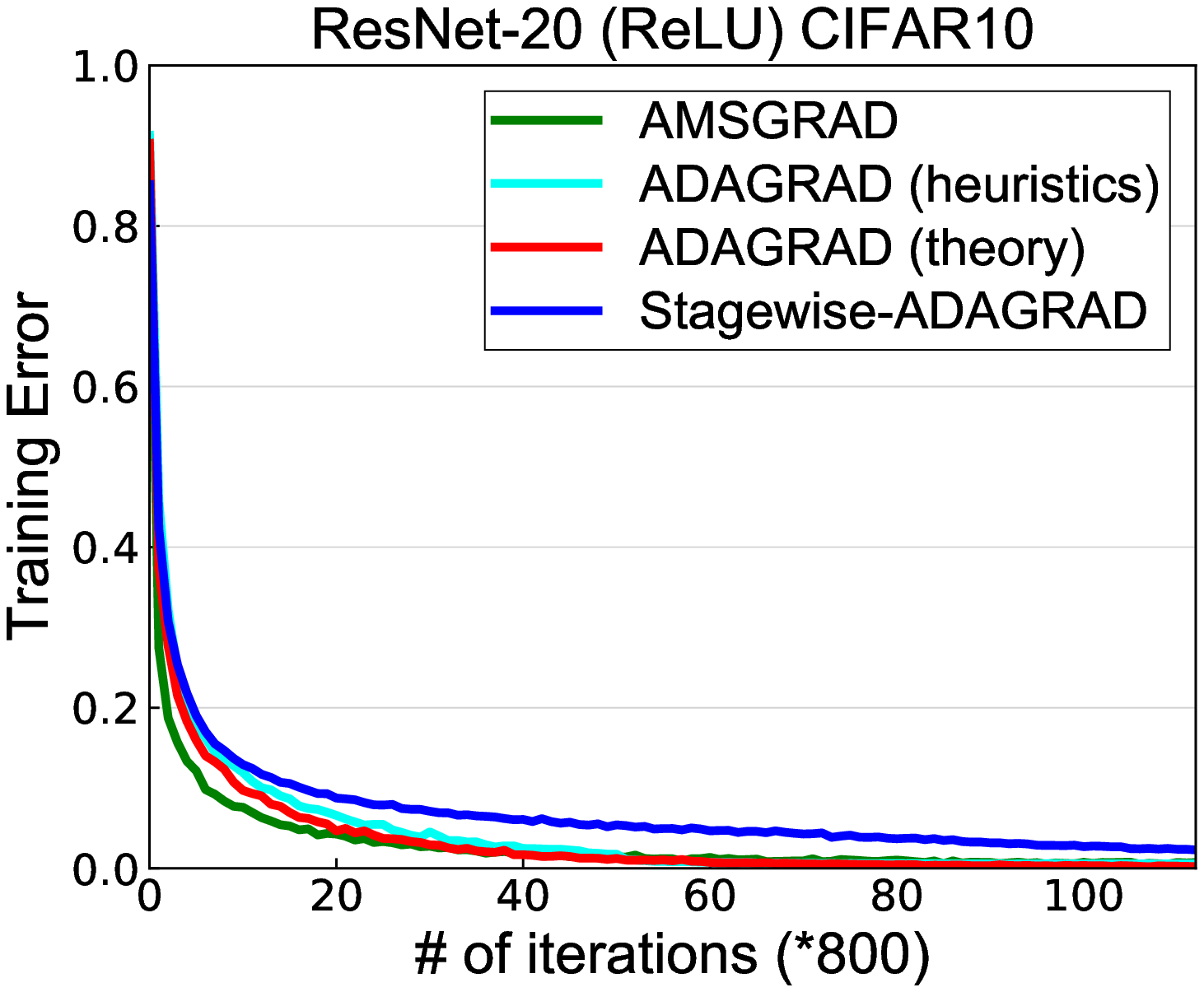}
    \hspace*{-0.1in} 
	\includegraphics[scale=0.22]{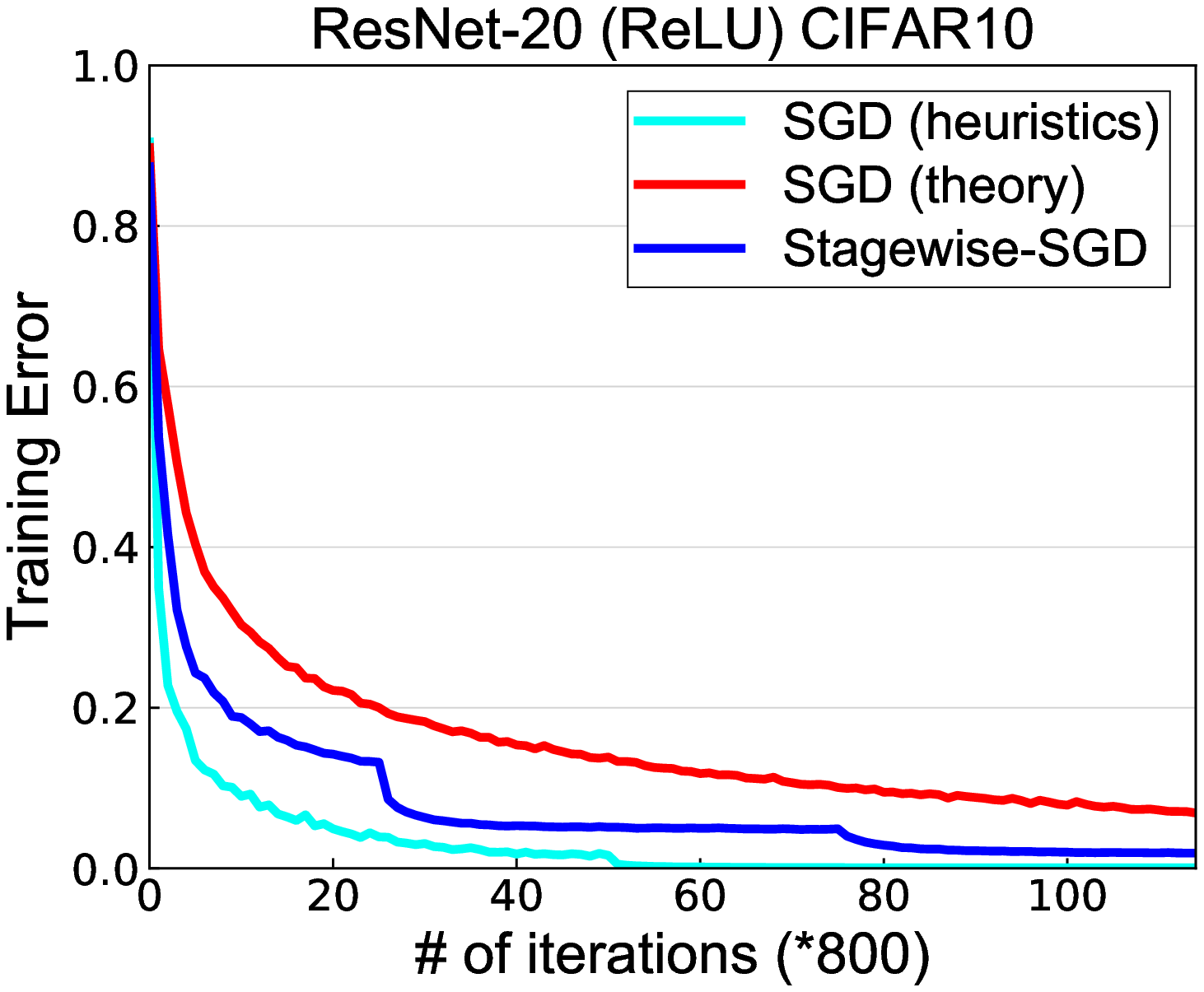}
    \hspace*{-0.1in} 
	\includegraphics[scale=0.22]{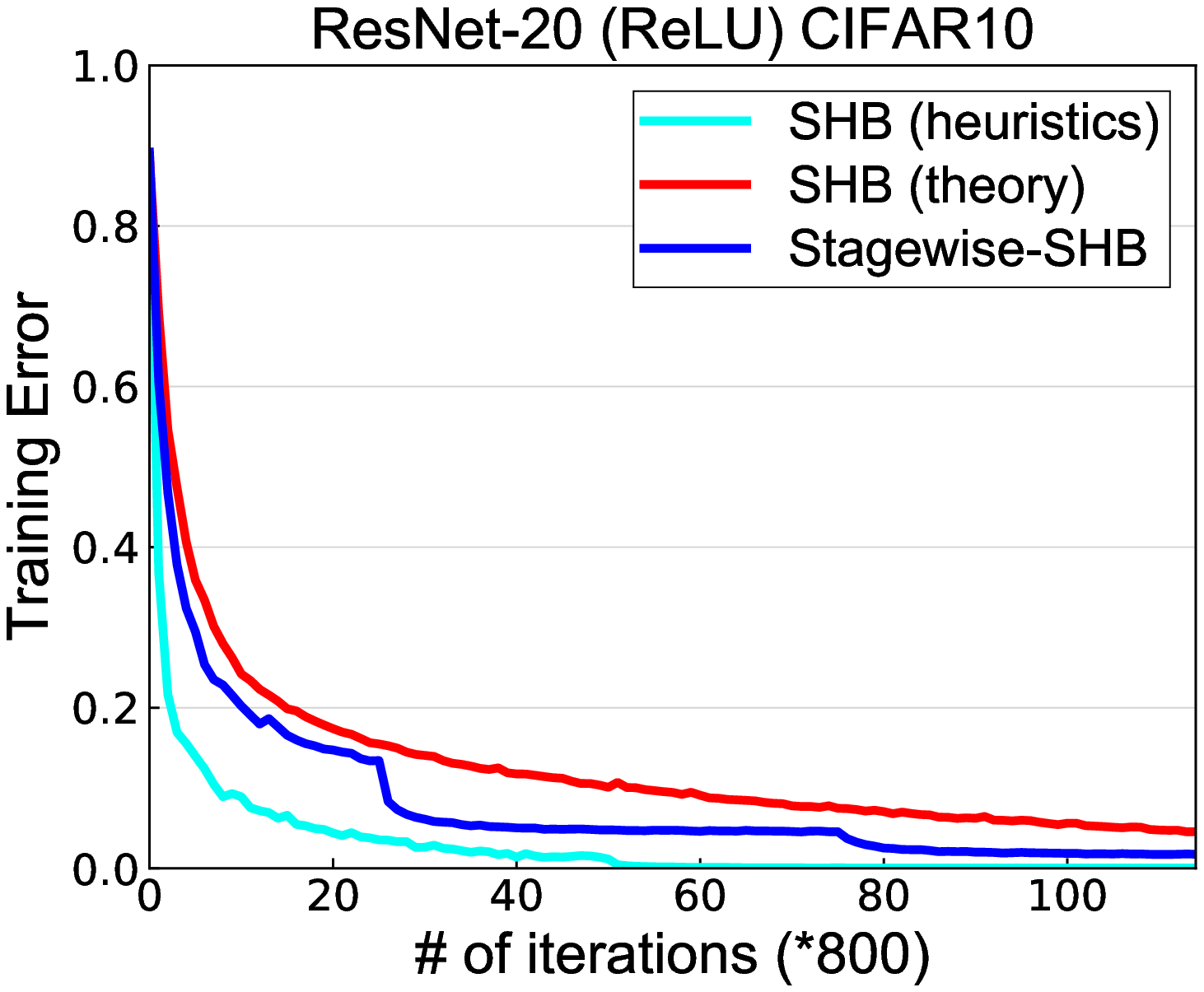}
	    \hspace*{-0.1in} 
	\includegraphics[scale=0.22]{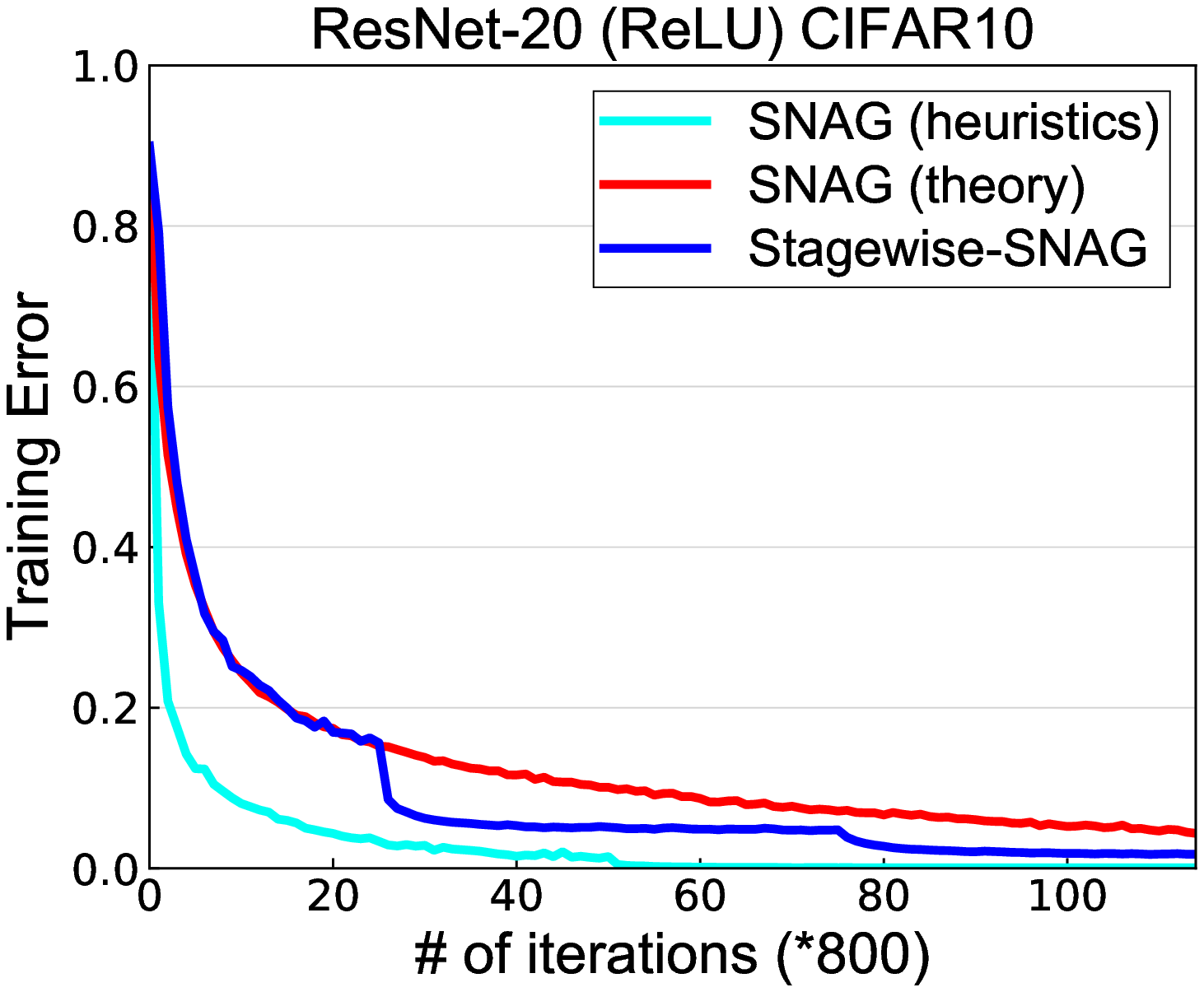}
	
	\hspace*{-0.1in}
		\includegraphics[scale=0.22]{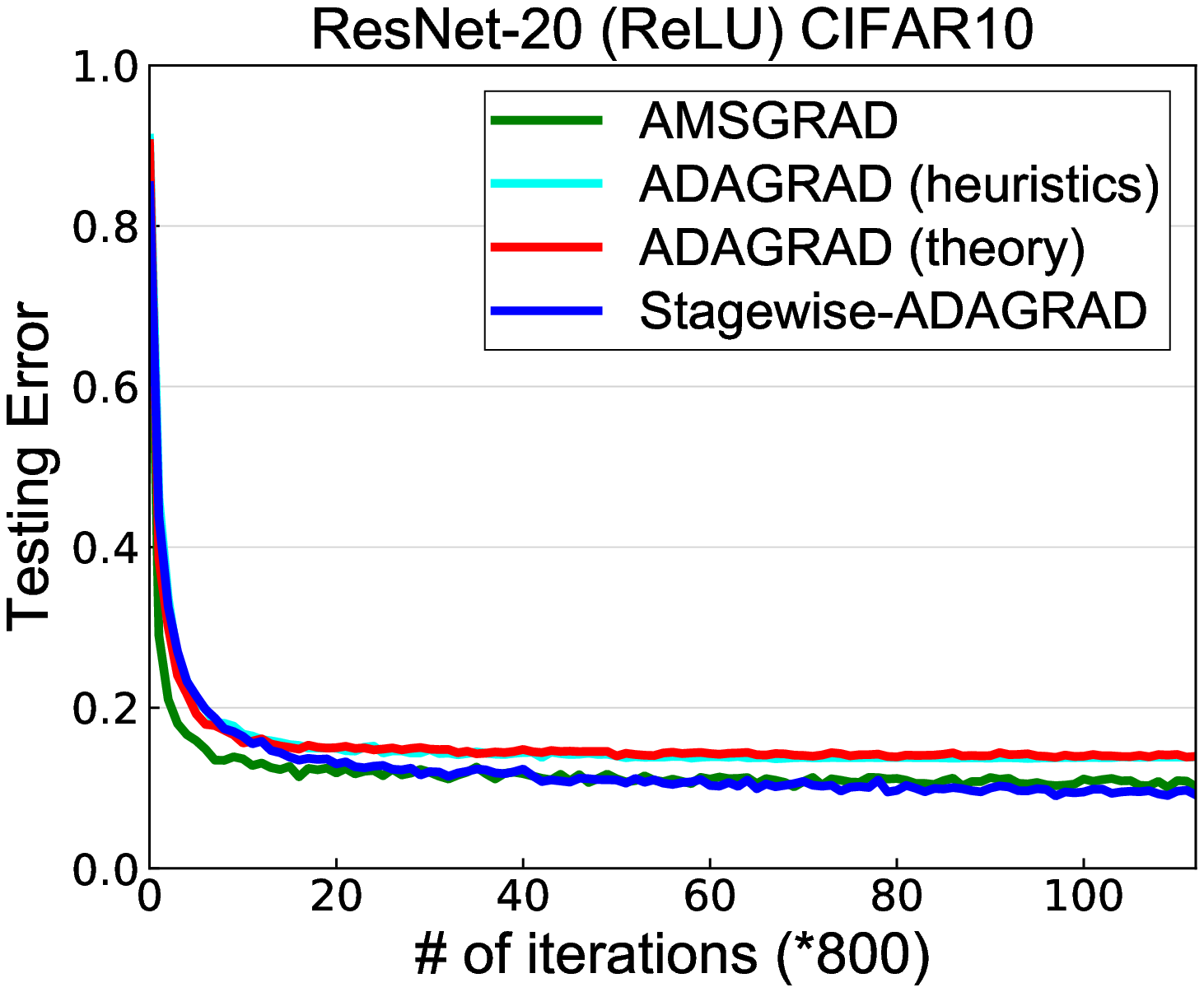}
    \hspace*{-0.1in} 
	\includegraphics[scale=0.22]{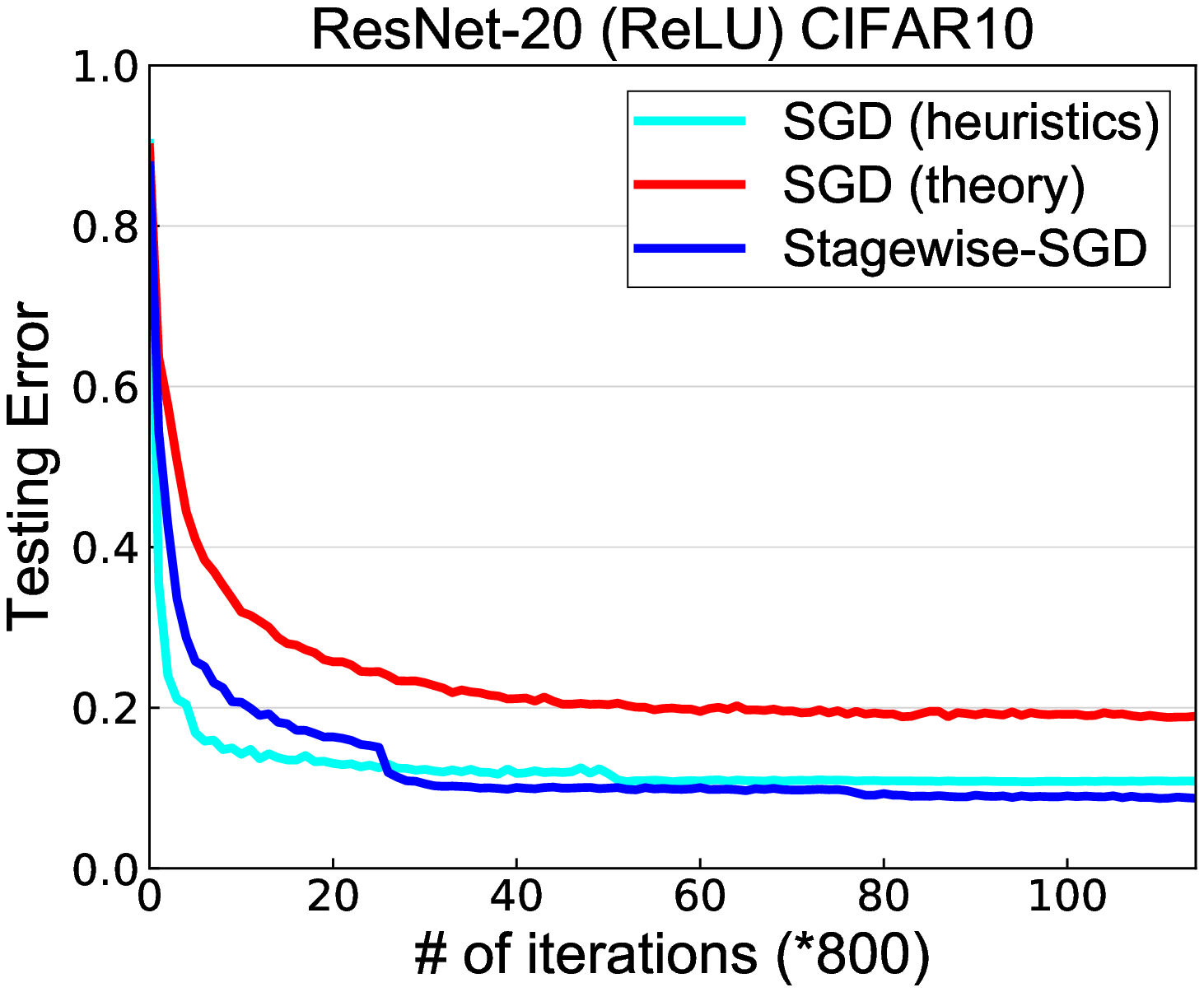}
    \hspace*{-0.1in} 
	\includegraphics[scale=0.22]{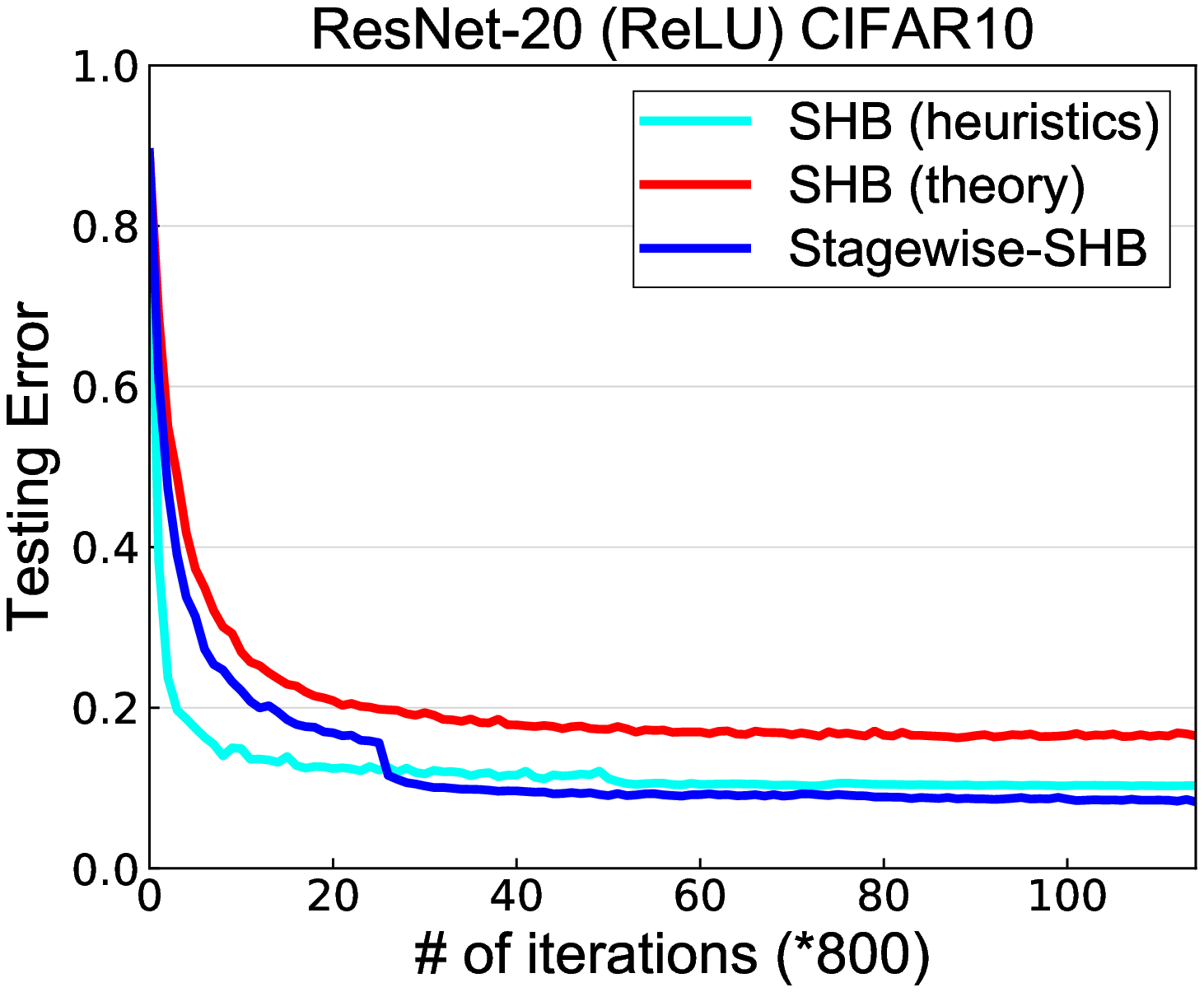}
	    \hspace*{-0.1in} 
	\includegraphics[scale=0.22]{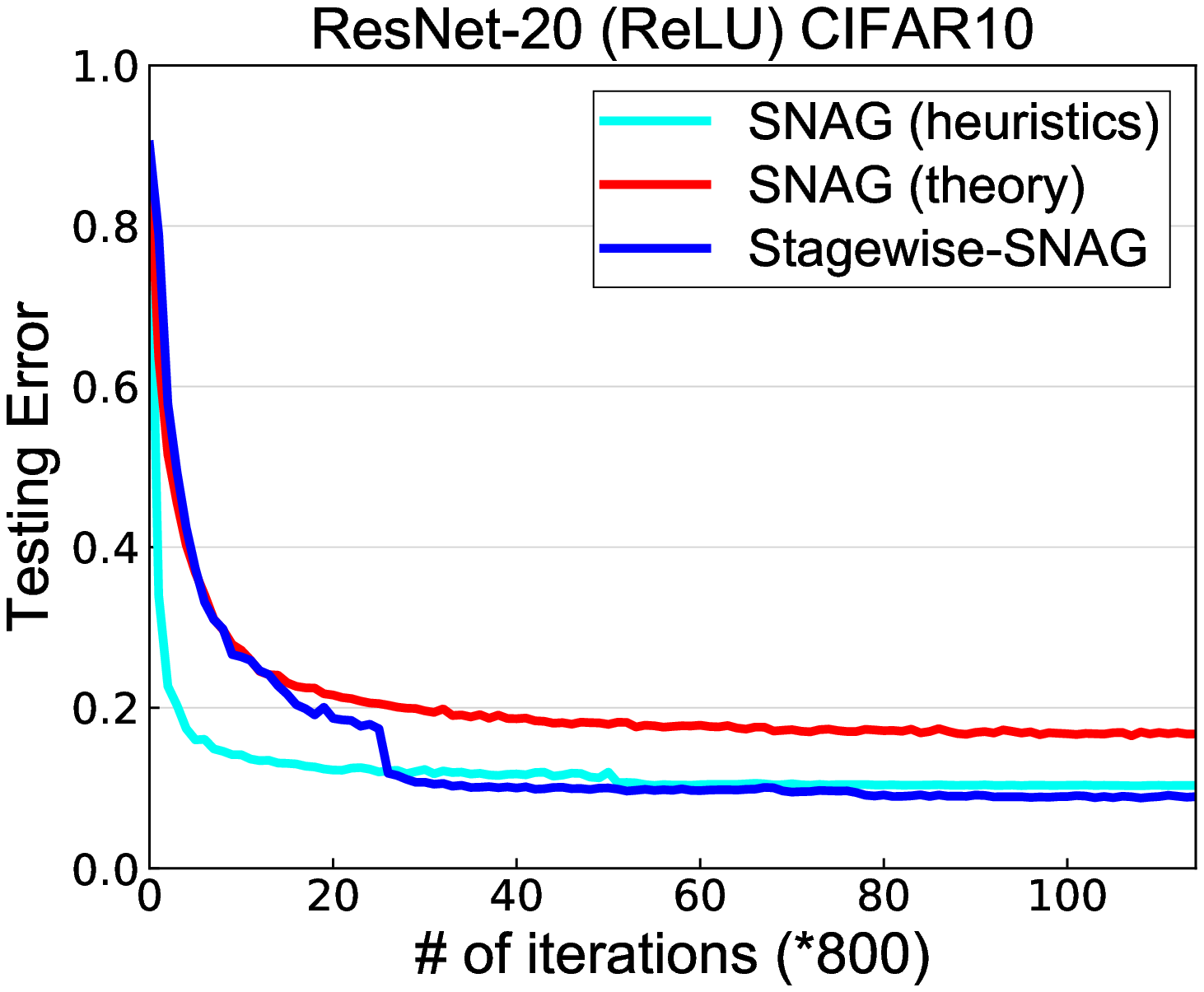}
    \vspace*{-0.15in}
\caption{Comparison of Training Error (Top) and Testing Error (bottom) on CIFAR-10 without Regularization.}
\label{fig:2}
\vspace*{0.15in}
	\hspace*{-0.2in} 
	\includegraphics[scale=0.22]{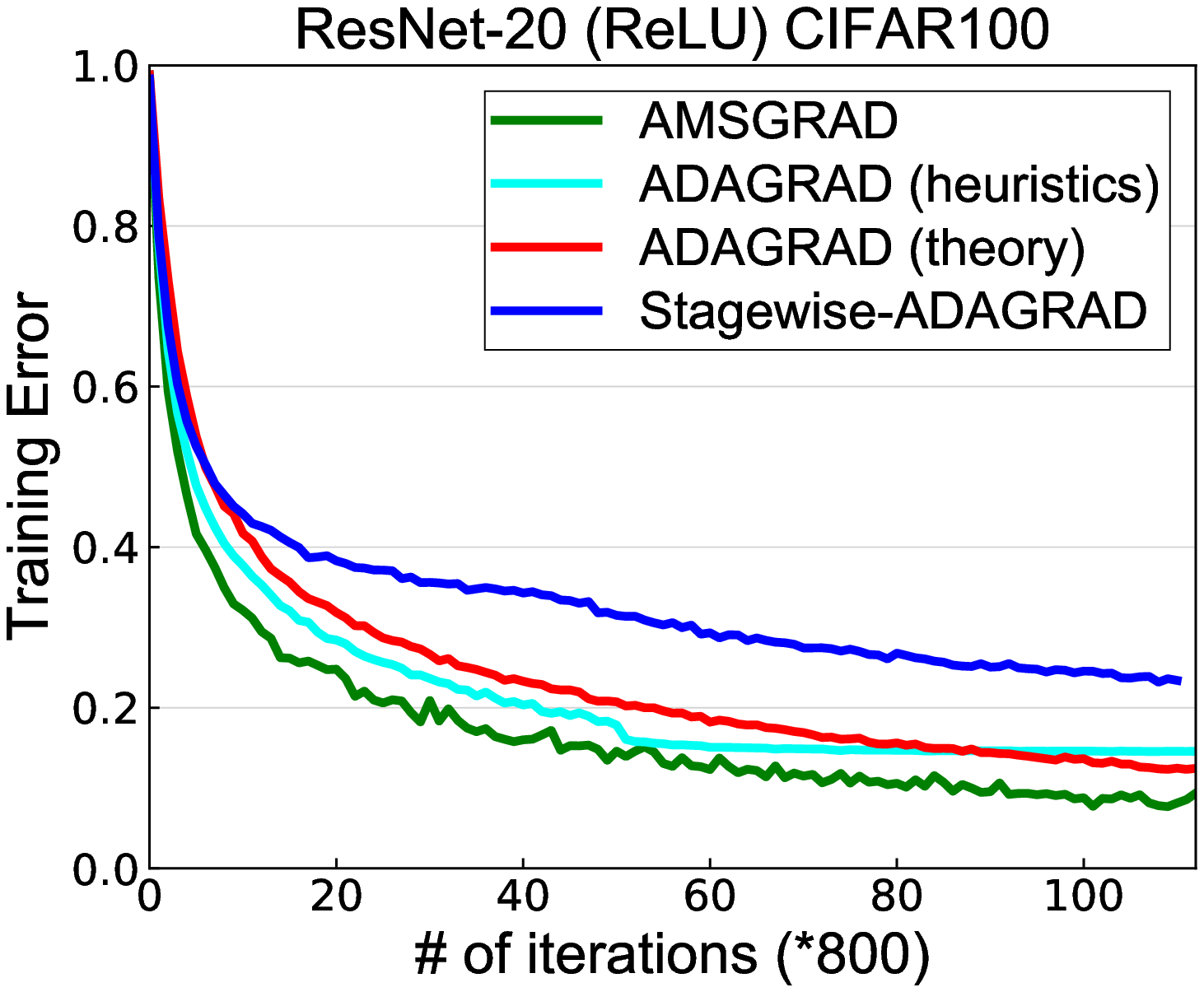}
    \hspace*{-0.1in} 
	\includegraphics[scale=0.22]{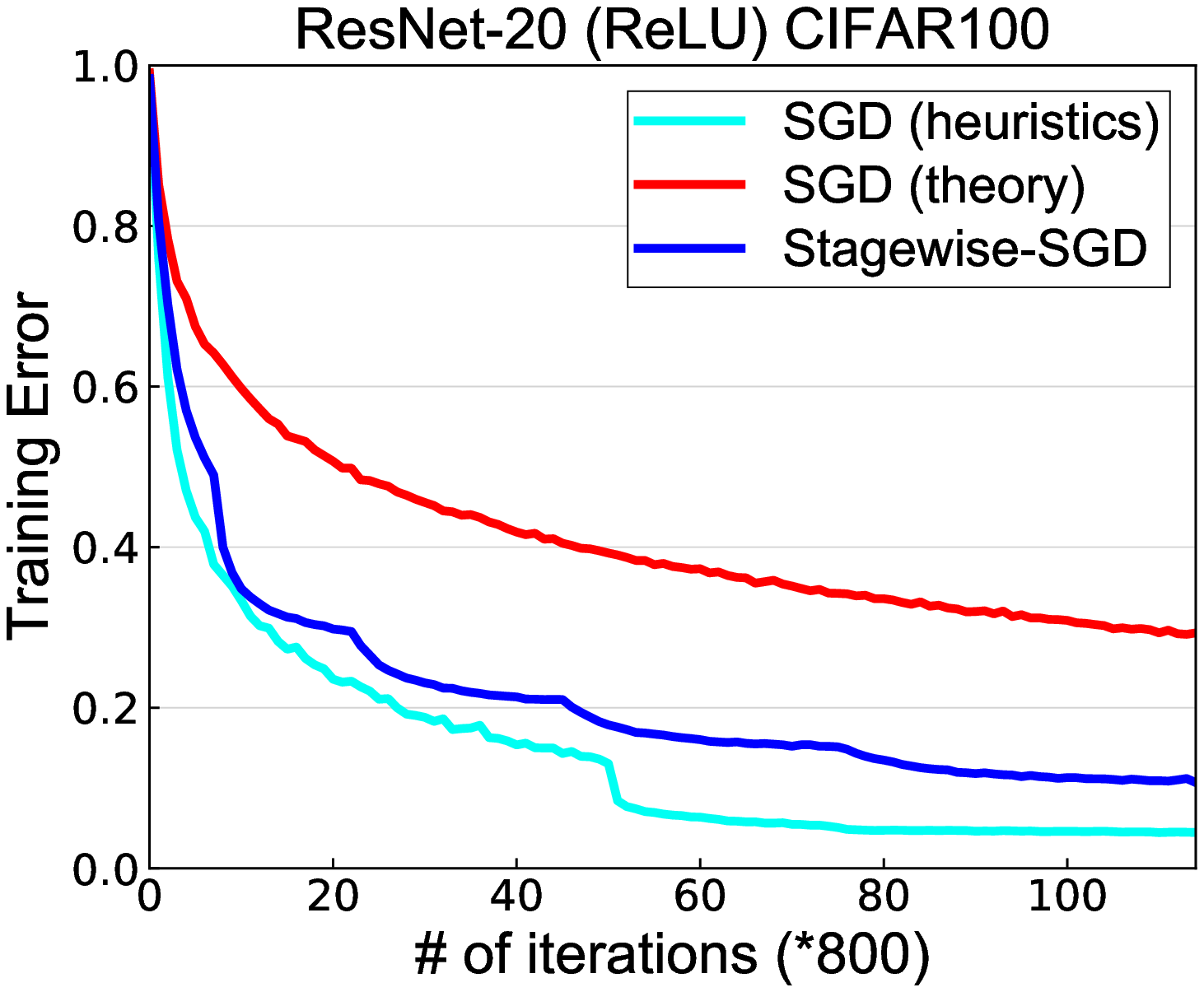}
    \hspace*{-0.1in} 
	\includegraphics[scale=0.22]{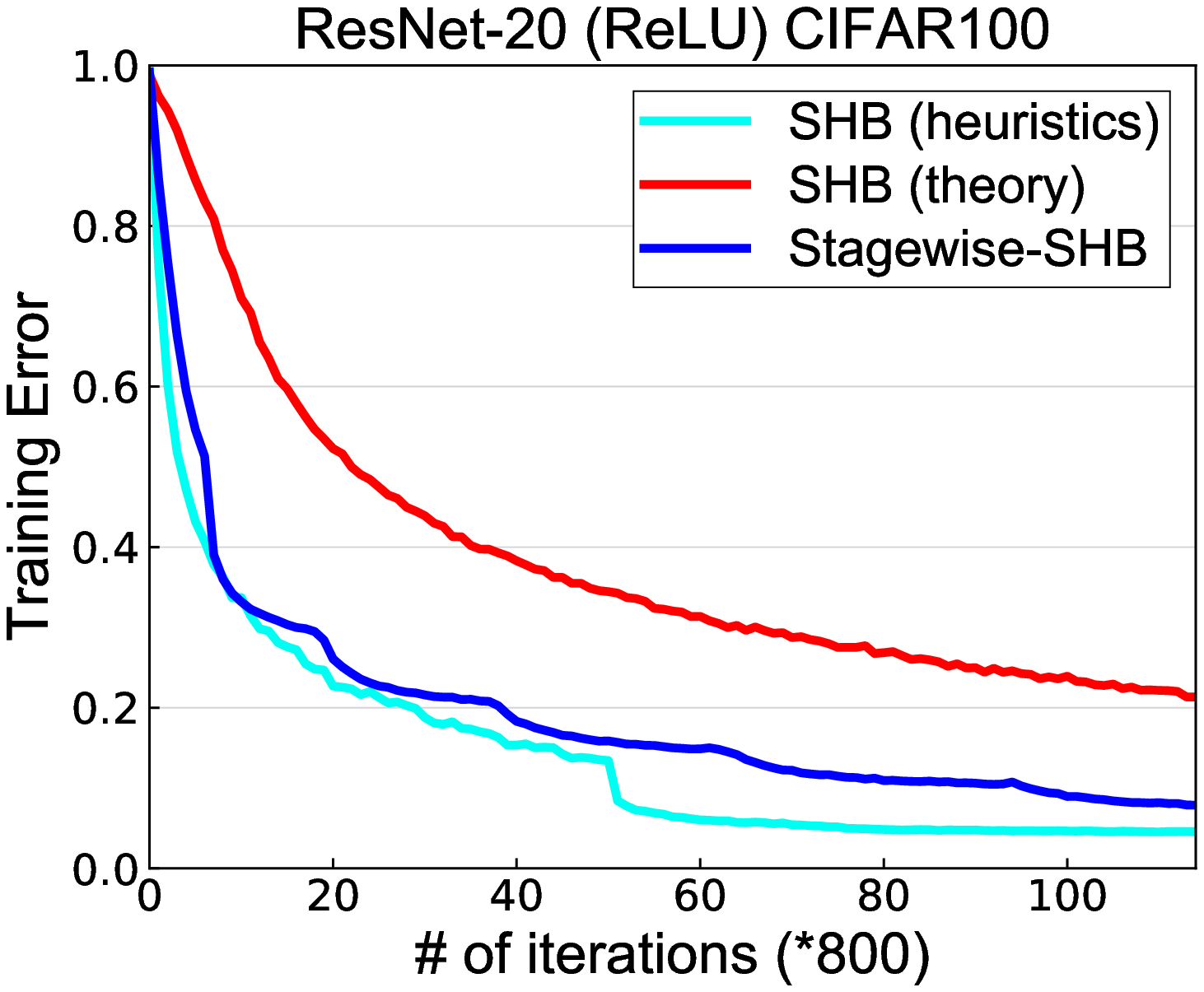}
	    \hspace*{-0.1in} 
	\includegraphics[scale=0.22]{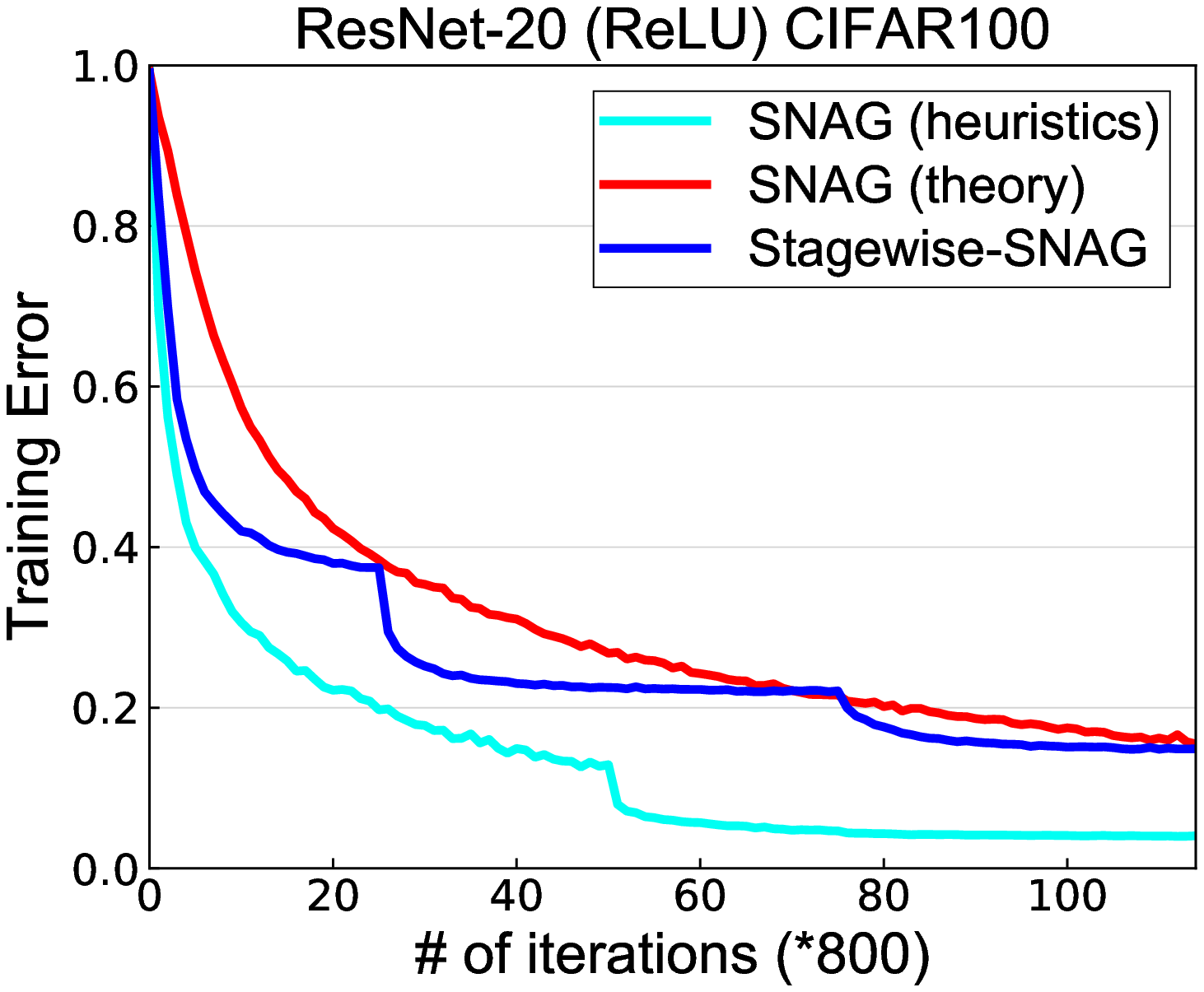}
	
	\hspace*{-0.1in} 
    \includegraphics[scale=0.22]{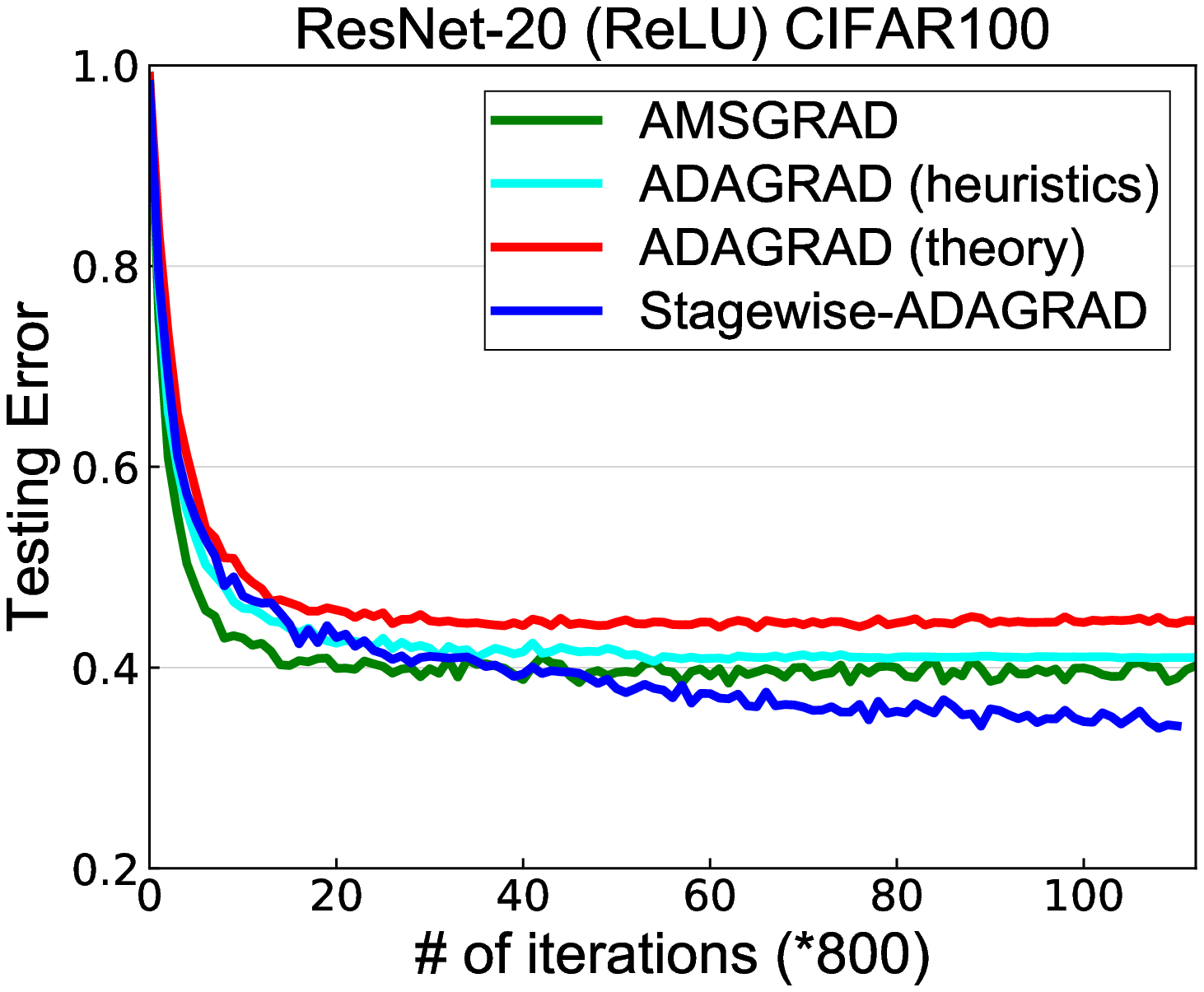}
    \hspace*{-0.1in} 
	\includegraphics[scale=0.22]{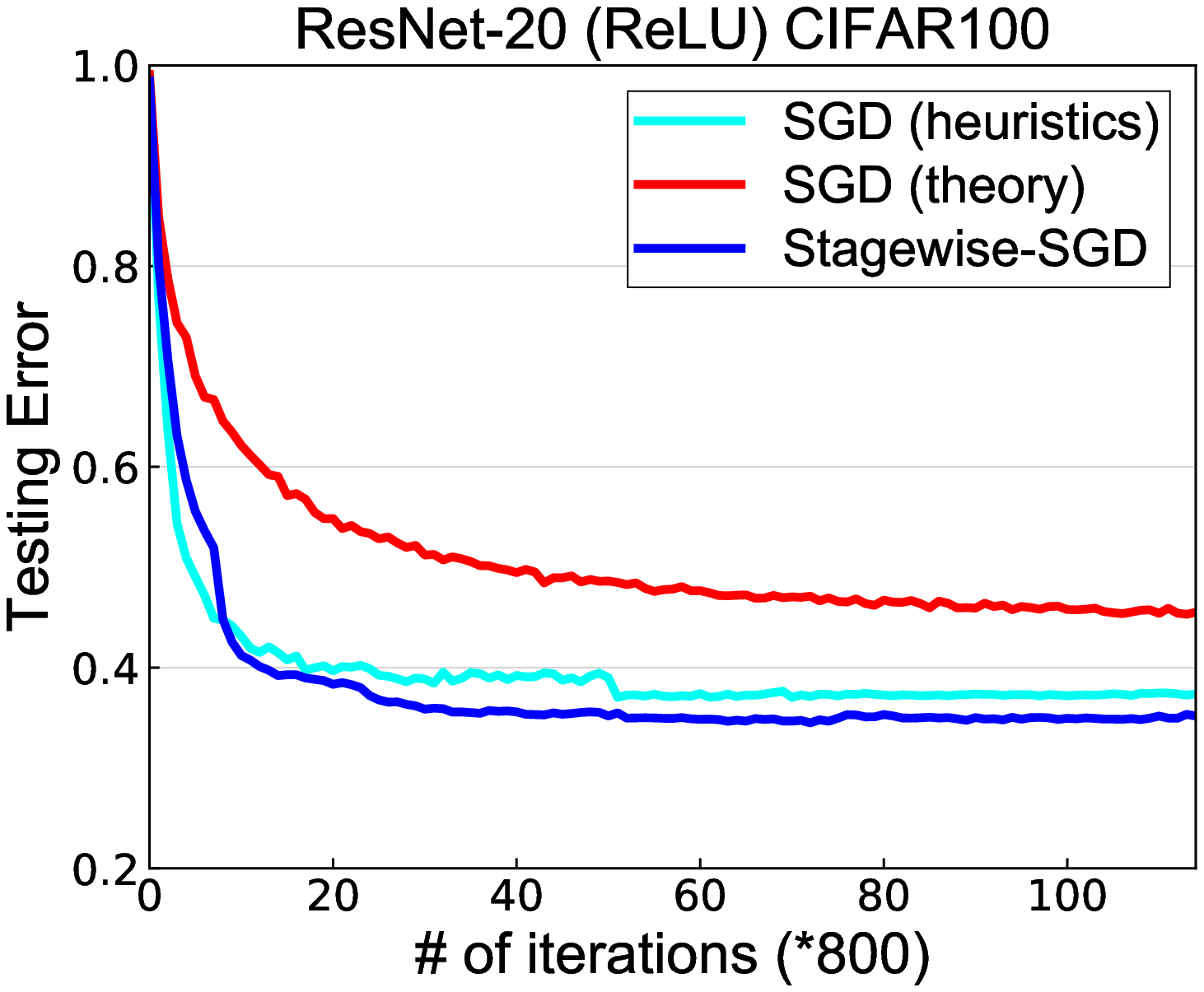}
    \hspace*{-0.1in} 
	\includegraphics[scale=0.22]{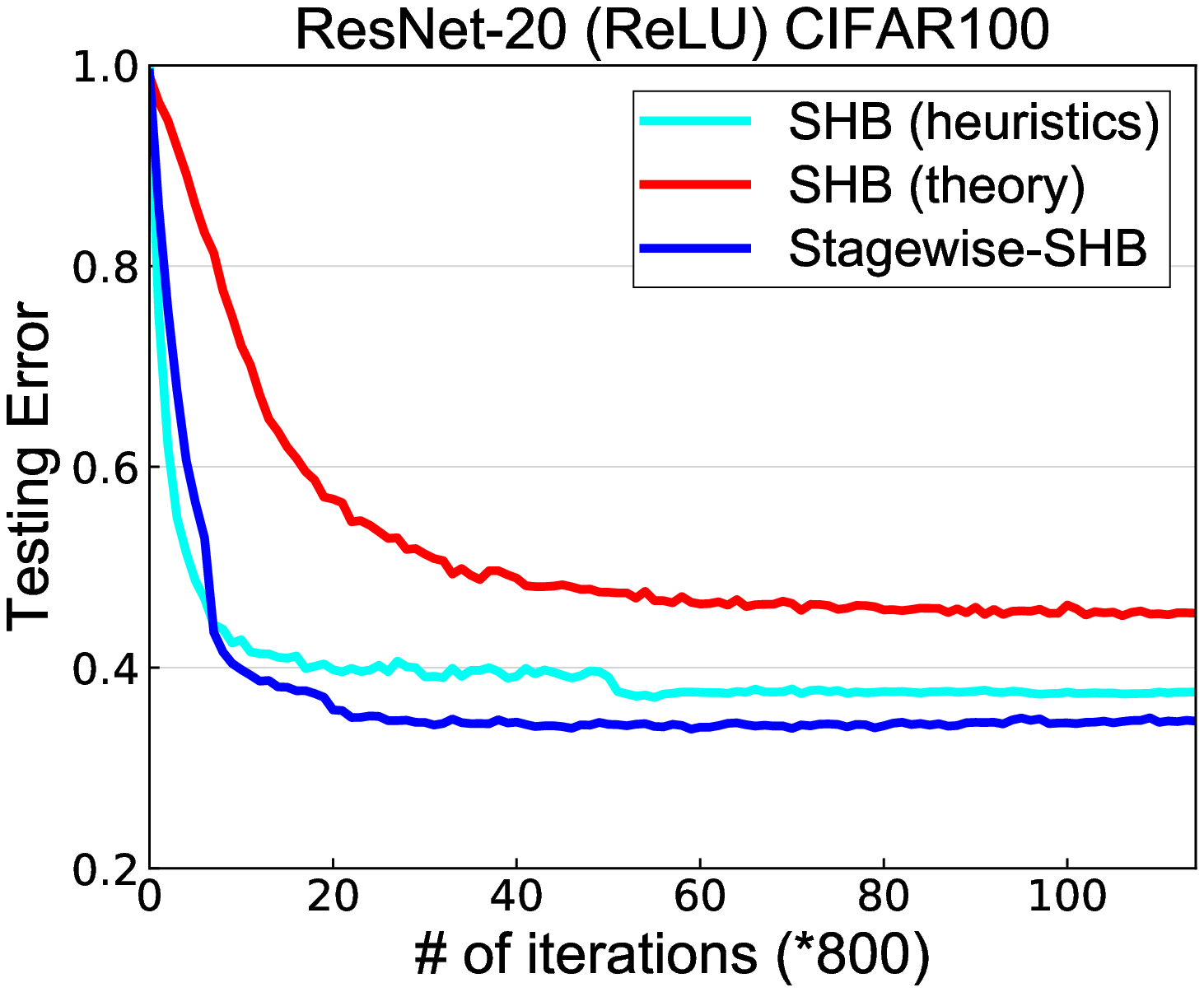}
	    \hspace*{-0.1in} 
	\includegraphics[scale=0.22]{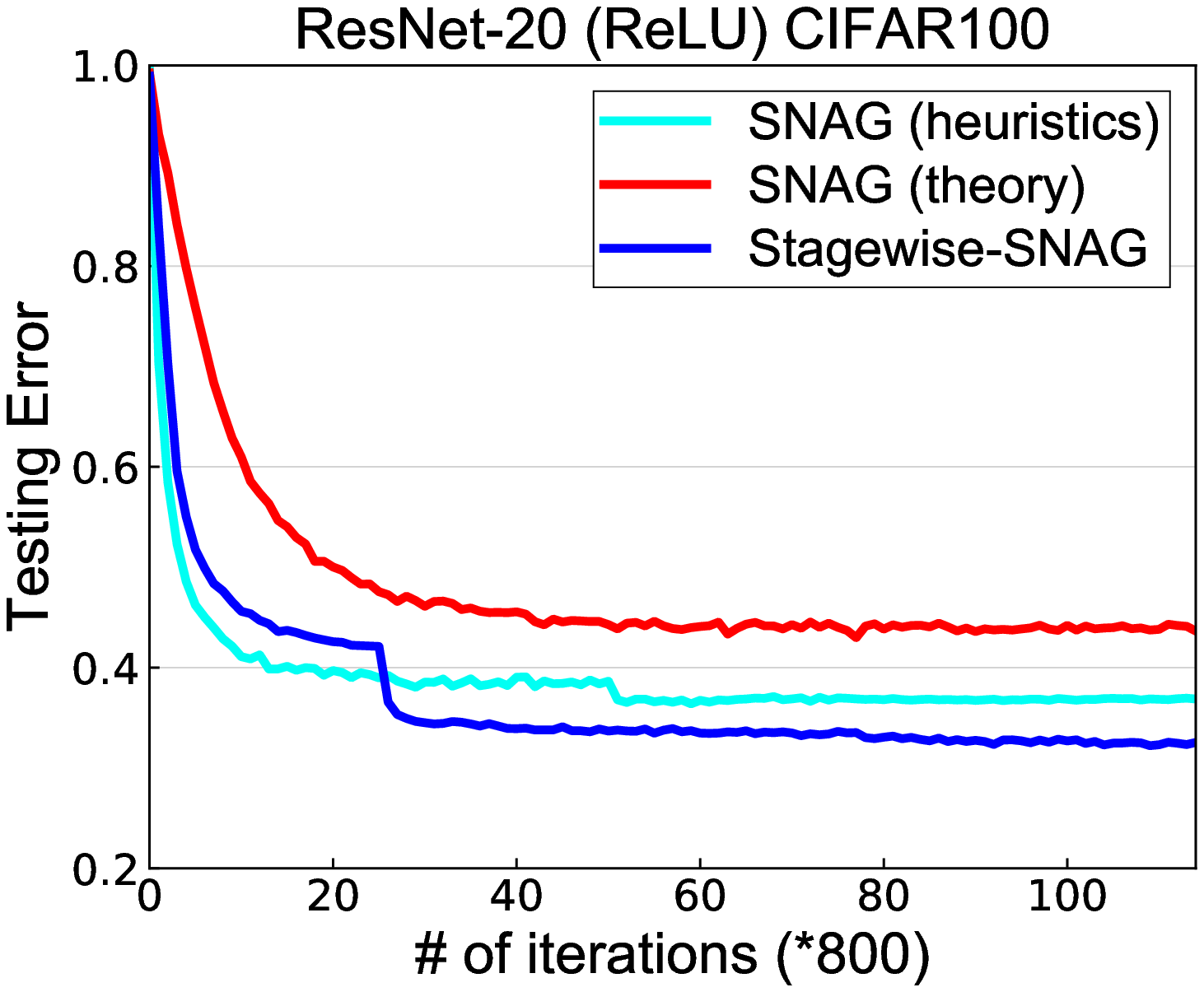}
    \vspace*{-0.15in}
\caption{Comparison of Training Error (Top) and Testing Error (bottom) on CIFAR-100 without Regularization.}

\label{fig:3}
\vspace*{0.15in}
\end{figure}

\begin{figure}[t] 
\centering

	\hspace*{-0.1in} 
	\includegraphics[scale=0.22]{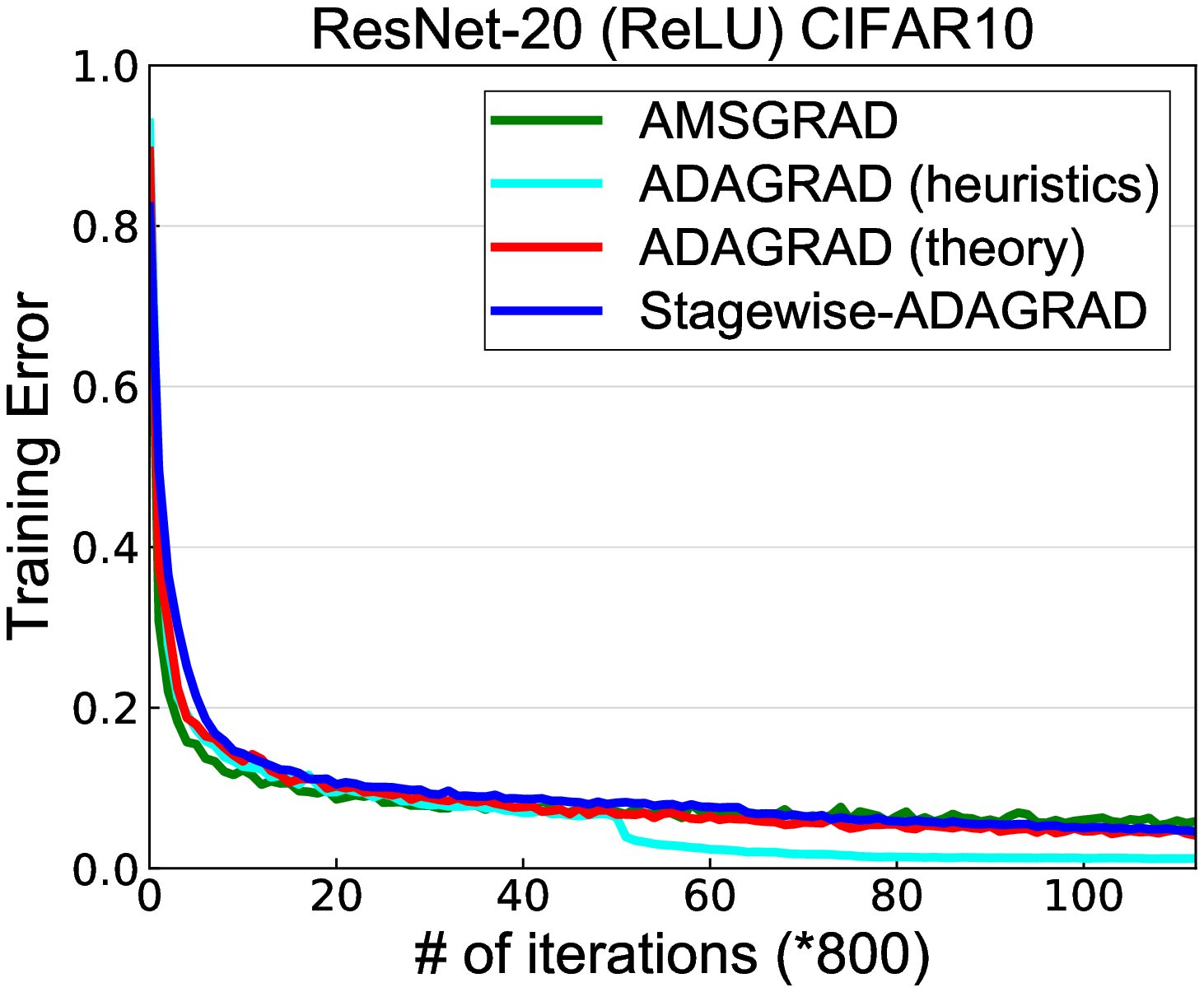}
    \hspace*{-0.1in} 
	\includegraphics[scale=0.22]{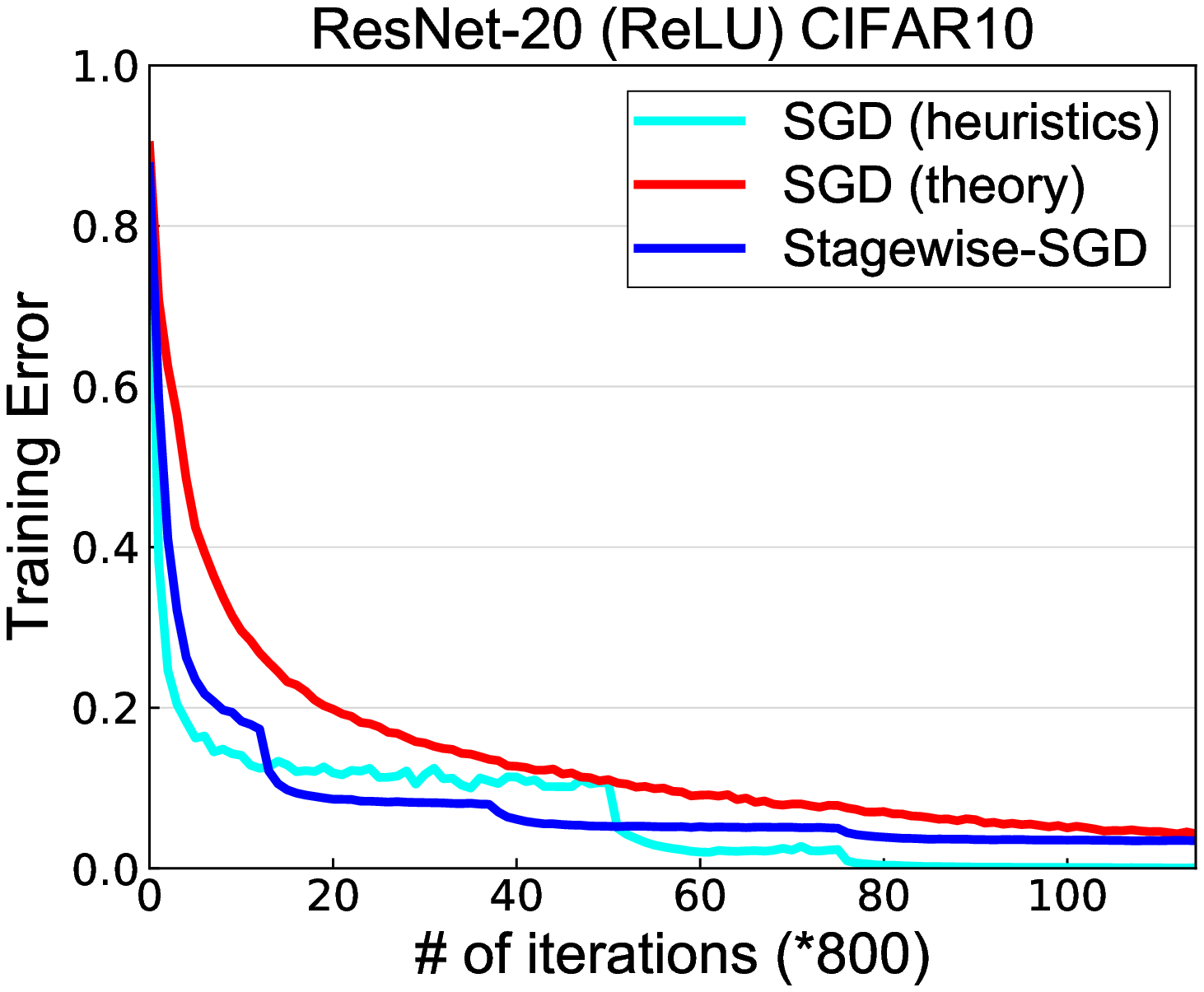}
    \hspace*{-0.1in} 
	\includegraphics[scale=0.22]{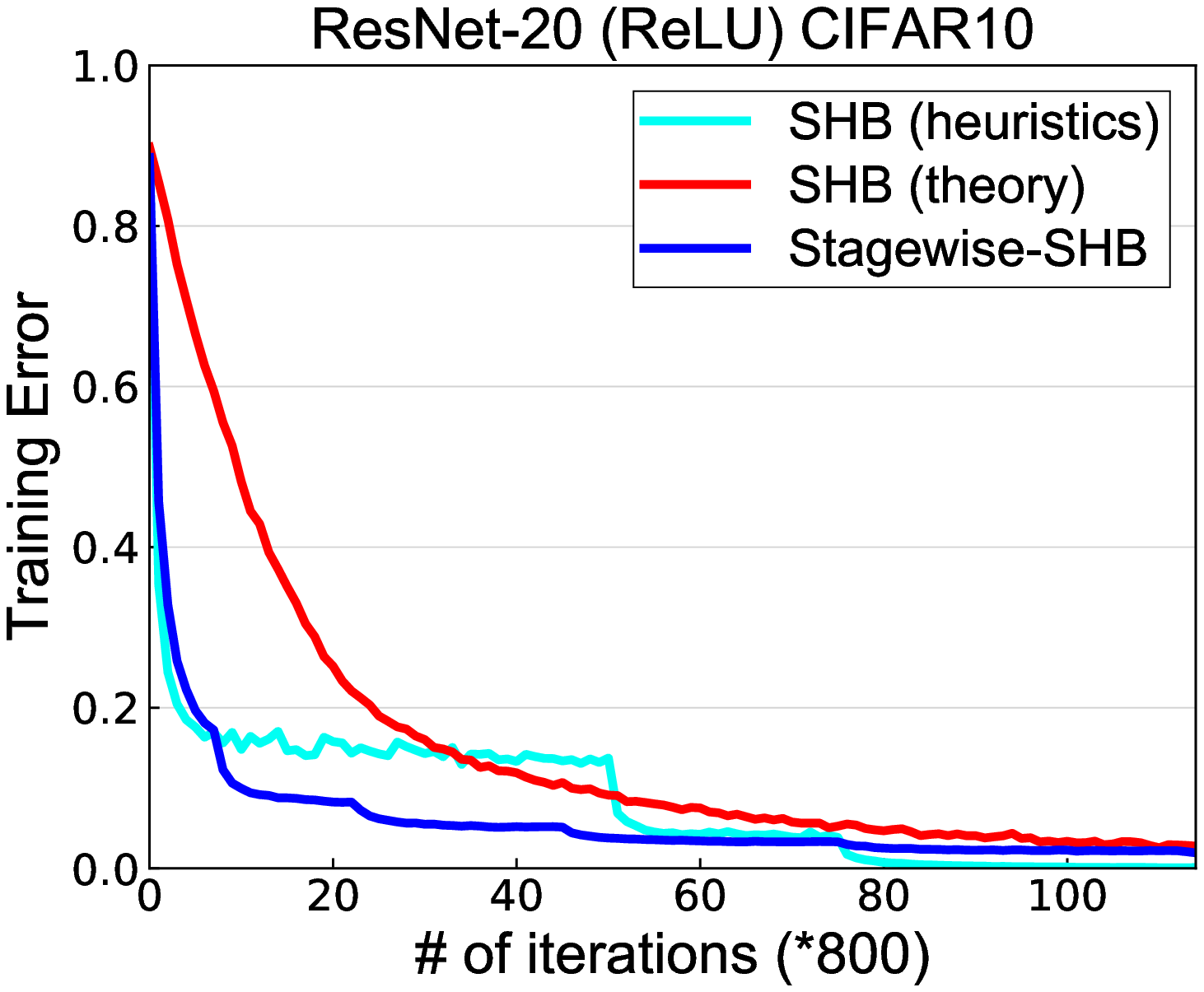}
	    \hspace*{-0.1in} 
	\includegraphics[scale=0.22]{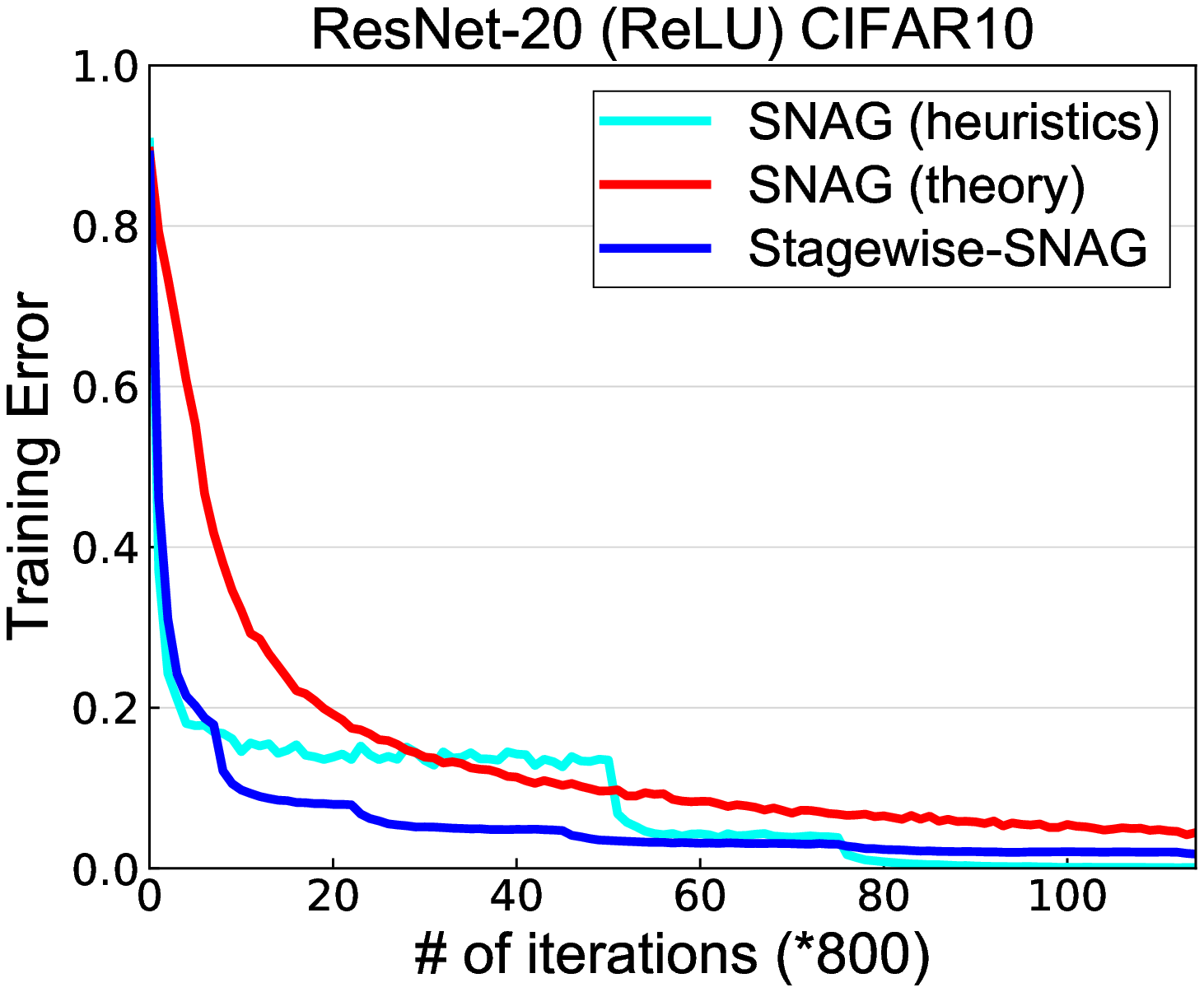}
	
\hspace*{-0.1in}	
	\includegraphics[scale=0.22]{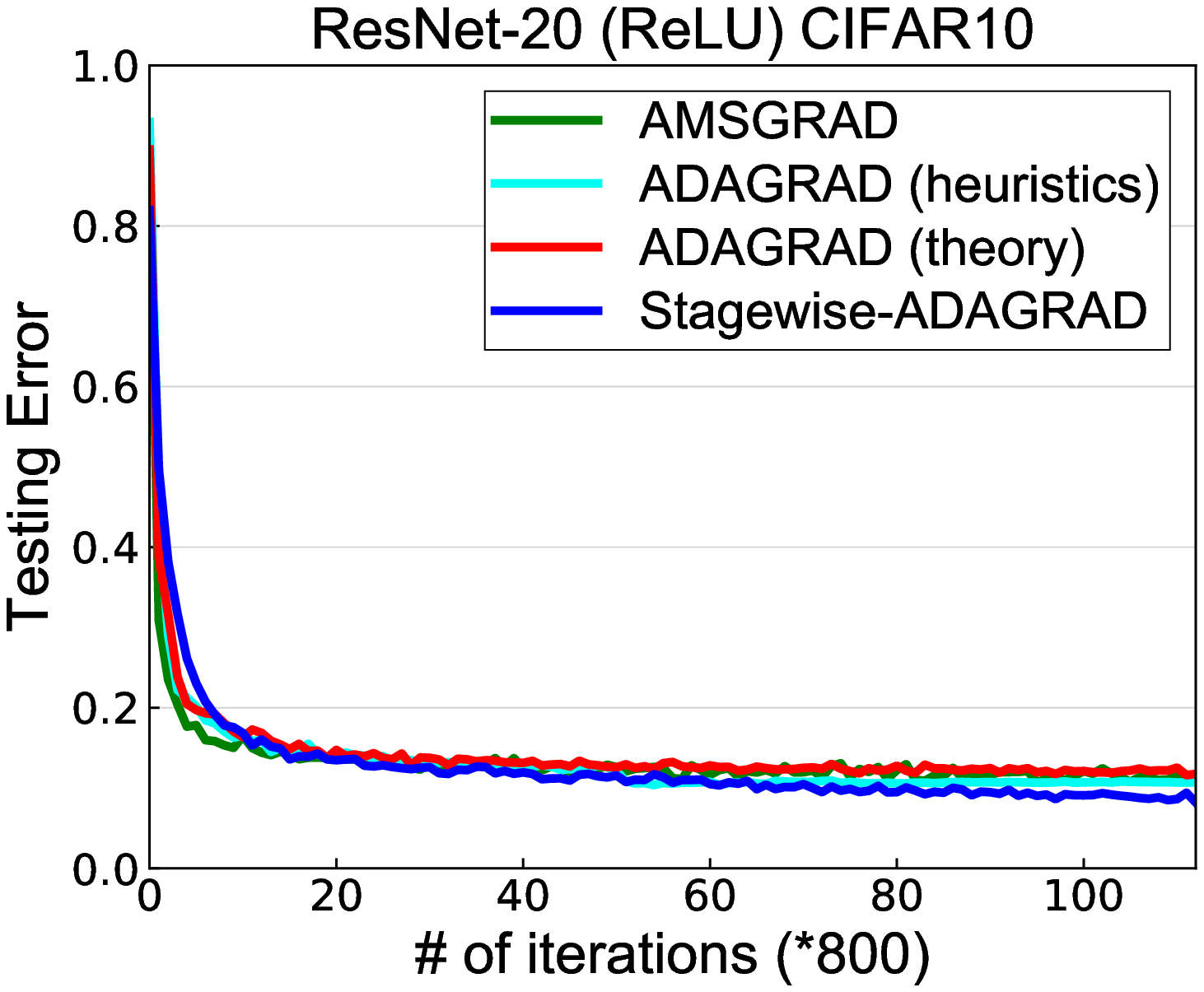}
    \hspace*{-0.1in} 
	\includegraphics[scale=0.22]{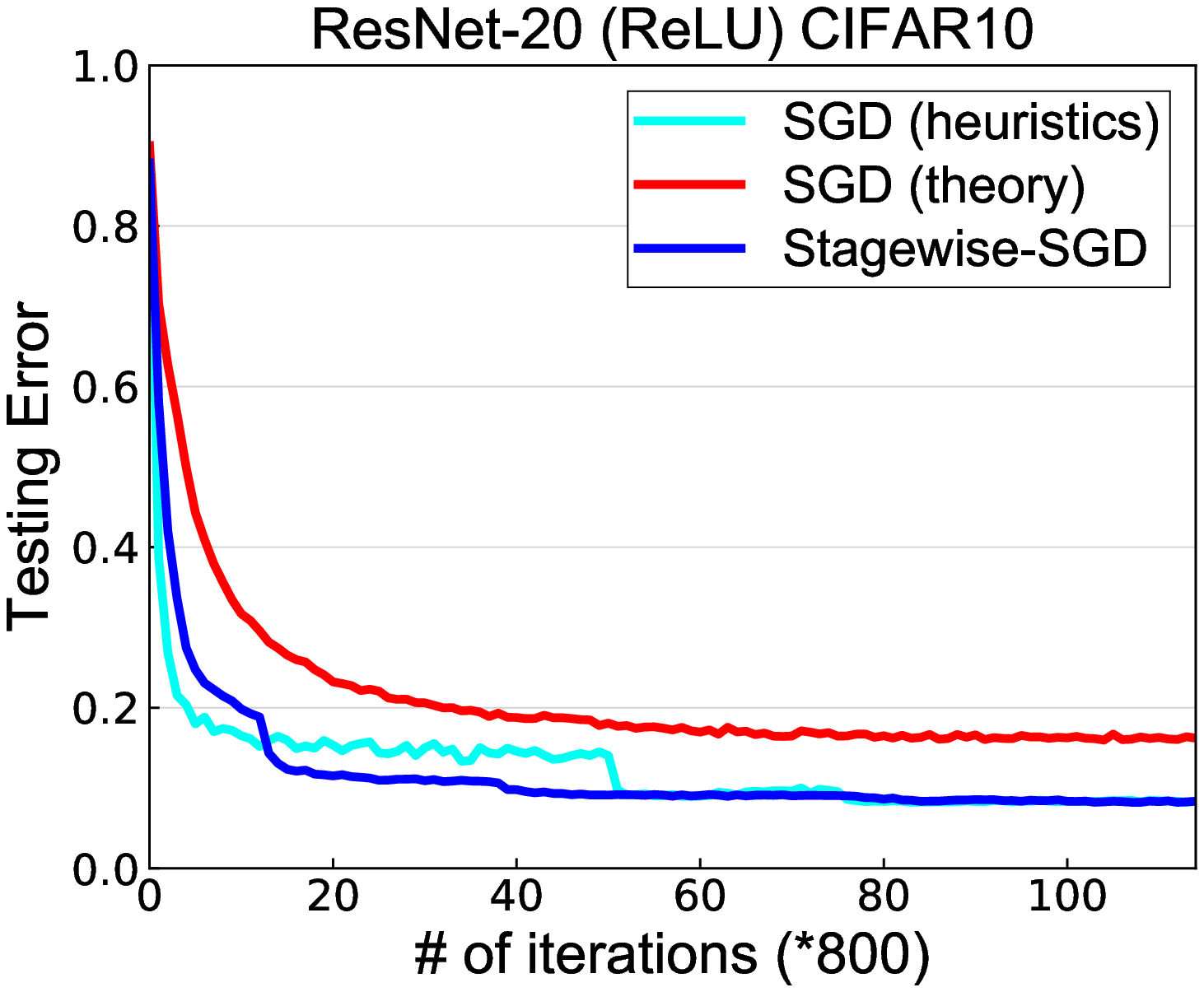}
    \hspace*{-0.1in} 
	\includegraphics[scale=0.22]{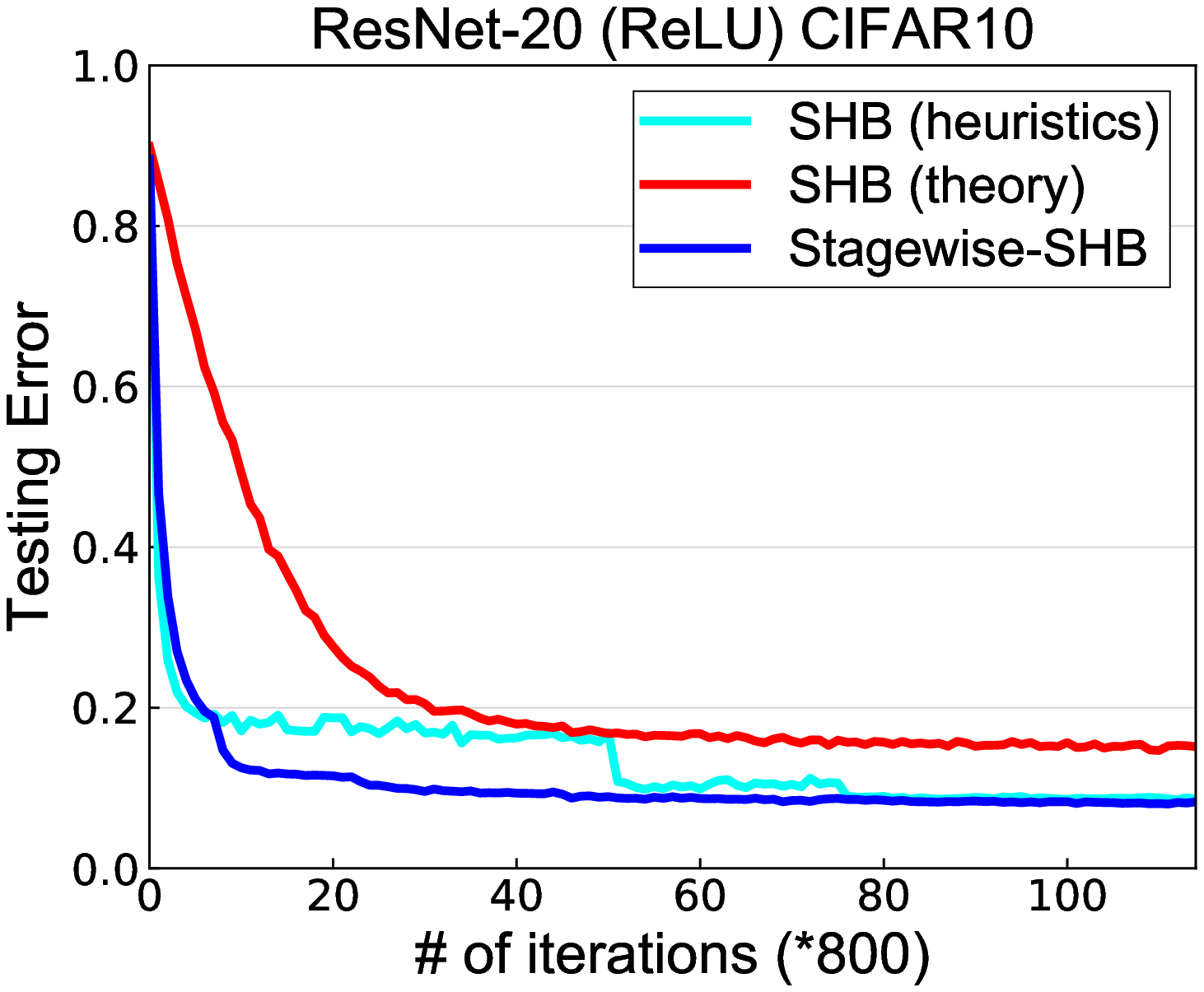}
	    \hspace*{-0.1in} 
	\includegraphics[scale=0.22]{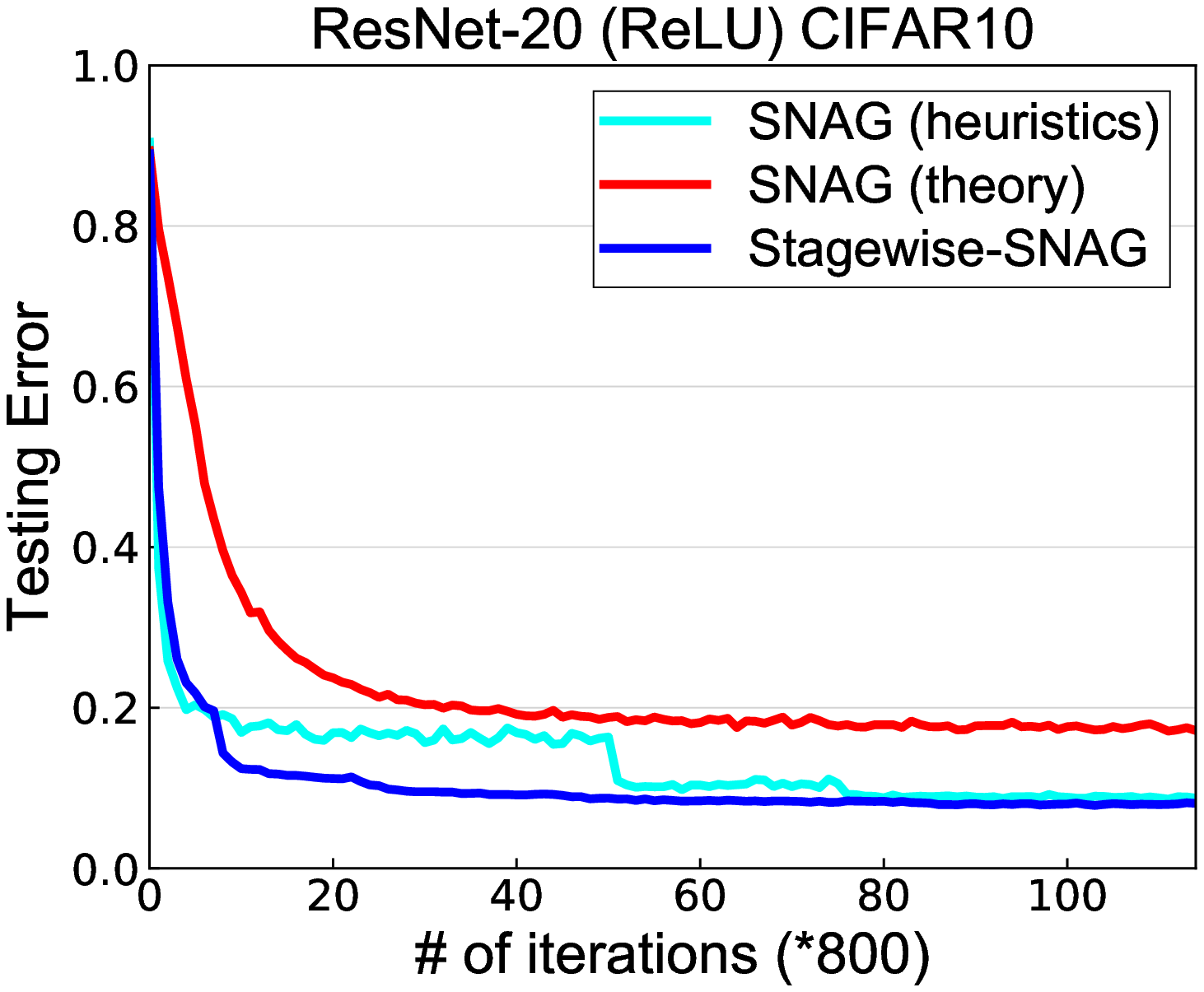}
    \vspace*{-0.15in}
\caption{Comparison of Training Error (Top) and Testing Error (bottom) on CIFAR-10 with Regularization. The regularization parameter is set $5e-4$. }

\label{fig:4}
\vspace*{0.15in}
	\hspace*{-0.1in} 
	\includegraphics[scale=0.22]{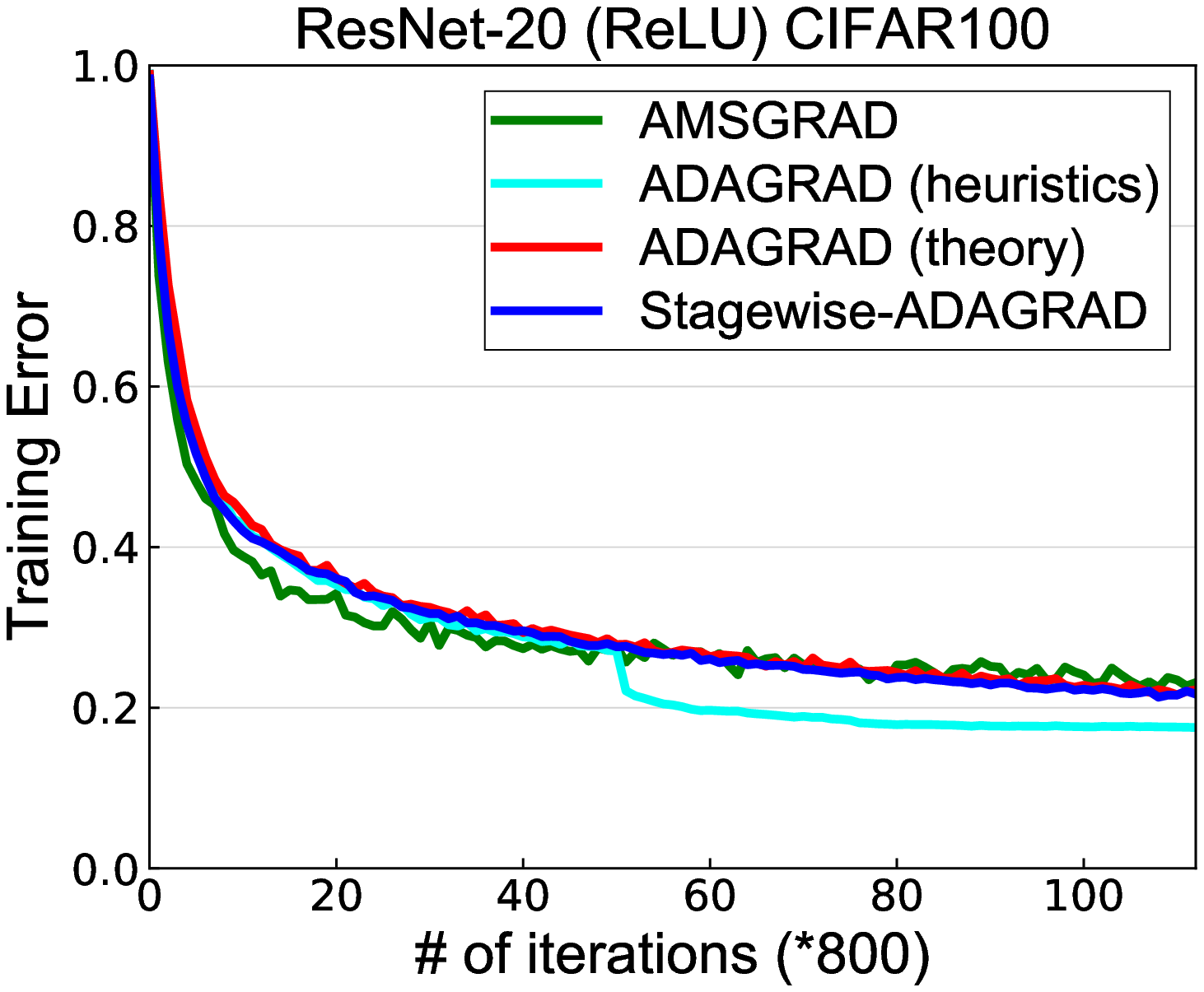}
    \hspace*{-0.1in} 
	\includegraphics[scale=0.22]{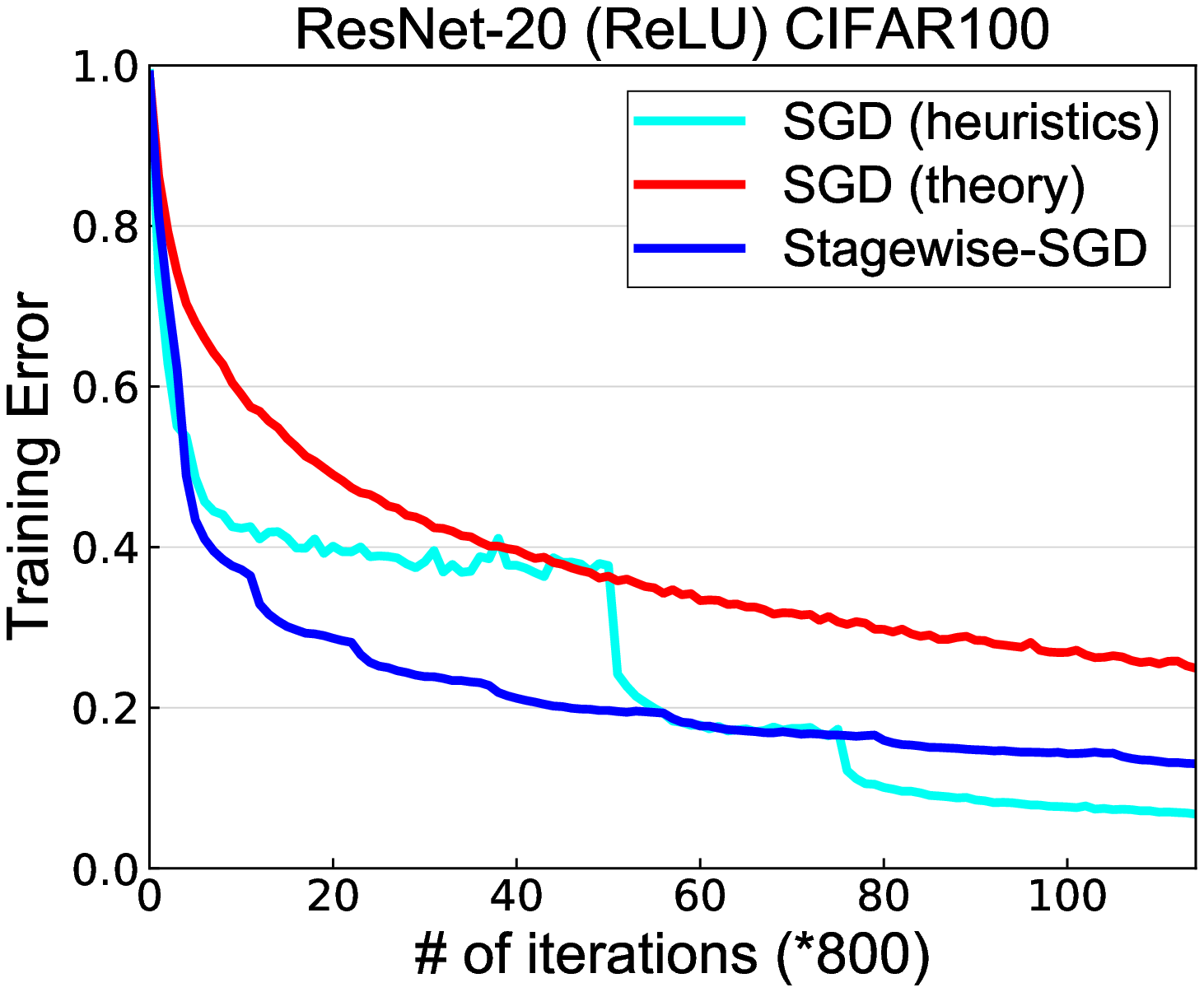}
    \hspace*{-0.1in} 
	\includegraphics[scale=0.22]{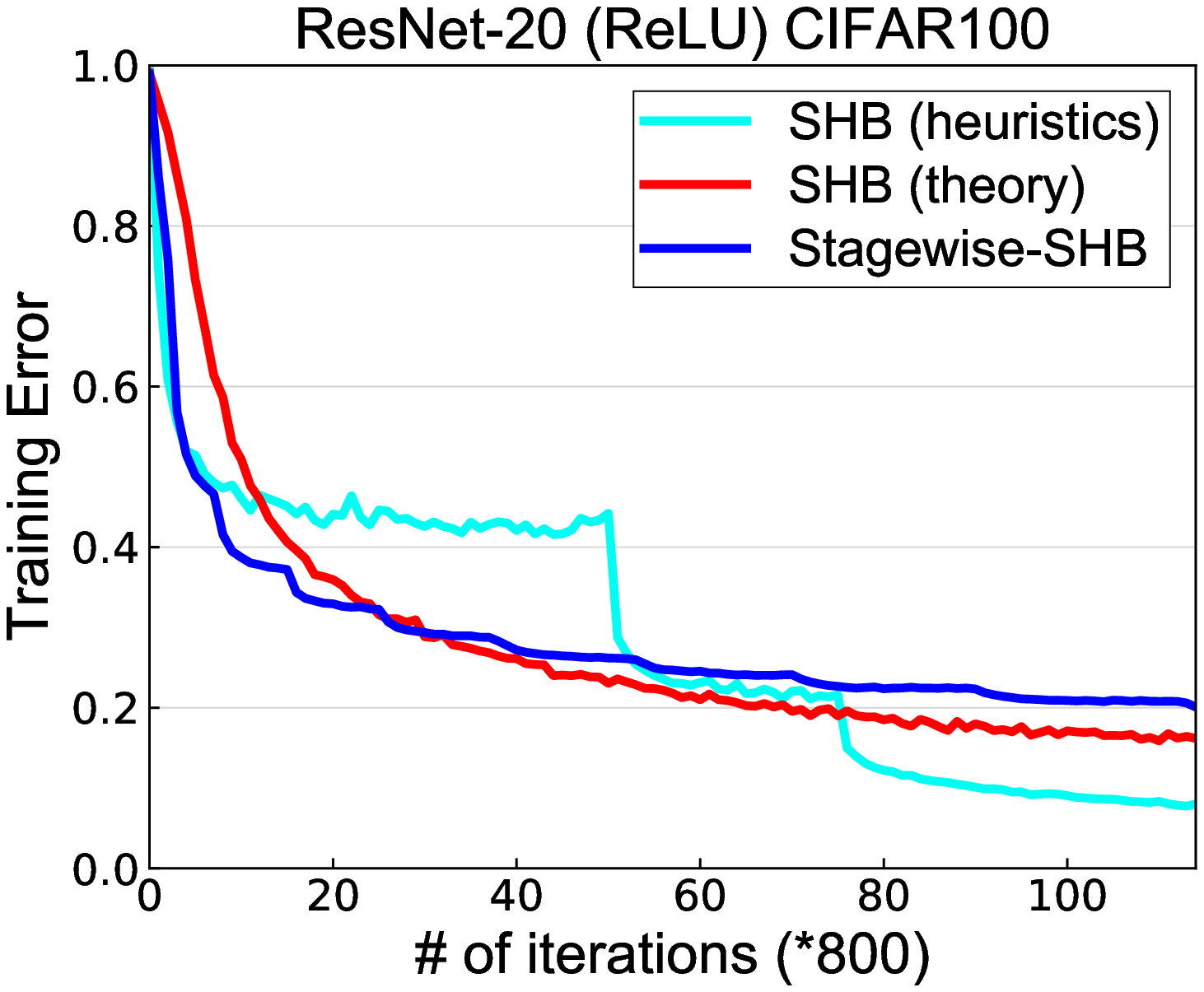}
	    \hspace*{-0.1in} 
	\includegraphics[scale=0.22]{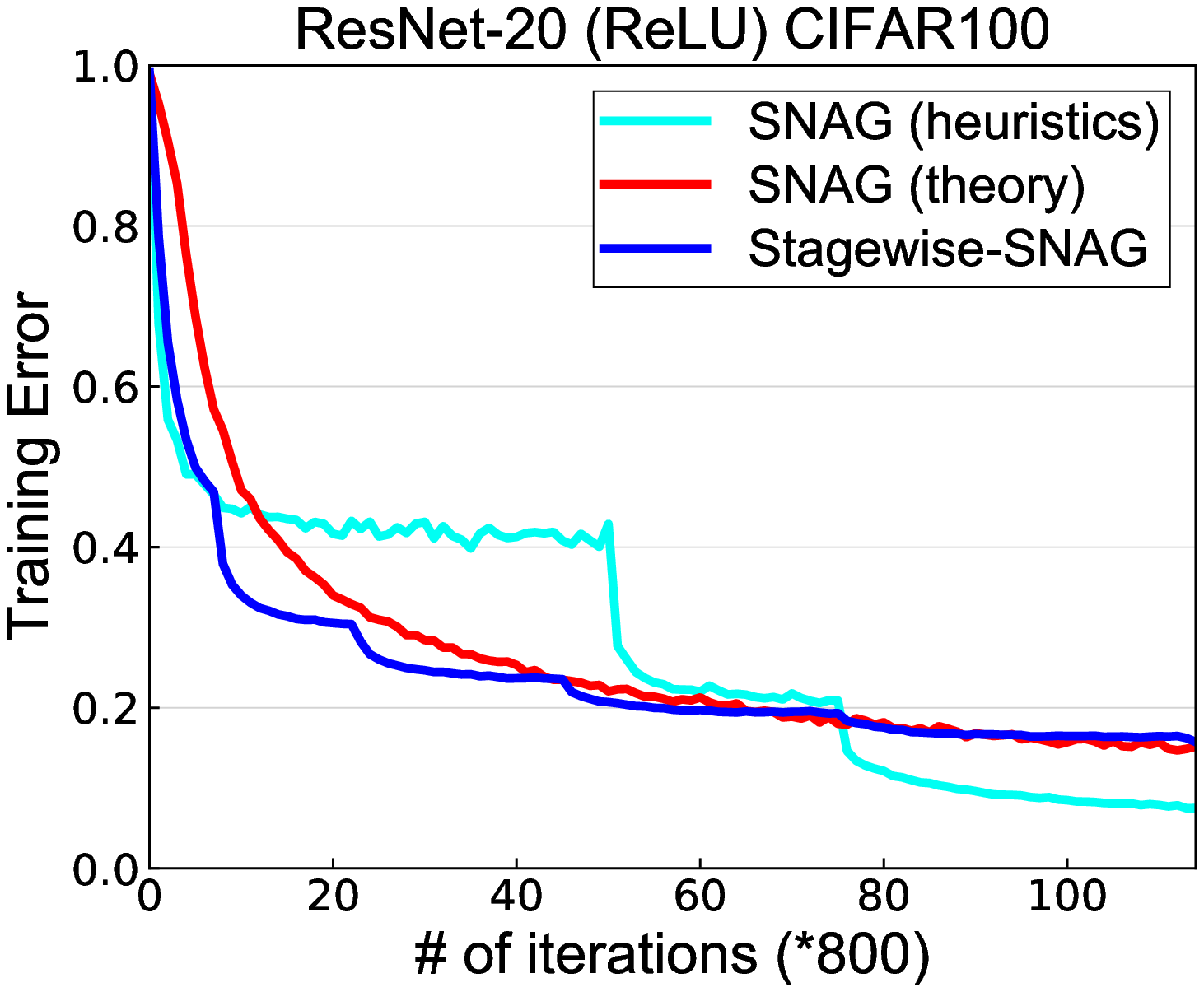}

	\includegraphics[scale=0.22]{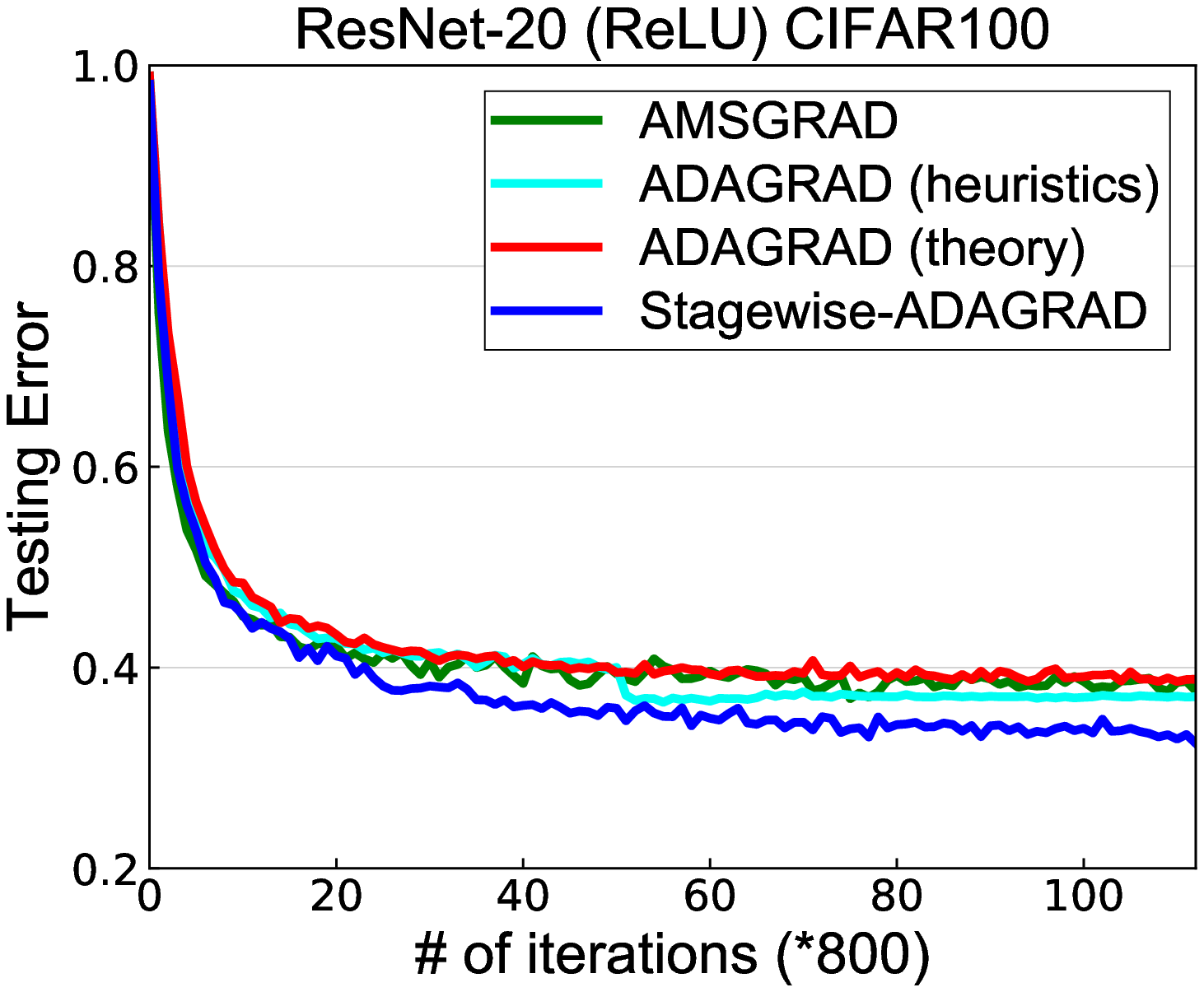}
    \hspace*{-0.1in} 
	\includegraphics[scale=0.22]{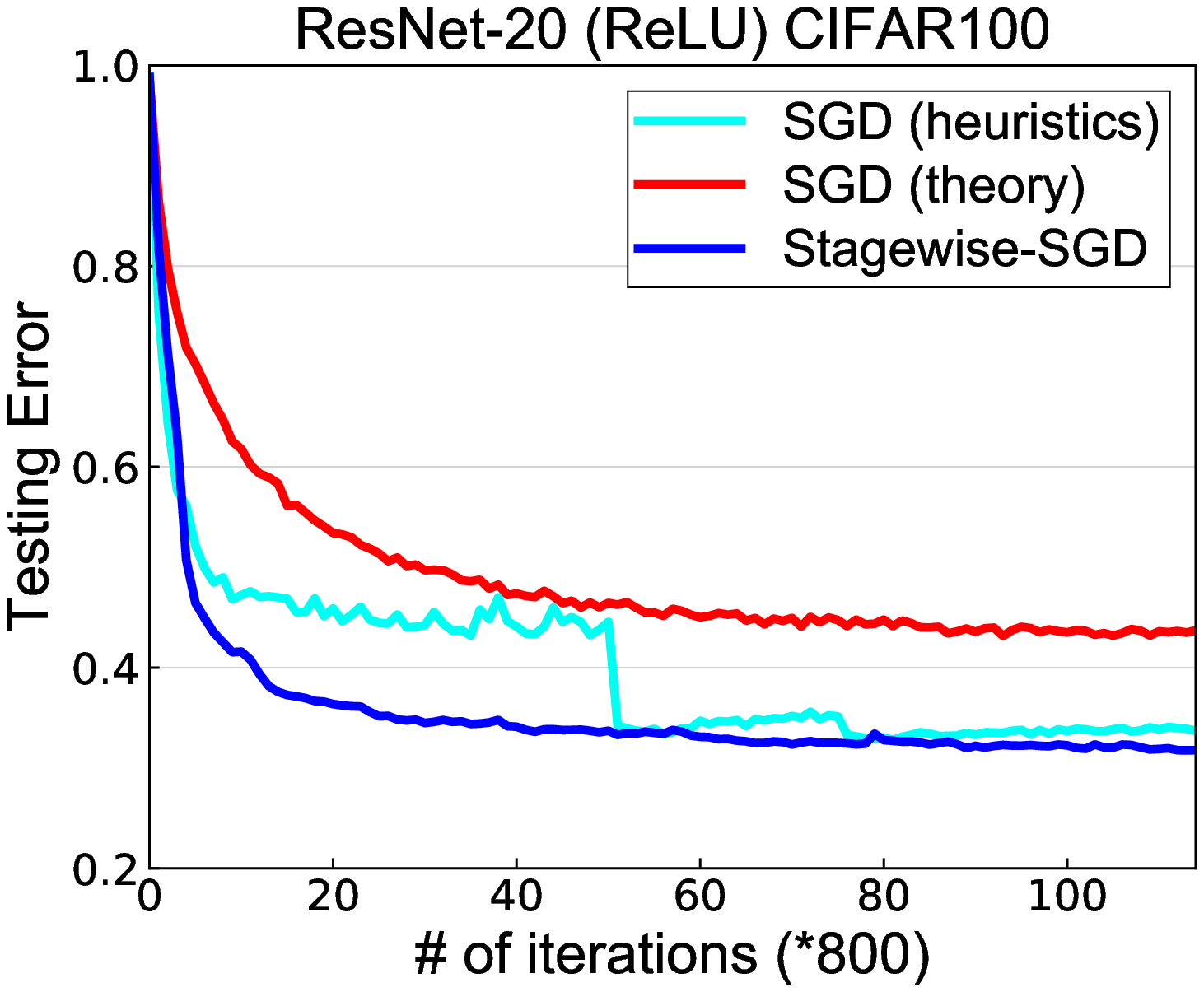}
    \hspace*{-0.1in} 
	\includegraphics[scale=0.22]{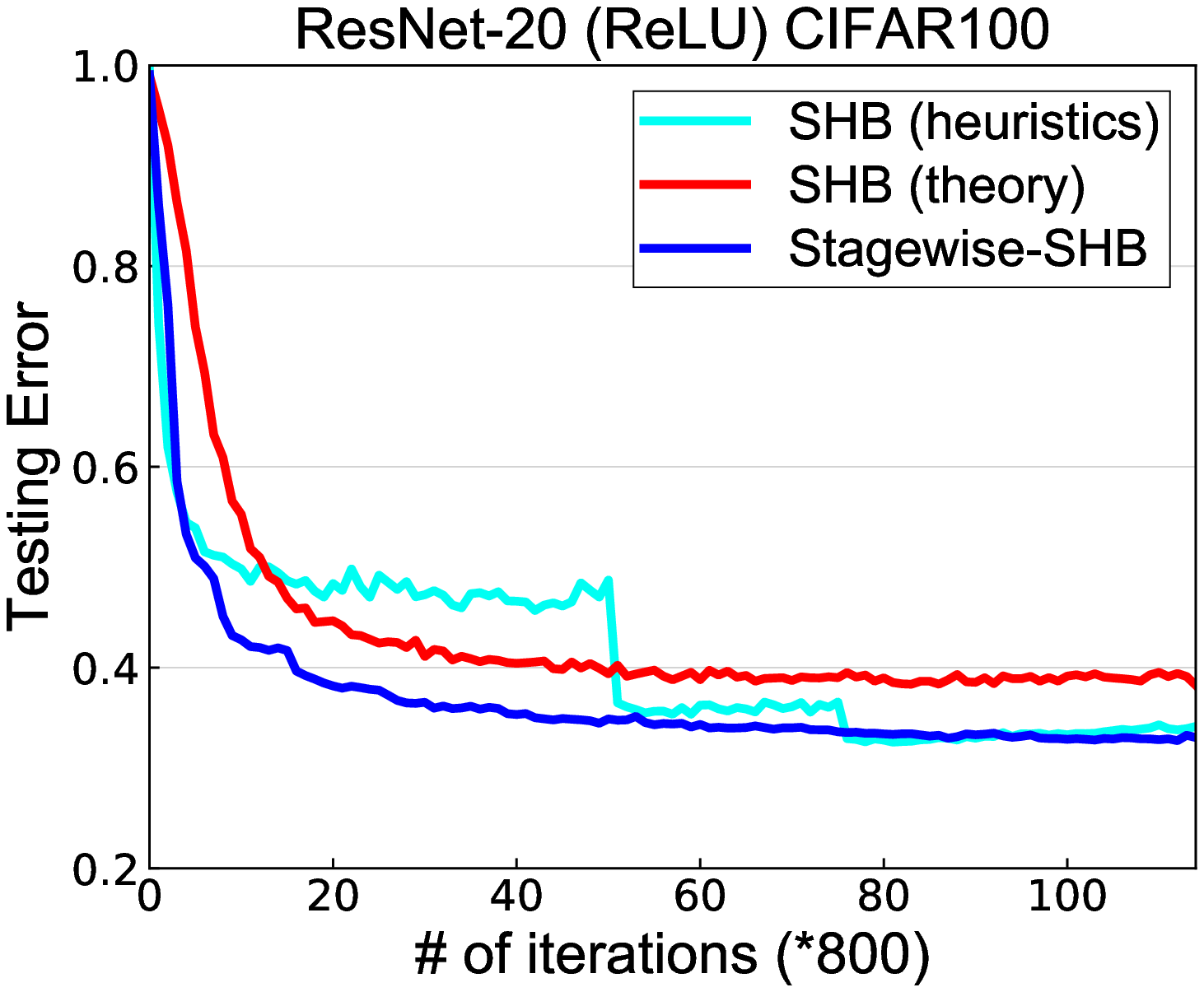}
	    \hspace*{-0.1in} 
	\includegraphics[scale=0.22]{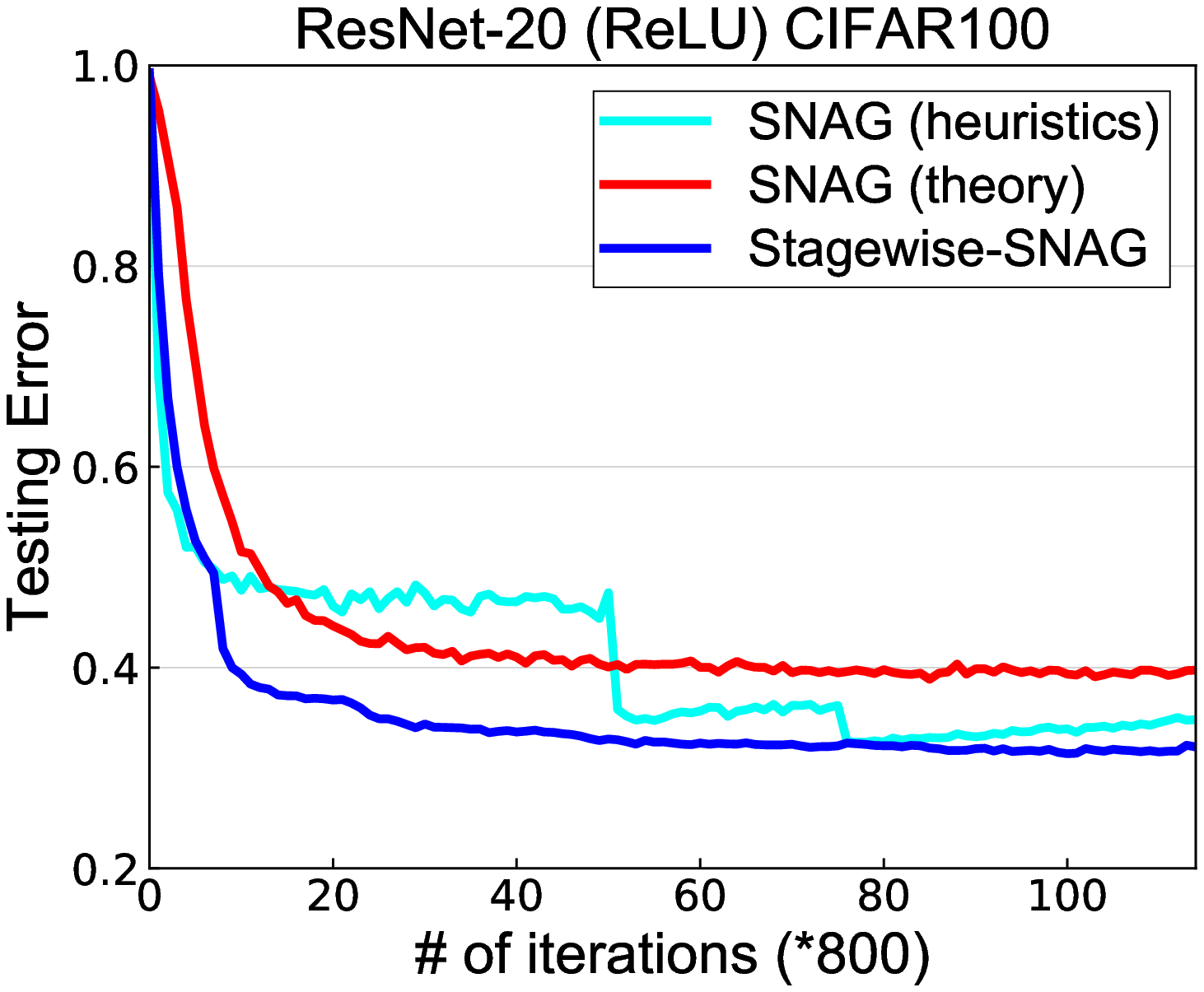}
    \vspace*{-0.15in}
\caption{Comparison of Training Error (Top) and Testing Error (bottom) on CIFAR-100 with Regularization.  The regularization parameter is set $5e-4$. }
\label{fig:5}
\vspace*{0.15in}
\end{figure}

{\bf Baselines.} We compare the proposed stagewise algorithms with  their variants implemented in TensorFlow. 
 It is notable that {\ada} has a step size (aka learning rate) parameter~\footnote{note it is not equivalent to the step size in {\sgd}.}, which is a constant in theory~\cite{ada18orabona,adamtypenoncx18,adagradmom18}. However, in the deep learning community a heuristic fixed frequency decay scheme for the step size parameter is commonly adopted~\cite{ren2018learning,DBLP:conf/nips/WilsonRSSR17}.  We thus compare two implementations of {\ada} - one with a constant learning rate parameter and another one with a fixed frequency decay scheme, which are referred to as {\ada} (theory) and {\ada} (heuristic).  For each baseline algorithms of {\sgd}, {\shb}, {\snag}, we also implement two versions - a theory version with iteratively decreasing size $\eta_0/\sqrt{t}$ suggested by previous theories and a heuristic approach with fixed frequency decay scheme used in practice, using (theory) and (heuristic) to indicate them.  The fixed frequency decay scheme used in the heuristic variants is similar to that in~\cite{he2016deep}, i.e.,  the step size parameter is decreased by 10 at 40k, 60k iterations.  We also compare stagwise {\ada} with AMSGrad~\cite{sashank2018adam} - a corrected version of Adam. 
 
 {\bf Parameters.} The stagewise step size $\eta_s =\eta_0/\sqrt{s}$ is used in stagwise {\ada},  the number of iterations $T_s$ in stagwise {\ada} is set according to Theorem~\ref{thm:nadagrad} with some simplifications for dealing with unknown $\hat G$, in particular we set $T_s$ to the smallest value larger than $T_0\sqrt{s\max_i\|g^s_{1:T_s,i}\|\sum_i\|g^s_{1:T_s,i}\|}$. For stagewise {\sgd}, {\shb}, {\snag}, the stagewise step size and iteration number is set to $\eta_s =\eta_0/s$ and $T_s=T_0s$, respectively. For parameter tuning, the initial step sizes of all algorithms are tuned in $\{0.1, 0.3, 0.5, 0.7, 0.9\}$. The value of  $\gamma$ of stagewise algorithms is tuned in $\{1, 10, 100, 500, 1000, 1500, 2000, 3000\}$. The initial value $T_0$ for stagewise {\sgd}, {\shb}, {\snag} is tuned in $\{10, 100, 1\text{k}, 5\text{k}, 6\text{k}, 7\text{k}, 10\text{k}, 20\text{k}\}$, and that for stagewise {\ada} is tuned in $\{1, 10, 15, 20, 25, 50, 100\}$. 
 
 

\begin{figure}[t] 
\centering
	\hspace*{-0.2in} 
		{\includegraphics[scale=0.22]{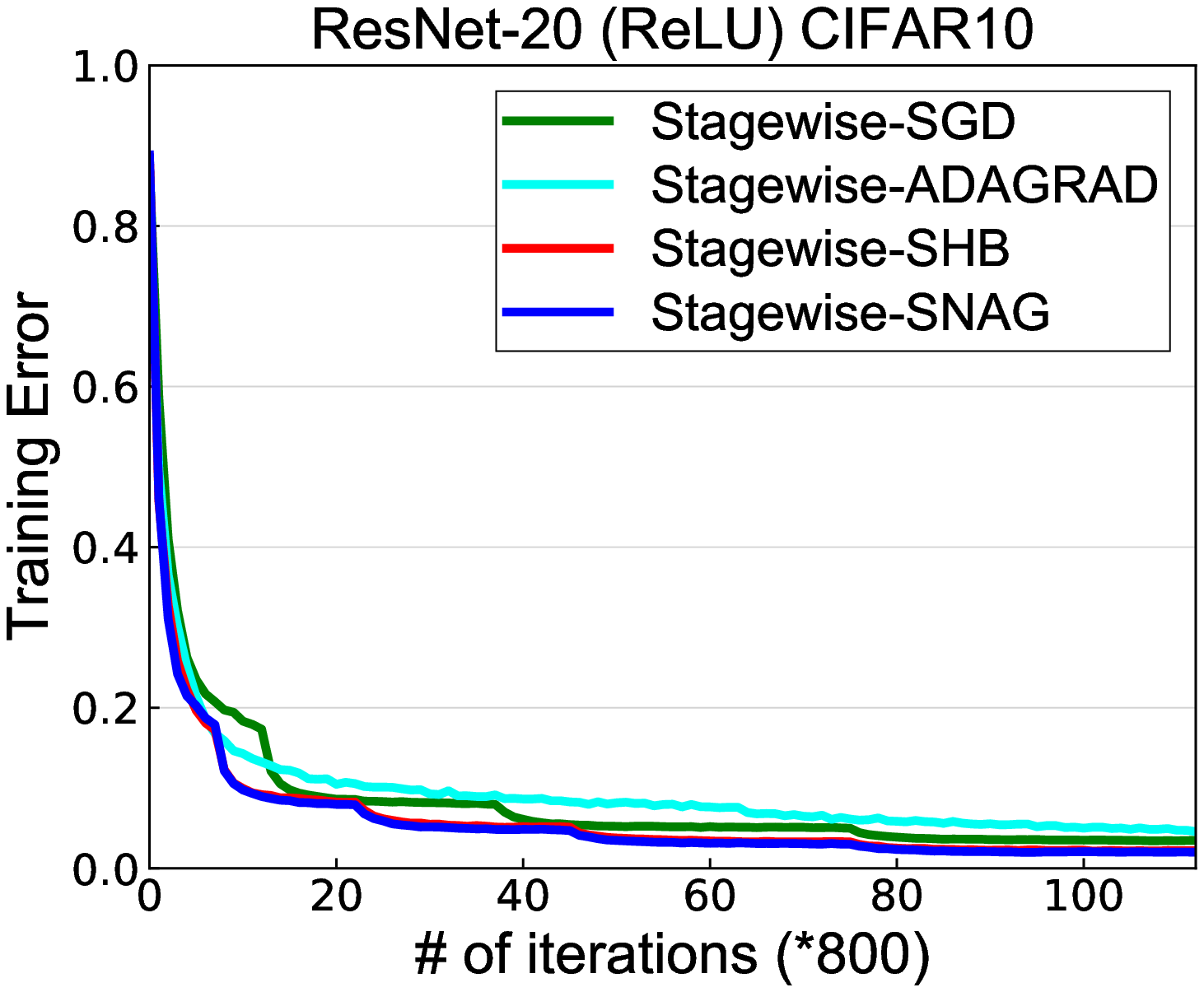}}
    \hspace*{-0.1in} 
	{\includegraphics[scale=0.22]{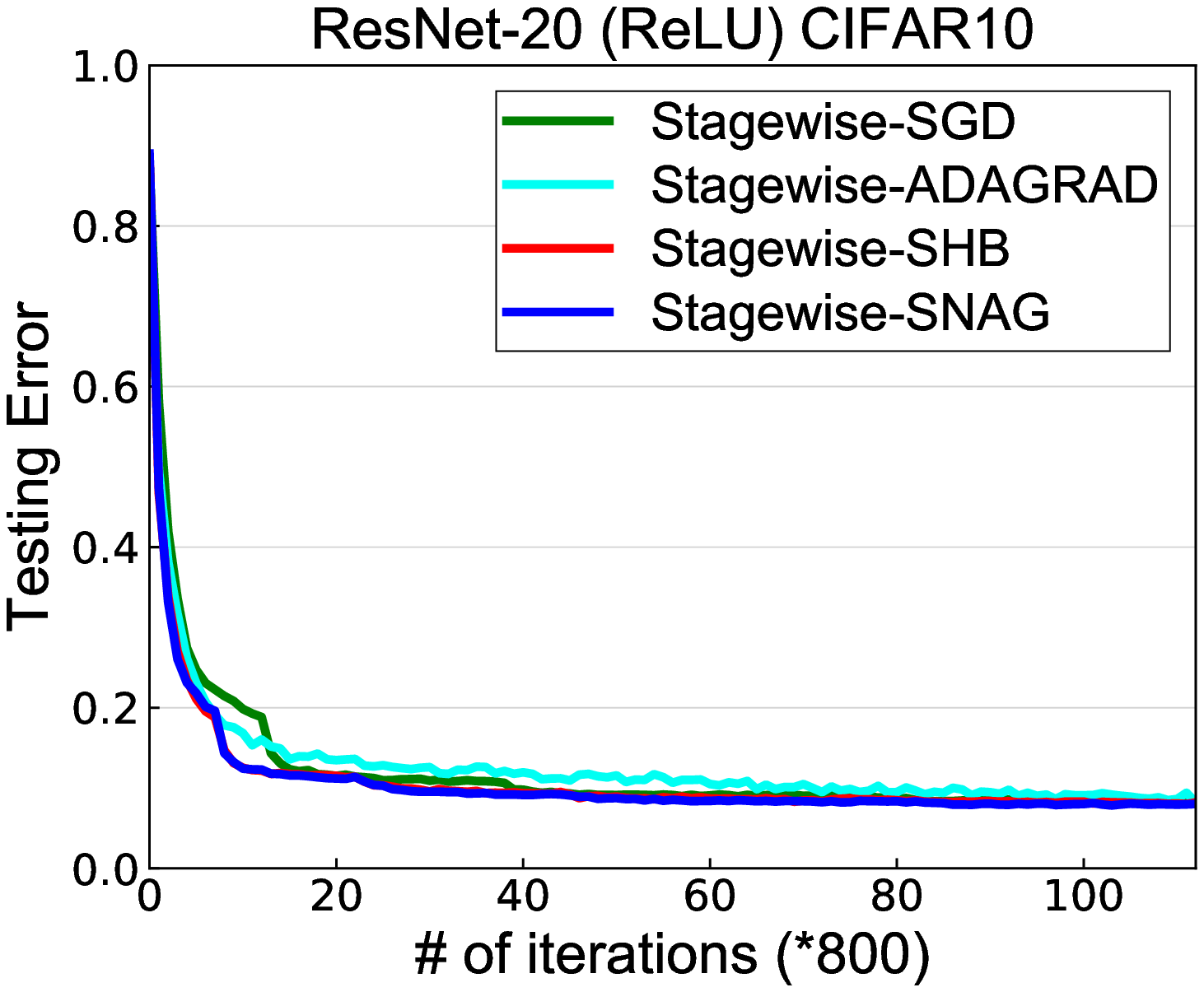}}
    \hspace*{-0.1in} 
	{\includegraphics[scale=0.22]{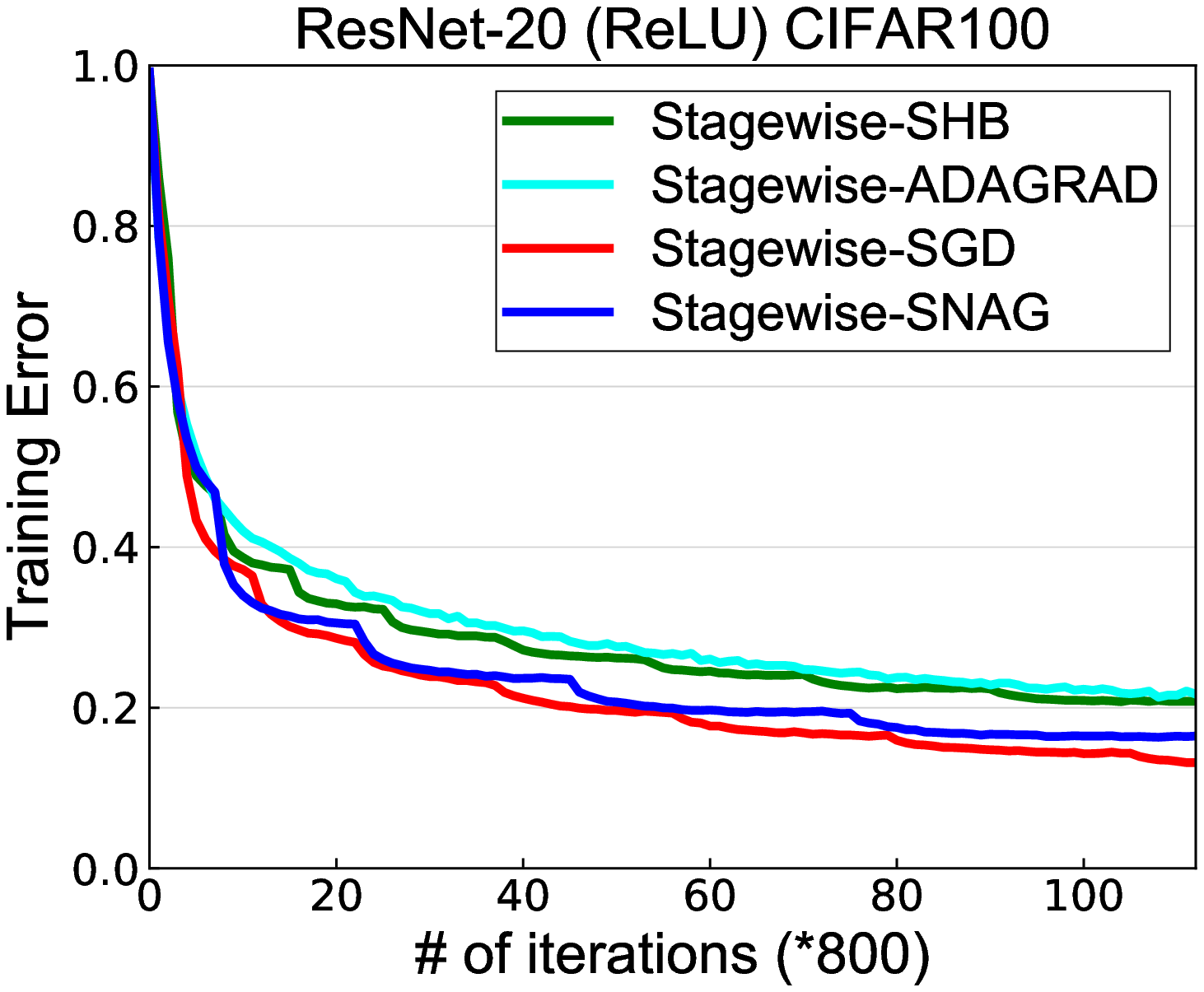}}
	    \hspace*{-0.1in} 
	{\includegraphics[scale=0.22]{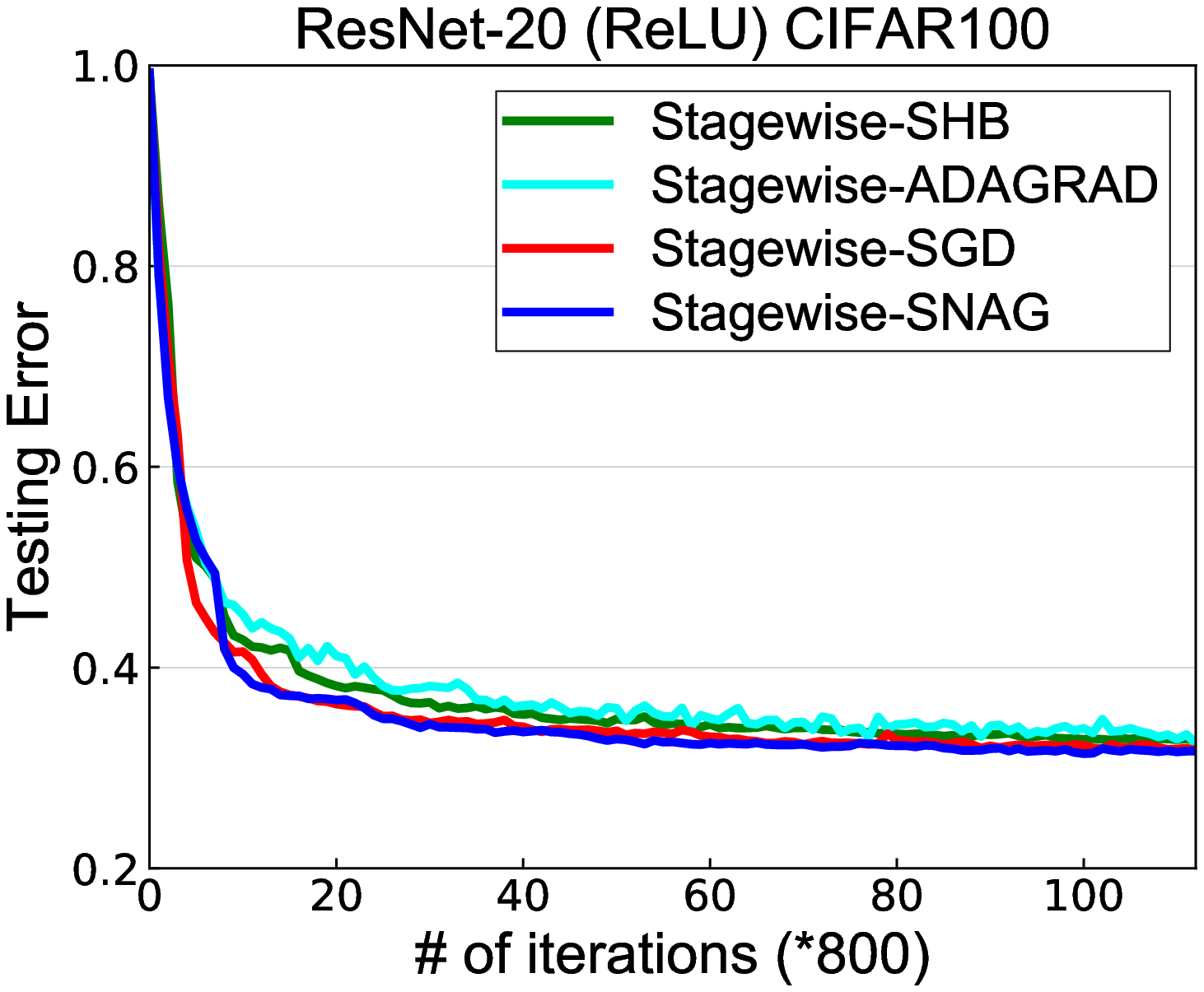}}
	
	\hspace*{-0.1in} 
    {\includegraphics[scale=0.22]{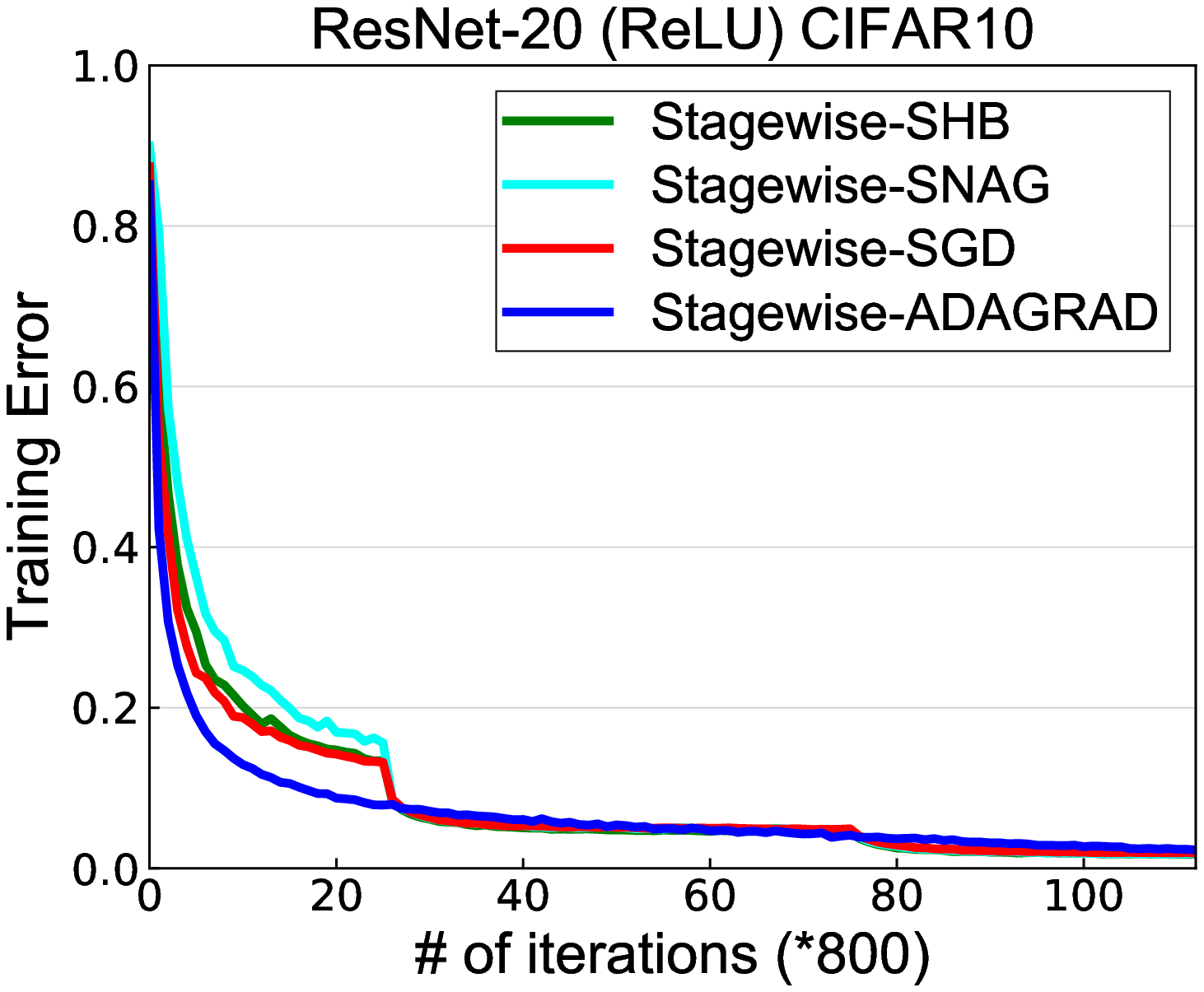}}
    \hspace*{-0.1in} 
	{\includegraphics[scale=0.22]{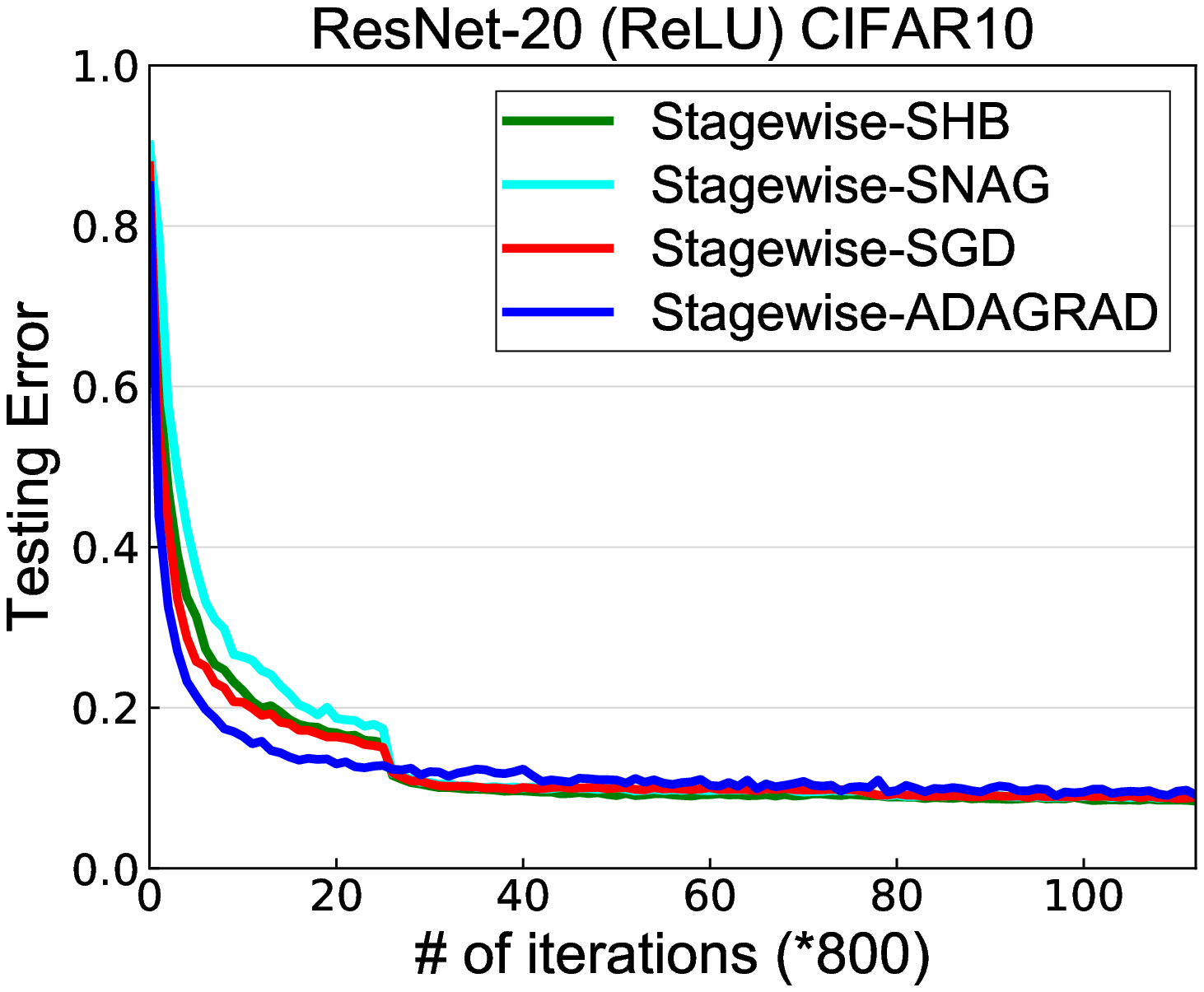}}
    \hspace*{-0.1in} 
	{\includegraphics[scale=0.22]{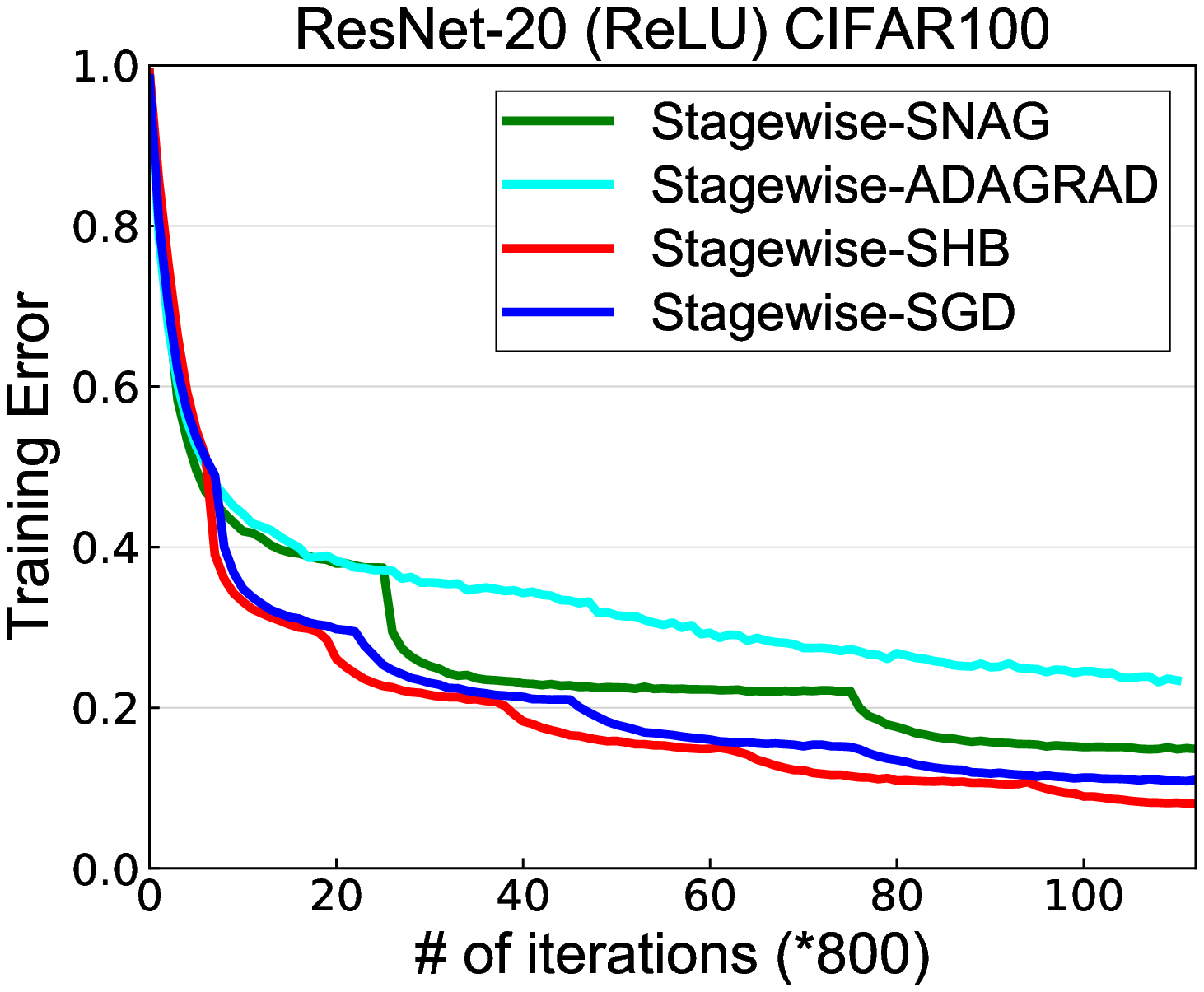}}
	    \hspace*{-0.1in} 
	{\includegraphics[scale=0.22]{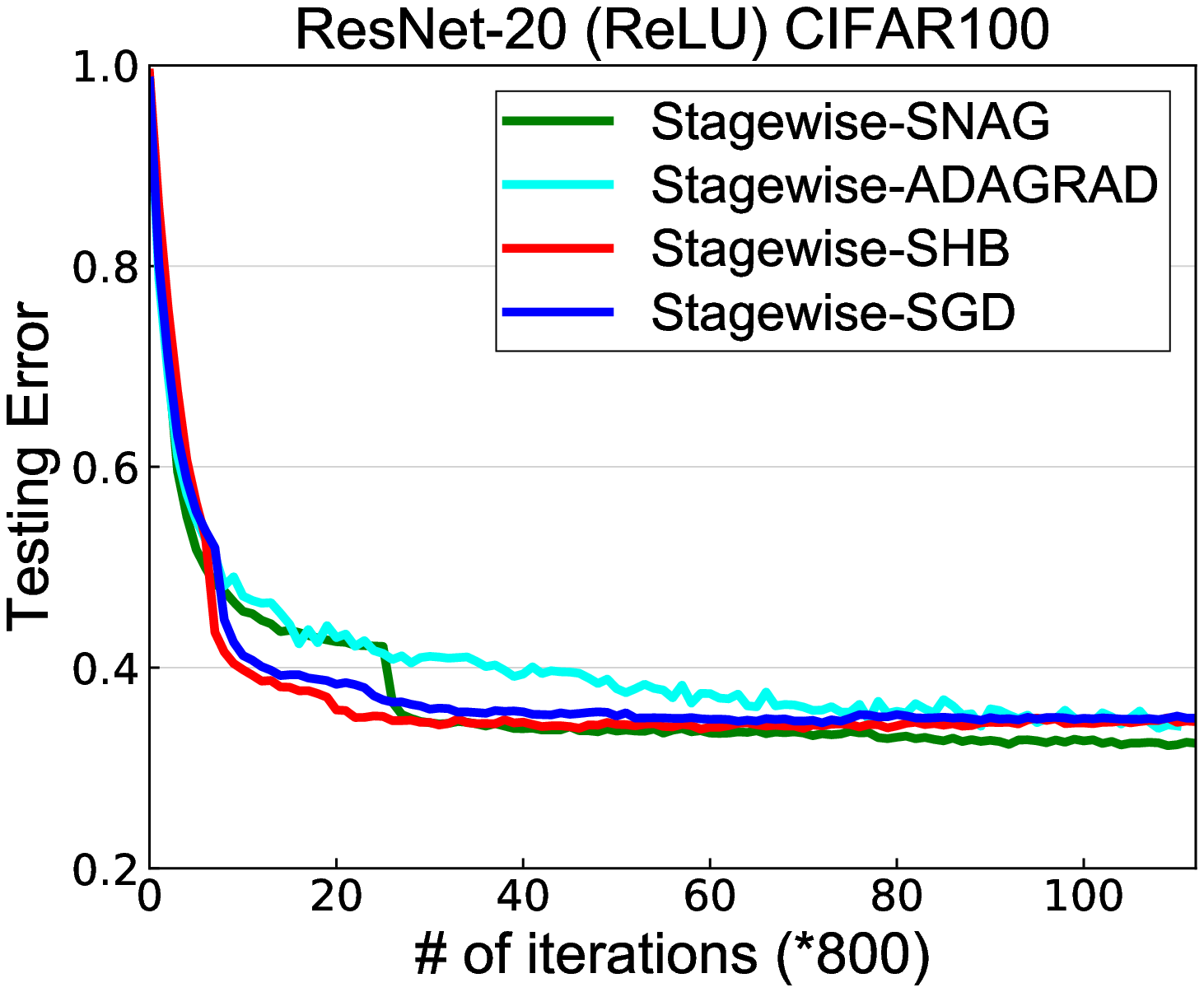}}
    \vspace*{-0.15in}
\caption{Comparison of different stagewise algorithms with regularization (Top) and without regularization (bottom).  The regularization parameter is set $5e-4$.} 
\label{fig:6}
\vspace*{0.15in}
\end{figure}

{\bf Results.} We consider two settings - with/without an $\ell_2$ norm regularization on weights.  For comparison, we evaluate the training error and testing error of obtained solutions in the process of training. For our stagewise algorithms, the evaluation is done based on the current averaged solution, and for other baselines the evaluation is done based on the current solution. The comparisons of training and testing error in the two settings (w/o regularization) on the two datasets are plotted in Figure~\ref{fig:2}, \ref{fig:3}, \ref{fig:4}, \ref{fig:5}. We also compare four stagewise methods as shown in Figure ~\ref{fig:6}. The final testing error (after running $80$k iterations) of different algorithms are reported in Table~\ref{tab:1}. From all results, we have several observations. (i) The proposed stagewise algorithms perform much  better in terms of testing error than the existing theoretical versions reported in the literature (marked with theory in the legend). This indicates the proposed stagewise step size scheme is better than iteratively decreasing step size scheme. (ii) The proposed stagewise algorithms achieve similar and sometimes even better testing error than the heuristic approaches with a fixed frequency decay scheme used in practice. However,  the heuristic approaches usually give smaller training error. This seems to indicate that the proposed algorithms are less vulnerable to the overfitting. In another word, the proposed algorithms have smaller generalization error, i.e., the difference between the testing error and the training error. (iii) The proposed stagewise algorithms (stagewise-{\sgd}, {\shb}, {\snag}, {\ada}) have comparable result and there is no clear winner depending on datasets and on whether regularization is added.

\begin{table*}[t]
		\caption{Comparison of Final Testing Error (\%) on CIFAR-10 and CIFAR-100 Datasets}
		\centering
		\label{tab:1}
		\scalebox{0.90}{\begin{tabular}{l|l|l|l|l}
			\toprule
			& \multicolumn{2}{c|}{CIFAR-10} &\multicolumn{2}{c}{CIFAR-100}\\
			\midrule
			Algorithms & with reg. &without reg. &with reg. &without reg.  \\
		\midrule
			SGD (theory)&16.25&19.18&43.51&45.78 \\
			SGD (heuristic)&8.34&10.81&33.67&37.19\\
			Stagewise-SGD & {\bf 8.34}&{\bf 9.01}&{\bf 32.25}&{\bf 34.95}\\
		    \midrule
			SHB (theory)&15.67&16.55&39.15&46.23\\
			SHB (heuristic)&8.58&10.28&33.30&37.56\\
			Stagewise-SHB &{\bf 8.30}&{\bf 8.61}&{\bf 32.85}&{\bf 34.49}\\
			\midrule
			SNAG (theory)&17.64&16.76&39.34&44.21\\
			SNAG (heuristic)&8.85&10.34&33.89&36.84\\
			Stagewise-SNAG & {\bf 8.00}&{\bf 8.93}&{\bf 31.42}&{\bf 33.29}\\
			\midrule
			AMSGrad& 10.76 & 11.13&38.62&39.96\\
			AdaGrad (theory)&12.11&13.96&39.09&44.49\\
			AdaGrad (heuristic)&10.71&13.80&37.04&41.06\\
			Stagewise-AdaGrad & {\bf 9.09}&{\bf 9.51}&{\bf 33.95}&{\bf 34.62}\\
				\bottomrule
		\end{tabular}}
	\end{table*}


\section{Conclusion}
In this paper, we have proposed a universal stagewise learning framework for solving stochastic non-convex optimization problems, which employs well-known tricks in practice that have not been well analyzed theoretically. We provided theoretical convergence results for the proposed algorithms for non-smooth non-convex optimization problems. We also established an adaptive convergence of a stochastic algorithm  using data adaptive coordinate-wise step size of {\ada}, and exhibited its faster convergence than non-adaptive stepsize when the growth of cumulative stochastic gradients is slow similar to that in the convex case. For future work, one may consider developing more variants of the proposed meta algorithm, e.g., stagewise AMSGrad, stagewise RMSProp, etc.


\section*{Acknowledgement} T. Yang are partially supported by National Science Foundation (IIS-1545995). Part of this work was done when Chen is interning at JD AI Research and Yang is visiting JD AI Research.

\appendix
\section{Proof of Theorem~\ref{thm:sgd}}
\begin{proof} Below, we use $\E_s$ to denote expectation over randomness in the $s$-th stage given all history before $s$-th stage. 
 Define
 \begin{align}\label{eqn:z}
 \z_{s} = \arg\min_{\x\in\Omega} f_s(\x) = \text{prox}_{\gamma (\phi+\delta_\Omega)}(\x_{s-1})
 \end{align}
Then $\nabla\phi_\gamma(\x_{s-1})  =\gamma^{-1}(\x_{s-1} - \z_s)$. Then we have $\phi(\x_{s})\geq \phi(\z_{s+1})+\frac{1}{2\gamma}\|\x_{s}-\z_{s+1}\|^2$. Next, we apply Lemma~\ref{lem:sgd} to each call of {\sgd} in stagewise {\sgd}, 
\begin{align*}
\E[f_s(\x_s) - f_s(\z_s)]\leq  \underbrace{\frac{\|\z_s - \x_{s-1}\|^2}{2\eta_s T_s} + \frac{\eta_s \hat G^2}{2}}\limits_{\mathcal E_s},
\end{align*}
where $\hat G^2$ is the upper bound of $\E[\|\g(\x; \xi) + \gamma^{-1}(\x - \x_{s-1})\|^2]$, which exists and can be set to $2G^2 + 2\gamma^{-2}D^2$ due to the Assumption~\ref{ass:0}-\ref{ass:1} and the bounded assumption of the domain. Then 
\begin{align*}
    \E_s\bigg[\phi(\x_{s}) + \frac{1}{2\gamma}\|\x_{s}-\x_{s-1}\|^2 \bigg] &\leq f_{s}(\z_{s}) + \mathcal{E}_{s}\leq f_{s}(\x_{s-1}) + \mathcal{E}_{s}\\
    &\leq \phi(\x_{s-1})+\mathcal{E}_{s}
\end{align*}
On the other hand, we have that 
\begin{align*}
    \|\x_{s}-\x_{s-1}\|^2 =& \|\x_{s}-\z_{s}+\z_{s}-\x_{s-1}\|^2\\
    =& \|\x_{s}-\z_{s}\|^2+\|\z_{s}-\x_{s-1}\|^2 + 2\langle \x_{s}-\z_{s}, \z_{s}- \x_{s-1}\rangle\\
    \geq& (1-\alpha_{s}^{-1})\|\x_{s}-\z_{s}\|^2 + (1-\alpha_{s})\|\x_{s-1}-\z_{s}\|^2
\end{align*}
where the inequality follows from the Young's inequality with $0<\alpha_{s}<1$. Thus we have that
\begin{align}\label{eqn:bound1}
   \E_s\bigg[ \frac{(1-\alpha_{s})}{2\gamma} \|\x_{s-1}-\z_{s}\|^2\bigg] \leq& \E_s\bigg[ \phi(\x_{s-1}) -\phi(\x_{s})+\frac{(\alpha_{s}^{-1}-1)}{2\gamma}\|\x_{s}-\z_{s}\|^2+\mathcal{E}_{s} \bigg]\nonumber\\
    \leq& \E\bigg[\phi(\x_{s-1}) -\phi(\x_{s})+\frac{(\alpha_{s}^{-1}-1)}{\gamma(\gamma^{-1}-\mu)}(f_s(\x_s) - f_s(\z_s))+ \mathcal{E}_{s} \bigg]\nonumber\\
    \leq &\E\bigg[ \phi(\x_{s-1}) -\phi(\x_{s}) + \frac{\alpha_{s}^{-1}-\gamma\mu}{(1-\gamma\mu)}\mathcal{E}_{s}\bigg]
\end{align}
Combining the above inequalities,  we have that 
\begin{align*}
    \bigg((1-\alpha_{s})\gamma - \frac{ \gamma^2(\alpha_{s}^{-1} - \mu \gamma)}{ (1-\mu\gamma)\eta_{s} {T_{s}}} \bigg)& \E_s[\|\nabla \phi_{\gamma}(\x_{s-1})\|^2]\leq \E_s\bigg[2\Delta_{s} + \frac{(\alpha_{s}^{-1} - \mu\gamma )\eta_{s}\hat G^2}{(1-\mu\gamma)}\bigg]
    \end{align*}
Multiplying both sides by $w_s$, we have that
\begin{align*}
   w_s \bigg((1-\alpha_{s})\gamma - \frac{ \gamma^2(\alpha_{s}^{-1} - \mu \gamma)}{ (1-\mu\gamma)\eta_{s} {T_{s}}} \bigg)& \E_s[\|\nabla \phi_{\gamma}(\x_{s-1})\|^2]\leq \E_s\bigg[2w_s\Delta_{s} + \frac{(\alpha_{s}^{-1} - \mu\gamma )w_s\eta_{s}\hat G^2}{(1-\mu\gamma)}\bigg]
    \end{align*}
By setting $\alpha_s = 1/2$ and $\gamma  = 1/(2\mu)$, $T_s\eta_s\geq 12\gamma$,  we have
\begin{align*}
   \frac{1}{4}w_s\gamma \E_s[\|\nabla \phi_{\gamma}(\x_{s-1})\|^2]\leq \E_s[2w_s\Delta_{s} +3w_s\eta_{s}\hat G^2]
    \end{align*}
 By summing over $s=1,\ldots,S+1$, we have
 \begin{align*}
 \sum_{s=1}^{S+1}w_s \E[\|\nabla \phi_{\gamma}(\x_{s-1})\|^2]\leq \E\bigg[16\mu\sum_{s=1}^{S+1}w_s\Delta_{s} + 24\mu\sum_{s=1}^{S+1}w_s\eta_{s}\hat G^2\bigg]
   \end{align*}
Taking the expectation w.r.t. $\tau\in \{0,\ldots, S\}$, we have that
\begin{align*}
   \E[\|\nabla\phi_{\gamma}(\x_{\tau})\|^2]] \leq \E\bigg[\frac{16\mu\sum_{s=1}^{S+1}w_s\Delta_s}{\sum_{s=1}^{S+1}w_s} + \frac{24\mu\sum_{s=1}^{S+1}w_s\eta_{s}\hat G^2}{\sum_{s=1}^{S+1}w_s}\bigg]
\end{align*}
For the first term on the R.H.S, we have that
\begin{align*}
\sum_{s=1}^{S+1} w_s\Delta_s &= \sum_{s=1}^{S+1}w_s (\phi(\x_{s-1}) - \phi(\x_s)) = \sum_{s=1}^{S+1} (w_{s-1}\phi(\x_{s-1}) - w_s\phi(\x_s)) + \sum_{s=1}^{S+1}(w_s - w_{s-1})\phi(\x_{s-1})\\
&\leq w_0 \phi(\x_0) - w_{S+1}\phi(\x_{S+1}) +\sum_{s=1}^{S+1}(w_s - w_{s-1})\phi(\x_{s-1})\\
& =\sum_{s=1}^{S+1}(w_s - w_{s-1})(\phi(\x_{s-1}) - \phi(\x_{S+1}))\leq \Delta_\phi\sum_{s=1}^{S+1}(w_s - w_{s-1}) = \Delta_\phi w_{S+1}
\end{align*}
Then, 
\begin{align*}
   \E[\|\nabla\phi_{\gamma}(\x_\tau)\|^2] \leq \frac{16\mu\Delta_\phi w_{S+1}}{\sum_{s=1}^{S+1}w_s} + \frac{24\mu\sum_{s=1}^{S+1}w_s\eta_{s}\hat G^2}{\sum_{s=1}^{S+1}w_s}
\end{align*}
The standard calculus tells that
\begin{align*}
\sum_{s=1}^Ss^{\alpha}&\geq \int_0^{S}x^\alpha d x= \frac{1}{\alpha+1}S^{\alpha+1}\\
\sum_{s=1}^Ss^{\alpha-1}&\leq SS^{\alpha-1} = S^\alpha, \forall\alpha\geq 1, \quad \sum_{s=1}^Ss^{\alpha-1}\leq \int_{0}^S x^{\alpha-1}d x = \frac{ S^\alpha}{\alpha}, \forall 0<\alpha<1
\end{align*}
Combining these facts, we have that
\begin{align*}
 \E[\|\nabla\phi_{\gamma}(\x_{\tau})\|^2]\leq \left\{\begin{array}{cc}\frac{16\mu\Delta_\phi(\alpha+1)}{S+1} +  \frac{24\mu \hat G^2(\alpha+1)}{S+1} & \alpha\geq 1\\\\ \frac{16\mu\Delta_\phi(\alpha+1)}{S+1} +  \frac{24\mu \hat G^2(\alpha+1)}{(S+1)\alpha} & 0<\alpha< 1
 \end{array}\right.
\end{align*}
In order to have $\E[\|\nabla\phi_{\gamma}(\x_{\tau})\|^2]\leq \epsilon^2$, we can set $S=O(1/\epsilon^2)$. The total number of iterations is
\begin{align*}
\sum_{s=1}^{S} T_s \leq \sum_{s=1}^{S}12\gamma s\leq 6\gamma S(S+1) = O(1/\epsilon^4)
\end{align*}
\end{proof}

\section{Proof of Theorem~\ref{thm:nsum}}
\begin{proof}
    According to the definition of $\z_s$ in (\ref{eqn:z}) and 
Lemma~\ref{lem:sum}, we have that
\begin{align*}
    &\hspace{6pt}\E_{s} \bigg[\phi(\x_{s}) + \frac{1}{2\gamma}\|\x_{s}-\x_{s-1}\|^2 \bigg]\\
     &\leq f_{s}(\z_{s}) + \underbrace{\frac{\beta (f_s(\x_{s-1})-f_s(\z_s))}{(1-\beta)(T_s+1)} + \frac{(1-\beta)\|\x_{s-1}-\z_s\|^2}{2\eta_s{(T_s+1)}} + \frac{2\eta_s {G}^2(2\rho\beta+1) }{1-\beta}  +\frac{1}{24\gamma}\|\x_{s-1} - \z_s\|^2}\limits_{\mathcal{E}_s}\\
    &\leq \phi(\x_{s-1})+\mathcal{E}_{s}.
    \end{align*}
Similar to the proof of Theorem~\ref{thm:sgd}, we have
\begin{align}\label{eqn:bound1}
    \frac{(1-\alpha_{s})}{2\gamma} \|\x_{s-1}-\z_{s}\|^2 \leq& 
    \E_s[\phi(\x_{s-1}) -\phi(\x_{s})] + \frac{\alpha_{s}^{-1}-\gamma\mu}{(1-\gamma\mu)}\mathcal{E}_{s}
\end{align}
Rearranging above inequality, we have that 
\begin{align*}
    &\bigg((1-\alpha_{s})\gamma - \frac{ \gamma^2(\alpha_{s}^{-1} - \mu \gamma)(1-\beta)}{ (1-\mu\gamma)\eta_{s} ({T_{s}+1})} - \frac{\alpha_{s}^{-1}-\gamma\mu}{(1-\gamma\mu)}\frac{\gamma}{24}\bigg) \|\nabla \phi_{\gamma}(\x_{s-1})\|^2\\
    \leq& 2\E_s[\Delta_{s}] + \frac{2(\alpha_{s}^{-1} - \mu\gamma )}{(1-\mu\gamma)}\bigg[\frac{\beta(f_s(\x_{s-1})-f_s(\z_{s}))}{(1-\beta)(T_{s}+1)} + \frac{2\eta_s\hat{G}^2(2\rho\beta+1) }{1-\beta}\bigg]
\end{align*}
The definition of $f_s$ gives that
\begin{align*}
f_s(\x_{s-1})-f_s(\z_{s}) =\phi(\x_{s-1})-\phi(\z_{s})-\frac{1}{2\gamma}\|\z_{s}-\x_{s-1}\|^2
\end{align*}
On the other hand, the $\mu$-weakly convexity of $\phi$ gives that
\begin{align*}
\phi(\z_s)\geq \phi(\x_{s-1}) + \langle \g(\x_{s-1}), \z_s-\x_{s-1}\rangle -\frac{\mu}{2} \|\z_s-\x_{s-1}\|^2,
\end{align*}
where $\g(\x_{s-1})\in \partial_F\phi(\x_{s-1})$. Combing these two inequalities we have that
\begin{align*}
    f_s(\x_{s-1})-f_s(\z_{s}) \leq& \langle \g(\x_{s-1}), \x_{s-1} - \z_s\rangle -\frac{\mu}{2} \|\z_s-\x_{s-1}\|^2\\
    \leq & \frac{G^2}{2\mu}+\frac{\mu-\mu}{2}\|\z_s-\x_{s-1}\|^2 = \frac{G^2}{2\mu}\end{align*} 
where the second inequality follows from Jensen's inequality for $\|\cdot\|$ and Young's inequality.
Combining above inequalities and multiplying both side by $w_s$, we have that
\begin{align}\label{eqn:bound2}
    &\hspace{-6pt}w_s\bigg((1-\alpha_{s})\gamma - \frac{ \gamma^2(\alpha_{s}^{-1} - \mu \gamma)(1-\beta)}{ (1-\mu\gamma)\eta_{s} ({T_{s}+1})}- \frac{\alpha_{s}^{-1}-\gamma\mu}{(1-\gamma\mu)}\frac{\gamma}{24}\bigg) \|\nabla \phi_{\gamma}(\x_{s-1})\|^2 \nonumber\\
    \leq& 2w_s\E_s[\Delta_{s}] + \frac{2w_s(\alpha_{s}^{-1} - \mu\gamma )}{(1-\mu\gamma )}\bigg[\frac{\beta G^2}{2\mu(1-\beta)(T_{s}+1)} + \frac{2\eta_s {G}^2(2\rho\beta+1)}{1-\beta}\bigg]
\end{align}
By setting $\alpha_{s} = 1/2$, $\eta_s(T_s+1)\geq 24(1-\beta)\gamma$, we have that
\begin{align*}
\frac{ w_s\gamma}{4}\|\nabla \phi_{\gamma}(\x_{s-1})\|^2\leq 2w_s\E_s[\Delta_s]+\frac{ w_s\eta_s\beta G^2}{8(1-\beta)^2} + \frac{12w_s\eta_s{G}^2(2\rho\beta+1) }{1-\beta}
\end{align*}
Summing over $s=1,\ldots,S+1$ and rearranging, we have
\begin{align*}
\sum_{s=1}^{S+1} w_s\|\nabla\phi_{\gamma}(\x_{s-1})\|^2 = \E\bigg[\sum_{s=1}^{S+1} \frac{8}{\gamma} w_s\Delta_s + \frac{w_s\eta_s G^2(\beta+ 96 (2\rho\beta+1)(1-\beta)) }{ 2\gamma(1-\beta)^2}\bigg]
\end{align*}
Following similar analysis as in the proof of Theorem~\ref{thm:sgd}, we can finish the proof. 
\end{proof}

\section{Proof of Theorem~\ref{thm:nadagrad}}
\begin{proof}
Applying Lemma~\ref{lem:adagrad} with $T_s\geq M_s\max\{\frac{\hat G + \max_i\|g^s_{1:T_s,i}\|}{2},  \sum_{i=1}^d\|g^s_{1:T_s,i}\|\}$ $M_s>0$, and the fact that  $\phi(\x_{s-1})\geq \phi(\z_{s})+\frac{1}{2\gamma}\|\x_{s-1}-\z_{s}\|^2$ in $s$th stage, we have that
\begin{align*}
    \E_{s} \bigg[\phi(\x_{s}) + \frac{1}{2\gamma_{s}}\|\x_{s}-\x_{s-1}\|^2 \bigg] &\leq f_{s}(\z_{s}) + \underbrace{\frac{1}{M_s\eta_s}\|\x_{s-1}-\z_s\|^2 +\frac{\eta_s}{M_s}}\limits_{\mathcal{E}_s}\\
    &\leq \phi(\x_s)+\mathcal{E}_{s}
    \end{align*}
According to (\ref{eqn:bound1}), we have that 
\begin{align*}
    \frac{(1-\alpha_{s})}{2\gamma} \E_s[\|\x_{s-1}-\z_{s}\|^2] \leq&  \phi(\x_{s-1}) -\phi(\x_{s})+\frac{(\alpha_{s}^{-1}-1)}{2\gamma}\|\x_{s}-\z_{s}\|^2+\mathcal{E}_{s} \nonumber\\
    \leq & \phi(\x_{s-1}) -\phi(\x_{s}) + \frac{\alpha_{s}^{-1}-\gamma\mu}{(1-\gamma\mu)} \bigg(\frac{1}{M_s\eta_s}\|\x_{s-1}-\z_s\|^2 +\frac{\eta_s}{M_s}\bigg)
\end{align*}
Rearranging above inequality then multiplying both side by $w_s$, we have that 
\begin{align*}
w_s\bigg((1-\alpha_{s})\gamma - &\frac{ 2\gamma^2(\alpha_{s}^{-1} - \mu \gamma)}{ (1-\mu\gamma)M_s\eta_{s} } \bigg) \|\nabla \phi_{\gamma}(\x_{s-1})\|^2\\
    \leq& 2w_s\E_s[\Delta_{s}] + \frac{2w_s\eta_s(\alpha_{s}^{-1} - \mu\gamma )}{M_s(1-\mu\gamma)}
\end{align*}
By using $M_s\eta_s \geq 24\gamma $ and summing over $s=1,\ldots, S+1$, we have that
\begin{align*}
    \sum_{s=1}^{S+1} w_s\|\nabla\phi_{\gamma}(\x_{s-1})\|^2 \leq \E\bigg[ \sum_{s=1}^{S+1} \frac{8w_s\Delta_s}{\gamma} + \frac{w_s\eta_s^2}{\gamma^2}\bigg]
\end{align*}
By the definition of $p_s$ in the theorem, taking expectation of $\|\nabla\phi_{\gamma}(\x_\tau)\|^2$ w.r.t. $\tau\in\{0,\ldots, S\}$ we have that
\begin{align*}
    \E[\|\nabla\phi_{\gamma}(\x_\tau)\|^2] =& \E\bigg[\frac{8}{\gamma}\sum_{s=1}^{S+1} \frac{w_s\Delta_s}{\sum_{i=1}^{S+1}w_i}\bigg] + \frac{c^2}{\gamma^2}\sum_{s=1}^{S+1}\frac{s^{\alpha-1}}{\sum_{i=1}^{S+1}w_i}\\
    \leq& \frac{8\Delta_\phi(\alpha+1)}{\gamma(S+1)} + \frac{c^2(\alpha+1)}{\gamma^2(S+1)\alpha^{\mathbb I(\alpha<1)}}
\end{align*}        



\end{proof}

\section{Proof of Lemma~\ref{lem:sum}}
\begin{proof}
Following the analysis in~\cite{yangnonconvexmo}, we directly have the following inequality, 
\begin{align*}
&\E[\|\x_{k+1} + \p_{k+1} - \x_*\|^2] =\\
& =  \E[\|\x_k + \p_k - \x_*\|^2] - \frac{2\eta}{1-\beta}\E[(\x_k  - \x_*)^{\top}\partial f(\x_k)]- \frac{2\eta\beta}{(1-\beta)^2}\E[(\x_k - \x_{k-1})^{\top}\partial f(\x_k)]\nonumber\\
&- \frac{2\rho\eta^2\beta}{(1-\beta)^2}\E[\g_{k-1}^{\top} \partial f(\x_k)]+ \left(\frac{\eta}{1-\beta}\right)^2\E[|\g_k\|^2]\nonumber 
\end{align*}
We also note that 
\begin{equation*}
\begin{aligned}
&f(\x_k ) - f(\x_*)\leq (\x_k - \x_*)^{\top}\partial f(\x_k) -  \frac{\lambda}{2}\|\x_k - \x_*\|^2\\
&f(\x_k) - f(\x_{k-1})\leq (\x_k - \x_{k-1})^{\top}\partial f(\x_k) -  \frac{\lambda}{2}\|\x_k - \x_{k-1}\|^2\\
&-\E[\g_{k-1}^{\top}\partial f(\x_k)]\leq \frac{\E[\|\g_{k-1}\|^2 + \|\partial f(\x_k)\|^2]}{2}\leq \frac{1}{\gamma^2}\|\x_{k-1} - \x_0\|^2 +  \frac{1}{\gamma^2}\|\x_{k} - \x_0\|^2 + 2G^2\\
&\E_k[\|\g_k\|^2]\leq \frac{2}{\gamma^2}\|\x_k - \x_0\|^2 + 2G^2
\end{aligned}
\end{equation*}
where the first two inequalities are due to the strong convexity of $f(\cdot)$ and the last three inequalities are due to the boundness assumption.  Thus
\begin{align*}
\E[\|\x_{k+1} + \p_{k+1} - \x\|^2&] \leq   \E[\|\x_k + \p_k - \x\|^2] - \frac{2\eta}{1-\beta}\E[(f(\x_k) - f(\x))]\\
 &- \frac{2\eta\beta}{(1-\beta)^2}\E[(f(\x_k) - f(\x_{k-1}))]  + \left(\frac{\eta}{1-\beta}\right)^2(2\rho\beta+1) 4G^2\\
 & -   \frac{\lambda\eta}{1-\beta}\|\x_k - \x_*\|^2 -    \frac{\lambda\eta\beta}{(1-\beta)^2}\|\x_k - \x_{k-1}\|^2\\
 &  + \frac{2\rho\beta}{(1-\beta)^2} \frac{\eta^2}{\gamma^2}\|\x_{k-1} - \x_0\|^2 + \frac{2\rho\beta + 2}{(1-\beta)^2}\frac{\eta^2}{\gamma^2}\|\x_k - \x_0\|^2
\end{align*}
By summarizing the above inequality over $k=0, \ldots, T$, we have
\begin{align*}
\frac{2\eta}{1-\beta}&\E\bigg[\sum_{k=0}^T(f(\x_k) - f(\x_*))\bigg]\leq \E[\|\x_0 - \x_*\|^2] + \frac{2\eta\beta}{(1-\beta)^2}\E[f(\x_0) - f(\x_*)]\\
&+ \left(\frac{\eta}{1-\beta}\right)^2(2\rho\beta+1) 4G^2(T+1)\\
& - \frac{\eta \lambda}{1-\beta}\sum_{k=0}^T\|\x_k - \x_*\|^2  + \frac{4\rho\beta}{(1-\beta)^2} \frac{\eta^2}{\gamma^2}\sum_{k=0}^{T}\|\x_{k-1} - \x_*\|^2 + \frac{4\rho\beta + 4}{(1-\beta)^2}\frac{\eta^2}{\gamma^2}\sum_{k=0}^T\|\x_k - \x_*\|^2\\
& +\frac{4\rho\beta+4}{(1-\beta)^2}\frac{\eta^2}{\gamma^2}(T+1)\|\x_0 - \x_*\|^2
\end{align*}
When $\eta\leq (1-\beta)\gamma^2\lambda/(8\rho\beta+4)$, we have
\begin{align*}
&\E\bigg[(f(\xh_T) - f(\x_*))\bigg]\leq \frac{(1-\beta)\|\x_0 - \x_*\|^2}{2\eta (T+1)} + \frac{\beta}{(1-\beta)}\frac{f(\x_0) - f(\x_*)}{T+1}  + \frac{\eta}{1-\beta}(2\rho\beta+1) 2G^2\\
& +\frac{4\rho\beta+4}{(1-\beta)}\frac{\eta}{\gamma^2}\|\x_0 - \x_*\|^2
\end{align*}
\end{proof}

    \section{Proof of Lemma~\ref{lem:adagrad}}
    The proof is almost a duplicate of the proof of Proposition 1 in~\cite{SadaGrad18}. For completeness, we present a proof here. 
    \begin{proof}
        Let $\psi_0(\x) = 0$ and $\|\x\|_{H}=\sqrt{\x^{\top}H\x}$. First, we can see that $\psi_{t+1}(\x)\geq \psi_t(\x)$ for any $t\geq 0$. Define $\zeta_t  = \sum_{\tau=1}^t\g_t$ and $\Delta_\tau  = (\partial F(\x_t) - \g_t)^{\top}(\x_t - \x)$. Let  $\psi_t^*$ be defined by
        \begin{align*}
        \psi^*_t(g) = \sup_{\x\in\Omega}g^{\top}\x - \frac{1}{\eta}\psi_t(\x)
        \end{align*}
        Taking the summation of objective gap in all iterations, we have
        \begin{align*}
        &\sum_{t=1}^T(f(\x_t) - f(\x)) \leq \sum_{t=1}^T\partial f(\x_t)^{\top}(\x_t - \x)  =  \sum_{t=1}^T\g_t^{\top}(\x_t - \x) + \sum_{t=1}^T\Delta_t\\
        & = \sum_{t=1}^T\g_t^{\top}\x_t - \sum_{t=1}^T\g_t^{\top}\x  - \frac{1}{\eta}\psi_T(\x) + \frac{1}{\eta}\psi_T(\x)+  \sum_{t=1}^T\Delta_t\\
        &\leq \frac{1}{\eta}\psi_T(\x) +  \sum_{t=1}^T\g_t^{\top}\x_t   + \sum_{t=1}^T\Delta_t + \sup_{\x\in\Omega}\left\{ - \sum_{t=1}^T\g_t^{\top}\x - \frac{1}{\eta}\psi_T(\x)\right\}\\
         & = \frac{1}{\eta}\psi_T(\x) + \sum_{t=1}^T\g_t^{\top}\x_t + \psi_T^*(-\zeta_T)  + \sum_{t=1}^T\Delta_t
        \end{align*}
        Note that
        \begin{align*}
        &\psi_T^*( - \zeta_T)  =  - \sum_{t=1}^T\g_t^{\top}\x_{T+1} - \frac{1}{\eta}\psi_T(\x_{T+1})\leq -\sum_{t=1}^T\g_t^{\top}\x_{T+1} - \frac{1}{\eta}\psi_{T-1}(\x_{T+1})\\
        &\leq \sup_{\x\in\Omega}- \zeta_T^{\top}\x - \frac{1}{\eta}\psi_{T-1}(\x) = \psi_{T-1}^*(-\zeta_T)\\
        &\leq \psi_{T-1}^*(-\zeta_{T-1}) - \g_T^{\top}\nabla\psi_{T-1}^*(-\zeta_{T-1}) + \frac{\eta}{2}\|\g_T\|_{\psi^*_{T-1}}^2
        \end{align*}
        where the last inequality uses the fact that $\psi_t(\x)$ is 1-strongly convex w.r.t $\|\cdot\|_{\psi_t} = \|\cdot\|_{H_t}$ and consequentially $\psi_t^*(\x)$ is $\eta$-smooth wr.t. $\|\cdot\|_{\psi^*_t} = \|\cdot\|_{H_t^{-1}}$.
        Thus, we have
        \begin{align*}
         &\sum_{t=1}^T\g_t^{\top}\x_t + \psi_T^*(-\zeta_T)\leq  \sum_{t=1}^T\g_t^{\top}\x_t + \psi_{T-1}^*(-\zeta_{T-1}) - \g_T^{\top}\nabla\psi_{T-1}^*(-\zeta_{T-1})+ \frac{\eta}{2}\|\g_T\|_{\psi^*_{T-1}}^2\\
         & = \sum_{t=1}^{T-1}\g_t^{\top}\x_t + \psi^*_{T-1}(-\zeta_{T-1})  + \frac{\eta}{2}\|\g_T\|_{\psi^*_{T-1}}^2
        \end{align*}
        By repeating this process, we have
        \begin{align*}
        &\sum_{t=1}^T\g_t^{\top}\x_t + \psi_T^*(-\zeta_T)\leq \psi_0^*(-\zeta_0) +\frac{\eta}{2} \sum_{t=1}^T\|\g_t\|^2_{\psi^*_{t-1}}=  \frac{\eta}{2}\sum_{t=1}^T\|\g_t\|^2_{\psi^*_{t-1}}
        \end{align*}
        Then
        \begin{align}\label{prop1:ineq1}
        \sum_{t=1}^T(f(\x_t) - f(\x))\leq& \frac{1}{\eta}\psi_T(\x) +\frac{\eta}{2}\sum_{t=1}^T\|\g_t\|^2_{\psi^*_{t-1}}+\sum_{t=1}^T\Delta_t
        \end{align}
        Following the analysis in~\cite{duchi2011adaptive}, we have
        \begin{align*}
        \sum_{t=1}^T\|\g_t\|_{\psi^*_{t-1}}^2\leq 2\sum_{i=1}^d\|\g_{1:T, i}\|_2
        \end{align*}
        Thus
        \begin{align*}
        &\sum_{t=1}^T(f(\x_t) - f(\x))\leq\frac{G\|\x - \x_1\|_2^2}{2\eta} + \frac{(\x - \x_1)^{\top}\diag(s_T)(\x-\x_1)}{2\eta}+ \eta\sum_{i=1}^d\|\g_{1:T, i}\|_2 + \sum_{t=1}^T\Delta_t\\
        &\leq \frac{G  +\max_i \|\g_{1:T,i}\|_2}{2\eta}\|\x - \x_1\|_2^2  + \eta\sum_{i=1}^d\|\g_{1:T,i}\|_2 +\sum_{t=1}^T\Delta_t
        \end{align*}
        Now by the value  of $T\geq M\max\{\frac{G+ \max_i\|g_{1:T,i}\|}{2},  \sum_{i=1}^d\|g_{1:T,i}\|\}$, we have
        \begin{align*}
        & \frac{(G + \max_i \|g_{1:T, i}\|_2)}{2\eta T}\leq \frac{1}{\eta M}\\
        & \frac{\eta\sum_{i=1}^d\|g_{1:T,i}\|_2}{T}   \leq \frac{\eta}{M}
        \end{align*}
        Dividing by $T$ on both sides and setting $\x=\x_*$, following the inequality (3) and the convexity of $f(\x)$ we have
        \begin{align*}
        f(\xh)-f_*\leq \frac{1}{M\eta}\|\x_{0}-\x_*\|^2 +\frac{\eta}{M} + \frac{1}{T}\sum_{t=1}^T\Delta_t
        \end{align*}
        Let $\{\mathcal F_t\}$ be the filtration associated with Algorithm 1 in the paper. Noticing that $T$ is a random variable with respect to $\{\mathcal F_t\}$, we cannot get rid of the last term directly. Define the Sequence $\{X_t\}_{t\in\N_+}$ as
        \begin{align}\label{eqn:sup}
        X_t=\frac{1}{t}\sum_{i=1}^t\Delta_i=\frac{1}{t} \sum_{i=1}^t \langle\g_i-\E[\g_i], \x_i-\x_*\rangle
        \end{align}
        where $\E[\g_i]\in\partial f(\x_i)$. 
        Since $\E\left[\g_{t+1}- \E[\g_{t+1}]\right]=0$ and $\x_{t+1}=\arg\min\limits_{\x\in\Omega}\eta \x^{\top}\left(\frac{1}{t}\sum_{\tau=1}^t\g_\tau\right) + \frac{1}{t}\psi_t(\x)$, which is measurable with respect to $\g_1,\ldots,\g_t$ and $\x_1,\ldots,\x_t$, it is easy to see $\{\Delta_t\}_{t\in N}$ is a martingale difference sequence with respect to $\{\F_t\}$, e.g. $E[\Delta_t|\F_{t-1}]=0$. On the other hand, since $\|\g_t\|_2$ is upper bounded by $G$, following the statement of $T$ in the theorem, $T\leq N= M^2\max\{ \frac{(G+1)^2}{4},d^2G^2 \}<\infty$ always holds. Then following Lemma~1 in \cite{SadaGrad18} we have that $\E[X_T]=0$.
        
        Now taking the expectation we have that
        \begin{align*}
        &\E[f(\xh)-f_*]\leq\frac{1}{M\eta}\|\x_{0}-\x_*\|^2 +\frac{\eta}{M}
        \end{align*}
        Then we finish the proof.
    \end{proof}

\bibliographystyle{abbrv}

\bibliography{main}
\end{document}